%% file: 0_main.tex
\documentclass[11pt,twoside]{amsart}
\input{includeNice3}
\usepackage{mabliautoref}
\usepackage{amssymb}
\usepackage{amscd}
\usepackage[abbrev,alphabetic]{amsrefs}
\usepackage{hyperref}
\usepackage[a4paper,margin=1in]{geometry}  
\usepackage{CJKutf8}
\usepackage{soul}
\usepackage{multicol}

\usepackage{xypic}
\usepackage{tikz-cd}
\usepackage{caption}

\usepackage{url}
\RequirePackage{booktabs, multirow}
\RequirePackage{pgf}

\title[Equisingular lifting of $F$-split slc CY surfaces]
{Equisingular lifting  of semi-log canonical $F$-split $K$-trivial surfaces} 
\author[F.~Bernasconi and Q.~Posva]{Fabio Bernasconi and Quentin Posva} 
\subjclass[2020]{Primary: 14G17, 14J32.}
\keywords{Calabi--Yau surfaces, global $F$-splitting, lifting, positive and mixed characteristic.}

\address{Dipartimento di Matematica “Guido
Castelnuovo”, SAPIENZA Università di Roma, Piazzale Aldo Moro 5, I-00185,
Roma, Italy} 
\email{fabio.bernasconi@uniroma1.it}

\address{Institut de Mathématiques, Université de Neuchâtel, Rue Émile-Argand 11, 2000 Neuchâtel, Switzerland} 
\email{quentin.posva@unine.ch}

\DeclareMathOperator{\Diff}{Diff}

\DeclareMathOperator{\Sing}{Sing}
\DeclareMathOperator{\Supp}{Supp}

\DeclareMathOperator{\id}{id}
\DeclareMathOperator{\Spec}{Spec}
\DeclareMathOperator{\PGL}{PGL}
\DeclareMathOperator{\GL}{GL}
\DeclareMathOperator{\Pic}{Pic}
\DeclareMathOperator{\Aut}{Aut}
\DeclareMathOperator{\Frac}{Frac}
\DeclareMathOperator{\Ext}{Ext}

\DeclareMathOperator{\NS}{NS}
\DeclareMathOperator{\Hom}{Hom}
\DeclareMathOperator{\Tr}{Tr}
\DeclareMathOperator{\can}{can}
\DeclareMathOperator{\Exc}{Exc}
\DeclareMathOperator{\pr}{pr}
\DeclareMathOperator{\sn}{sn}
\DeclareMathOperator{\Ram}{Ram}
\DeclareMathOperator{\Ex}{Ex}
\DeclareMathOperator{\Stab}{Stab}
\DeclareMathOperator{\rk}{rk}

\DeclareMathOperator{\Isom}{Isom}
\DeclareMathOperator{\Fix}{Fix}

\DeclareMathOperator{\Eff}{Eff}
\DeclareMathOperator{\Bl}{Bl}
\DeclareMathOperator{\NE}{NE}
\DeclareMathOperator{\res}{res}
\DeclareMathOperator{\im}{im}

\DeclareMathOperator{\codim}{codim}
\DeclareMathOperator{\Tot}{Tot}

\newcommand{\sL}{\mathcal{L}}
\newcommand{\sA}{\mathcal{A}}
\newcommand{\sHom}{\mathcal{H}om}
\newcommand{\bP}{\mathbb{P}}
\newcommand{\bA}{\mathbb{A}}
\newcommand{\sT}{\mathcal{T}}

\newcommand{\sO}{\mathcal{O}}

\newcommand{\fm}{\mathfrak{m}}

\newcommand{\sS}{\mathcal{S}}
\newcommand{\sD}{\mathcal{D}}

\newcommand{\sE}{\mathcal{E}}

\newcommand{\bQ}{\mathbb{Q}}

\newcommand{\sF}{\mathcal{F}}
\newcommand{\sC}{\mathcal{C}}
\newcommand{\sX}{\mathcal{X}}
\newcommand{\sU}{\mathcal{U}}

\newcommand{\bR}{\mathbb{R}}
\newcommand{\sY}{\mathcal{Y}}

\begin{document}
	
\maketitle

\begin{abstract}
    We show that projective globally $F$-split semi-log canonical $K$-trivial surfaces over an algebraically closed field of characteristic $p>0$ admit equisingular liftings over the ring of Witt vectors.
\end{abstract}

\tableofcontents

\input{1_introduction}

\input{2_prelims}

\input{2.5_Canonical_liftings}

\input{3_Log_liftings}

\input{4_Canonical_liftings}

\input{6_case_of_4A1}

\bibliographystyle{amsalpha}
\bibliography{refs}

\end{document}

%% file: 1_introduction.tex
\section{Introduction}

Given a projective variety $X$ over an algebraically closed field $k$ of positive characteristic, it is a natural question to ask whether $X$ can be lifted to a projective variety defined over a field of characteristic 0.
This is not always the case as Serre constructed an example of a non-liftable variety in \cite{Ser61}.
Since then many efforts have been dedicated to either construct new examples (see e.g.\ \cite{Hirokado_Non_liftable_CY_3fold, Schroeer_CY_3fold_with_obstructed_def, Cynk_van_Straten_Small_resolutions_and_non_liftable_CY, Cascini_Tanaka_KVV_Failure, AZ-nonlift, vDB21, RS22}), identifying necessary conditions for the existence of a lifting (as in \cite{DI87, Lan15, Pet23}) or determining special classes of varieties that do lift, especially in the Fano and $K$-trivial cases (see e.g.\ \cite{Del81, Nyg83, MS87, Sch03, 
Liedtke_Arithmetic_moduli_and_liftings_for_Enriques,
AZ19, LT19, KTTWYY22}). The present paper lies in this last category.

 It is expected that lifts exist for the varieties $X$ whose Frobenius morphism $F$ satisfy some arithmetic condition. For example, it is common to study varieties with trivial canonical bundle which are weakly ordinary, i.e.\ the action of the Cartier operator is surjective on top forms.
 A generalisation of weak ordinarity, coming from commutative algebra, is the following: we say that $X$ is \emph{globally $F$-split} if the natural $\mathcal{O}_X$-module homomorphism $\sO_X\to F_*\sO_X$ is split.
 In \cite{Zda18}, it is shown that globally $F$-split variety admit natural lifts to $W_2(k)$.
 Motivated by Zdanowicz's result, in \cite[Directions 1.7]{AZ19} the authors ask whether globally $F$-split varieties always lift over the ring of Witt vectors.  
 A striking recent result seems to provide support to this speculation: if $X$ is a smooth variety with trivial canonical sheaf that is Bloch--Kato ordinary (which is stronger than just being globally $F$-split), then $X$ lifts to $W(k)$ \cite{brantner2025liftcy}.
 Another evidence supporting this conjecture  is given by a recent result of Petrov which states that the de Rham complex $F_*\Omega_{X/k}^{\bullet}$ of a smooth globally $F$-split variety decomposes \cite{Petrov_Decomposition_dR_complex_for_qFsplit}.

In dimension 2, the lifting of normal globally $F$-split surfaces was settled affirmatively in \cite{BBKW24}. In fact, the authors proved a much stronger result: given a globally $F$-split normal surface $X$ over $k$, it is possible to find an \emph{equisingular} lifting of $X$ over the Witt vectors. 

In this article, we generalise their result to the case of non-normal semi-log canonical (slc) $K$-trivial surfaces.
Our motivation comes from a future investigation of the moduli theory of $K$-trivial surfaces in mixed characteristic, as several foundational questions have been settled in recent years \cite{7authors, Posva_Gluing_stable_families_surfaces_mix_char, ABP24}. 
From the point of view of the moduli theory of higher dimensional algebraic varieties \cite{k-modulibook}, non-normal slc $K$-trivial varieties appear naturally as degenerations of smooth $K$-trivial surfaces obtained by running suitable Minimal Model Program algorithms, and thus lie in the boundary of a possible moduli space of $K$-trivial surfaces. 
Our generalization of the main result of \cite{BBKW24} would then show that the open part of the boundary describing globally $F$-split $K$-trivial surfaces in characteristic $p$ belongs to the specialisation of the boundary of the moduli space in characteristic 0. 

Our main result reads as follows:

\begin{theorem}[\autoref{thm:main_with_technical_statement}]\label{main_thm}
   Let $k$ be an algebraically closed field of characteristic $p>0$.
   Let $X$ be a projective semi-log canonical globally $F$-split projective surface with $K_X \equiv 0$.
   Then there exists a projective strong semi-log canonical lifting $\mathcal{X}$ of $X$ to $W(k)$.
\end{theorem}

The technical meaning of $\sX$ being a strong slc lifting of $X$ is given in \autoref{section:def_log_lifts}: it should be thought as an equisingular deformation of $X$ in mixed characteristic. The core idea is the following: there should exist a Cartesian diagram
        $$\begin{tikzcd}
        (X^\nu, D) \arrow[d, "\nu_X"]\arrow[r, hook] &
        (\sX^\nu,\sD) \arrow[d, "\nu_\sX"] \\
        X\arrow[r, hook] & \sX
        \end{tikzcd}$$
where $(X^\nu,D)$ is the normalization couple of $X$ (here $D$ is a reduced divisor: see \autoref{section:gluing_theory} for a review of normalization of slc varieties), $\sX$ is demi-normal with special fiber $X$, and the normalization $(\sX^\nu,\sD)$ of $\sX$ is a strong log lifting of $(X^\nu,D)$ in the sense of \cite{BBKW24}.

\subsection*{Sketch of the proof}
Let us give a quick overview of the proof of \autoref{main_thm}. For now, assume that $p>2$ and let $X$ be as in the statement: we may assume that it is not normal, since the normal case is treated in \cite{BBKW24}. As $X$ is slc, the normalization pair $(X^\nu,D)$ comes with a non-trivial log involution $\tau$ of $(D^\nu,\Diff(0))$. It turns out that finding an equisingular lift of $X$ is equivalent to finding a strong log lifting $(\sX^\nu,\sD)$ of $(X^\nu, D)$ with $K_{\sX^\nu}+\sD$ being  $\mathbb{Q}$-Cartier, together with an involution of $(\sD^\nu,\Diff(0))$ lifting $\tau$ (see \autoref{lem: suff-cond-log-lift-nonnormal}).

As \cite{BBKW24} provides us with strong log liftings of $(X^\nu,D)$, our task is to lift $\tau$. The difficulty of doing so depends on the components of the log CY curve $(D^\nu,\Diff(0))$. The ones that are rational and whose boundary contains at most three points pose no difficulty (\autoref{cor:lift_involutions_on_P1}): in fact, most of this paper is dealing with the remaining cases. They are of two kinds: the components of $D^\nu$ of genus one, and the rational components whose boundaries contain four points each. An example of the later is given in \autoref{example:4A_1}.

Let $D_1$ be a component of $D^\nu$ which is of either of these two types. By adjunction, the pair $(D_1,\Diff(0))$ is globally $F$-split. Our strategy is to show that we can find a strong log lifting $(\sX^\nu,\sD)$ where the component $\sD_1$ of $\sD^\nu$ lifting $D_1$ is the \emph{canonical lift} of $D_1$ over the Witt vectors. If $D_1$ has genus one, then it is ordinary, and we use ``canonical lift" in the sense of \cite[Appendix]{MS87}. If $D_1$ has genus zero, then the meaning of ``canonical lift" is given in \autoref{section:can_lift_P1_4pts}: in that case, the key observation is that there is an ordinary genus one curve $E$ and a finite morphism $E\to D_1$ branched over the support of $\Diff(0)$---and the canonical lift of $D_1$ is induced by the one of $E$. In both cases, general properties of the canonical lifts ensure that $\tau$ can be uniquely lifted along them.

The construction of strong log liftings $(\sX^\nu,\sD)$ that induce the canonical liftings of the special components of $D^\nu$ is established by a case-by-case analysis.

In case $D_1$ is a component of genus one, we observe that either $K_{X^\nu}+D\sim 0$; or the crepant birational equivalence class of $(X^\nu,D)$ contains a surface pair of very particular type, and in that case we call $(X^\nu,D)$ \emph{of Enriques-type} (\autoref{prop:genus_1_forces_Z_linear_CY_cond}). In both cases we classify the smooth minimal models: when $K_{X^\nu}+D\sim 0$ this is done in \autoref{lem: class_logsnc_CY} (in which we fix a mistake of \cite[Proposition 5.9]{BBKW24} where a case in characteristic $2$ was overlooked); when $(X^\nu,D)$ is of Enriques--type, the classification is done in \autoref{prop: classification-enriques-type}. Using these classifications we construct adequate strong log liftings in \autoref{prop: gen_canonical_lift_BBKW}.


The case where $D_1$ is a rational curve with four boundary points is more delicate, and occupies \autoref{section:4A_1}. The crucial point is: if $Z\to X^\nu$ is the cyclic cover given by the $\bQ$-Cartier divisor $K_{X^\nu}+D$, then the induced branched cover of $D_1$ is the morphism $E\to D_1$ mentioned above in the definition of the canonical lift of $D_1$. It turns out that $K_{X^\nu}+D$ is $2$-Cartier, so we reduce to lift the pair $(Z,E)$ together with its $\mu_2$-action over $W(k)$, in such a way that the lift of the ordinary curve $E$ is the canonical one. To do so, we classify the $\mu_2$-equivariant minimal models of $(Z,E)$ and lift them one-by-one over $W(k)$ in an appropriate way. The price of having to work $\mu_2$-equivariantly is balanced out by a reduction to very explicit geometrical arguments and cohomological vanishings on conic bundles and del Pezzo surfaces.

The case of $p=2$ is treated separately in \autoref{section:p=2}. In this case it turns out that $K_X$ is always Cartier (\autoref{lemma:GFS_when_p=2}): so the method of \autoref{thm: main_KX-trivial} applies, once some tools from \cite{Posva_Gluing_stable_families_surfaces_mix_char} used in \autoref{lem: suff-cond-log-lift-nonnormal} are suitably adapted. 
We note that the liftability over $W(k)$ of semi-smooth Enriques surfaces in characteristic $2$ has been thoroughly studied by Schröer \cite{Sch03}.

\subsection*{Organisation of the article}
The paper is organized as follows. In \autoref{section:preliminaries} we gather useful results for the rest of the paper. In \autoref{section:can_lift} we extend the canonical lifting theory of \cite[Appendix]{MS87} in two directions and present some further construction of (log) liftings.  
In \autoref{section:def_log_lifts} we define a notion of equisingular liftings for demi-normal varieties. In \autoref{section:lift_index_1} we construct equisingular liftings in case where the conductor does not contain a $\mathbb{P}^1$ passing through $4A_1$-singularities. 
That remaining case is studied in \autoref{section:4A_1}: since it is quite long, we have included a description of its content at its beginning.

\subsection*{Acknowledgments}
FB is supported by the grant PZ00P2-21610 from the Swiss National Science Foundation and QP is supported by the PostDoc Mobility grant P500PT-210980 from the Swiss National Science Foundation. 
We would like to thank  T. Kawakami, A. Petracci, R. Svaldi, E. Yasinsky
and J. Witaszek for useful conversations on the content of this article.
We are especially grateful to J. Baudin, for several conversations on abelian varieties.
We also thank the anonymous referees for several useful comments that improved the article.

%% file: 2_prelims.tex
\section{Preliminaries}\label{section:preliminaries}

\subsection{Notations and conventions}

\begin{enumerate} 
    \item Throughout the article, $k$ denotes an algebraically closed field of prime characteristic $p\geq 2$. At some places we will allow $p=0$ or restrict to $p>2$.
    \item A $k$-variety, resp.\ a $W(k)$-variety, is a flat separated, possibly disconnected, reduced equidimensional scheme of finite type over $k$, resp.\ over $W(k)$. 
    When it's clear from the context, we just say variety.
     We say $X$ is a curve (resp.~surface) over $k$ if $X$ is a $k$-variety of dimension one (resp.~two). 
    \item For a $k$-variety, we denote by $F \colon X \to X$ its (absolute) Frobenius morphism, and we use the same notation for $\mathbb{F}_p$-algebras.
    \item We say $(X,\Delta)$ is a \emph{couple} if $X$ is a normal variety (over $k$ or over $W(k)$) and $\Delta$ is an effective $\mathbb{Q}$-divisor on $X$. We say a couple $(X,\Delta)$ is a \emph{pair} if $K_X+\Delta$ is a $\mathbb{Q}$-Cartier $\mathbb{Q}$-divisor. 
    We say $\Delta$ is a \emph{boundary} if its coefficients lie in the interval $[0,1]$, and that $\Delta$ is reduced if its coefficients belong to $\{0,1\}$.
     \item Given a variety $X$, we say that $f \colon Y \to X$ is a log resolution if $f$ is a proper birational morphism such that $Y$ is regular and the exceptional locus $\text{Ex}(f)$ is purely divisorial and snc.  
     \item For the definition of the singularities of pairs appearing in the MMP (such as canonical, klt, dlt, log canonical and semi-log canonical) we refer to \cite[Definition 2.8]{kk-singbook}.
     \item A projective pair $(X, \Delta)$ is \emph{log Calabi--Yau} (or log CY) if $(X,\Delta)$ has log canonical singularities and $K_X+\Delta \equiv 0$.
     \item We say a variety $X$ is \emph{globally $F$-split} (or GFS) if $\mathcal{O}_X \to F_*\mathcal{O}_X$ splits as a $\mathcal{O}_X$-module homomorphism.
     We say a pair $(X,\Delta)$ is \emph{globally sharply $F$-split} 
     if there exists $e>0$ such that 
     $\mathcal{O}_X \to F^e_*\mathcal{O}_X(\lceil(p^e-1)\Delta \rceil)$ splits as a $\mathcal{O}_X$-module homomorphism. 
     \item Let us recall the construction of the related to $F$-splittings, called the trace map of Frobenius (see e.g. \cite[\S 4.2]{SS10}). Let $Y$ be a demi-normal variety and $\Delta$ be a reduced Mumford divisor on $Y$. If one applies $\Hom(\bullet,\sO_Y)=\Hom(\bullet,\sO(K_Y))[\otimes] \sO(-K_Y)$ to the natural map $t\colon\sO_Y\to F_*^e\sO_Y(\lceil{(p^e-1)\Delta\rceil})$, one obtains the trace map
        $$\Tr_{(Y,\Delta)}\colon F^e_*\sO_Y(\lfloor{(1-p^e)(K_Y+\Delta)\rfloor})\longrightarrow \sO_Y.$$
    Then left-splittings of $t$ (i.e. maps $s\colon \sO_Y(\lceil{(p^e-1)\Delta \rceil})\to \sO_Y$ satisfying $s\circ t=\id$) correspond under $\Hom(-,\sO_Y)$ to right-splittings of $\Tr_{(Y,\Delta)}$. Moreover, the right-splittings of $\Tr_{(Y,\Delta)}$ correspond to the elements of $H^0(Y,\mathcal{O}_Y{\lfloor{(1-p^e)(K_Y+\Delta)\rfloor}})$ which map to $1\in H^0(Y,\sO_Y)$ under $H^0(Y,\Tr_{(Y,\Delta)})$.
\end{enumerate}

\subsection{Adjunction on surfaces}
We highlight some basic facts about adjunction on log surfaces.

\begin{lemma}\label{lemma:everything_is_Q_Cartier}
Let $(S,B+\Delta)$ be a log canonical surface pair over $k$, where $B$ is a reduced curve. 
Then $S$ is $\mathbb{Q}$-factorial in a neighbourhood of $B$.
\end{lemma}

\begin{proof}
If $s\in S$ belongs to the support of $B$, then by \cite[Proposition 2.28.(2)]{kk-singbook} the local ring $\sO_{S,s}$ has rational singularities. By \cite[Proposition B.2]{Tan15} it follows that $\sO_{S,s}$ is $\mathbb{Q}$-factorial.
\end{proof}

\begin{lemma}\label{lemma:adjunction_on_surfaces}
Let $(S,B+\Delta)$ be a quasi-projective log canonical surface pair over $k$, where $B$ is a reduced curve. Then:
    \begin{enumerate}
        \item $B$ is a nodal curve (hence a Gorenstein curve);
        \item $\Delta\cap B$ does not contain any node of $B$;
        \item On $B$ we have the adjunction formula
                $$(K_S+B+\Delta)|_B=K_B+
                \underbrace{\sum_i \frac{n_i-1}{n_i}P_i
                +\sum_j P_j}_{\Diff_B(0)}
                +\Delta|_B, $$
            and on $B^\nu$ we have the adjunction formula
                $$(K_S+B+\Delta)|_{B^\nu}=K_{B^\nu}+
                \underbrace{ \sum_i \frac{n_i-1}{n_i}P_i
                +\sum_j P_j
                +\nu_B^{-1}\Sing(B)}_{\Diff_{B^\nu}(0)}
                +\Delta|_{B^\nu}$$
        where $n_i\in \mathbb{Z}_{>0}$ and the points $P_i,P_j$ are contained in the regular locus of $B$. Here  $\Delta|_{B^\nu}$ is the 
        pullback of $\Delta$ to $B^\nu$.
    \end{enumerate}
\end{lemma}
\begin{proof}
The first and second assertions are \cite[Theorem 2.31]{kk-singbook}. The third one is essentially a consequence of \cite[Proposition 2.36 and Theorem 3.36]{kk-singbook}, as we shall explain. It is assumed in \emph{op. cit.}, Theorem 3.36, that $B$ is proper, but this is not a necessary assumption. Indeed, we can realize $(S,B+\Delta)$ as an open subset of a projective lc surface pair $(\bar{S},\bar{B}+\bar{\Delta})$ and the adjunction maps commute with restrictions on open subsets, so the projective case implies the quasi-projective one.

The adjunction formulas $(K_S+B)|_B=K_B+\Diff_{B}(0)$ and $(K_S+B)|_{B^\nu}=K_{B^\nu}+\Diff_{B^\nu}(0)$ follow from \cite[Theorem 3.36 and preceding paragraph]{kk-singbook}. The fact that the $n_i$ are positive integers follows from the positive-definiteness of the matrices $\Gamma_p$ appearing in the statement of \cite[Theorem 3.36]{kk-singbook} (see \cite[2.26]{kk-singbook}). It remains to add the contribution of $\Delta$. By the second assertion and the local nature of the adjunction formulas, we may assume that $B$ is regular. By \autoref{lemma:everything_is_Q_Cartier} we may assume that $\Delta$ is $\mathbb{Q}$-Cartier. Then as explained in \cite[Proposition 2.36.(2) and its proof]{kk-singbook}, we have
        $$(K_S+B+\Delta)|_{B}=K_{B}+\Diff_{B}(0)+\Delta|_B. $$ This completes the proof.
\end{proof}

\begin{corollary}\label{corollary:possible_CY_log_curves}
Let $(S,B'+C)$ be a log CY surface pair with reduced boundary, where $C$ is irreducible. Then either:
    \begin{enumerate}
        \item $C$ is a rational curve disjoint from $B'$ with exactly one node, or
        \item $C$ is a regular curve, and $(C,\Diff_{C}(B'))$ is one of the following types:
            \begin{multicols}{2}
            \begin{itemize}
                \item $g(C)=1$ and $\Diff_C(0)$ is empty;
                \item $(\mathbb{P}^1,P+Q)$;
                \item $(\mathbb{P}^1,P+\frac{1}{2}Q_1+\frac{1}{2}Q_2)$;
                \item $(\mathbb{P}^1,\sum_{i=1}^3\frac{2}{3}Q_i)$;
                \item $(\mathbb{P}^1,\frac{1}{2}P+\frac{2}{3}Q+\frac{5}{6}R)$;
                \item $(\mathbb{P}^1,\frac{1}{2}P+\frac{3}{4}Q_1+\frac{3}{4}Q_2)$;
                \item $(\mathbb{P}^1,\sum_{i=1}^4\frac{1}{2}Q_i)$. 
            \end{itemize}
            \end{multicols}
        \noindent (Here it is understood that the points $P,Q,Q_i, R$ are always pairwise distinct when appearing in the same boundary.)
    \end{enumerate}
\end{corollary}
\begin{proof}
The CY condition implies that the restriction of $K_S+B'+C$ to $C^\nu$ has degree $0$. It follows from \cite[Theorem 2.31]{kk-singbook} that the term $(B'\cdot C)$ appearing in the adjunction formula is a reduced divisor on $C$. The cases where $C$ is a singular or regular curve of genus one follow easily and it remains to classify the possible coefficients of the pairs $(\mathbb{P}^1,\sum_i \frac{n_i-1}{n_i}P_i+\sum_j P_j)$ where $\sum_i \frac{n_i-1}{n_i}+\sum_j1=2$ for $n_i \in \mathbb{Z}_{>0}$. It is easily seen that the possibilities are the ones given in the statement.
\end{proof}


\begin{remark}\label{rmk:when_is_log_curve_GFS}
Every log CY structure on $\mathbb{P}^1$ appearing in \autoref{corollary:possible_CY_log_curves} above may or may not be globally sharply $F$-split, depending on the characteristic $p$ and on the boundary: see \cite[Theorem 4.2]{Watanabe_F_reg_and_F_pure_graded_rings} for the precise necessary and sufficient conditions. 
Notice that $(\bP^1,\sum_{i=1}^4\frac{1}{2}Q_i)$ is never globally sharply $F$-split if $p=2$.
\end{remark}

\subsection{Demi-normal varieties and $F$-splittings}\label{section:gluing_theory}
We give a quick review of gluing theory for demi-normal varieties (see \cite{Posva_Gluing_stable_families_surfaces_mix_char} and \cite[\S 5.1]{kk-singbook} for more details), and then we study how $F$-splittings ascend and descend normalizations.

Let $X$ be a demi-normal variety over $k$ or $W(k)$. 
If $p=2$ we assume that $X$ has only \emph{separable nodes}: see \cite[\S 3.1]{Posva_Gluing_for_surfaces_and_threefolds} for the definition.
Let $\nu_X\colon X^\nu\to X$ be the normalization, and let $C\subset X$ and $D\subset X^\nu$ be the conductor sub-schemes of $\nu_X$. Then:
    \begin{itemize}
        \item The commutative diagram
                \begin{equation}\label{eqn:pushout_diagram}
                \begin{tikzcd}
                D\arrow[r, hook]\arrow[d, "\nu_X|_D"] & X^\nu \arrow[d, "\nu_X"] \\
                C \arrow[r, hook] & X.
                \end{tikzcd}
                \end{equation}
        is both cartesian (as $I_C\sO_{X^\nu}=I_D$ by general properties of the conductor ideal) and co-cartesian (see \cite[Theorem 9.30, 10.18]{kk-singbook});
        \item $D$ is reduced of pure codimension one in $X^\nu$, and similarly for $C\subset X$;
        \item $\nu_X|_D\colon D\to C$ is separable of degree $2$ over each generic point of $C$: this induces a non-trivial involution $\tau$ of $D^\nu$.
    \end{itemize}
We refer to $(X^\nu,D,\tau)$ as the \emph{normalization triple} of $X$. We record:
\begin{lemma}\label{lemma:normal_conductors}
There is a big dense open subset $U$ of $X$ such that $C\cap U$ and $D \cap U^\nu$ are normal.
\end{lemma}
\begin{proof}
Let $Z=\nu_X(\Sing(D))\cup \Sing(C)$. Then $\codim_XZ\geq 2$ as $D$ and $C$ are reduced and $\nu_X$ is finite. So $U=X\setminus Z$ works.
\end{proof}

If furthermore $X$ has semi-log canonical singularities, then we also have:
    \begin{itemize}
        \item $(X^\nu,D)$ is a log canonical pair and $\nu^*K_X=K_{X^\nu}+D$;
        \item $\tau$ is a log involution of the pair $(D^\nu, \Diff_{D^\nu}(0))$.
    \end{itemize}

Let $\upsilon\colon D^\nu\to X^\nu$ be the composition of the normalization $D^\nu\to D$ with the embedding $D\hookrightarrow X^\nu$. Then the (closed) image of 
        $$(\upsilon, \upsilon\circ \tau)\colon D^\nu \longrightarrow 
        X\times X$$
generates a finite set-theoretic equivalence relation $R(\tau)\rightrightarrows X^\nu$ and $X=X^\nu/R(\tau)$ is the geometric quotient by this equivalence relation. 

\begin{remark}\label{rmk:index_change}
If $X$ is slc, then follows from \cite[Proposition 3.1.7]{Posva_Gluing_for_surfaces_and_threefolds} that the Cartier index of $K_X$ is either equal to, or twice, the Cartier index of $K_{X^\nu}+D$. 
\end{remark}

Next we explain how $F$-splittings interact with normalization of demi-normal varieties. We fix a demi-normal variety $X$ over $k$, with normalization $X^\nu\to X$ and associated conductor sub-schemes $D\subset X^\nu$ and $C\subset X$.
\begin{proposition}\label{prop:GFS_ascends_normalization}
If $X$ is globally $F$-split, then so is $(X^\nu,D)$. 
\end{proposition}
\begin{proof}
Let $\varphi$ be an $\sO_{X}$-linear splitting of the natural map $\sO_X\to F_*^e\sO_X$. Let $Q$ be the sheaf of total fractions of $X^\nu$: then $\varphi$ extends naturally and uniquely to a $Q$-linear map $\varphi\colon F_*^eQ\to Q$ by \cite[Proposition 7.11]{Schwede_Centers_F_purity}. We regard $\sO_{X^\nu}$ as a sub-sheaf of $Q$ and $F^e_*\sO_{X^\nu}((p^e-1)D)$ as a sub-sheaf of $F^e_*Q$, and we claim that 
    \begin{equation}\label{eqn:restriction_splitting}
    \varphi\left( F^e_*\sO_{X^\nu}((p^e-1)D)\right)\subseteq \sO_{X^\nu}.
    \end{equation}
Suppose that this is the case: then the composition 
        $$\sO_{X^\nu}\hookrightarrow F^e_*\sO_{X^\nu}((p^e-1)D)
        \overset{\varphi}{\longrightarrow} \sO_{X^\nu}$$
is an $\sO_{X^\nu}$-linear map that sends $1$ to $1$, so it must be the identity map.

Since the sheaves $\sO_{X^\nu}$ and $F^e_*\sO_{X^\nu}((p^e-1)D)$ are reflexive on $X^\nu$, it is sufficient to check \autoref{eqn:restriction_splitting} at every codimension one point $\eta$ of $X^\nu$. So let $s\in \sO((p^e-1)D)_\eta$. If $t$ is a local equation for $D$ at $\eta$, then we can write
        $$s=\frac{a}{t^{p^e-1}}, \quad a\in \sO_{X^\nu,\eta}.$$
Then
        $$\varphi(F^e_*s) =
        \varphi\left( \frac{1}{t}\cdot F^e_*(at)\right)
        = \frac{1}{t}\varphi(F^e_*(at)).$$
Now $at$ is an element of the conductor ideal $\sO(-D)$, and $\varphi$ sends $F_*^e\sO(-D)$ to $\sO(-D)$ by \cite[Proposition 7.10]{Schwede_Centers_F_purity}. As $\sO(-D)_\eta=t\sO_{X^\nu,\eta}$, we obtain that $\varphi(s)\in \sO_{X^\nu,\eta}$, as desired.
\end{proof}

\begin{corollary}
Assume that $X$ is a demi-normal globally $F$-split surface. If $K_X$ is $\bQ$-Cartier, then $X$ is semi-log canonical.
\end{corollary}
\begin{proof}
It follows from the above proposition that $(X^\nu,D)$ is globally $F$-split. By \cite[Proposition 2.2]{BBKW24} the couple $(X^\nu,D)$ is a log canonical pair. Since $\nu^*K_X=K_{X^\nu}+D$ and $K_X$ is $\bQ$-Cartier, we obtain that $X$ is slc.
\end{proof}

We note that when $p\neq 2$, then in dimension two $K_X$ is $\bQ$-Cartier if and only if $\Diff_{D^\nu}(0)$ is $\tau$-invariant \cite[Proposition 5.18]{kk-singbook}.

\begin{proposition}\label{prop:descent_1/p_linear_map}
Let $\varphi\colon F^e_*\sO_{X^\nu}((p^e-1)D)\to \sO_{X^\nu}$ be a $\sO_{X^\nu}$-linear map (resp.\ a global $F$-splitting of $(X^\nu,D)$). Then $\varphi$ restricts to a $\sO_{X}$-linear map $\varphi_0\colon F^e_*\sO_{X}\to\sO_{X}$ (resp.\ to a global $F$-splitting of $X$) if and only if the induced map $\bar{\varphi}^\nu\colon F^e_*\sO_{D^\nu} \to \sO_{D^\nu}$ sends $F_*^e(\sO_{D^\nu})^{F_*^e\tau}$ to $(\sO_{D^\nu})^\tau$.
\end{proposition}
\begin{proof}
The induced map $\bar{\varphi}^\nu\colon F^e_*\sO_{D^\nu} \to \sO_{D^\nu}$ is constructed as follows. There is an $\sO_D$-linear map $\bar{\varphi}\colon F_*^e\sO_D\to \sO_D$ making the diagram
    $$\begin{tikzcd}
    0 \arrow[r] & F_*^e\sO_{X^\nu}(-D)\arrow[r]\arrow[d] &
    F_*^e\sO_{X^\nu} \arrow[r]\arrow[d, "\varphi"] & F_*^e\sO_D \arrow[r]\arrow[d, "\bar{\varphi}"] & 0 \\
    0 \arrow[r] & \sO_{X^\nu}(-D) \arrow[r] &
    \sO_{X^\nu} \arrow[r] & \sO_D \arrow[r] & 0
    \end{tikzcd}$$
commutative, where the vertical map on the left exists by \cite[Proposition 7.10]{Schwede_Centers_F_purity}. By \cite[Proposition 7.11]{Schwede_Centers_F_purity}, $\bar{\varphi}$ extends to the desired $\sO_{D^\nu}$-linear map $\bar{\varphi}^\nu$.

Assuming that $\varphi_0$ exists, it is clear that $\varphi_0$ is a splitting of $\sO_X\to F_*^e\sO_X$ if $\varphi$ is a splitting of $\sO_{X^\nu}\to F_*^e\sO_{X^\nu}((p^e-1)D)$.

Now we prove the ``if and only if" statement. First assume that $\bar{\varphi}^\nu$ sends the sub-sheaf of $F^e_*\tau$-invariant elements to the sub-sheaf of $\tau$-invariants, and denote by $\bar{\varphi}^{\nu,\tau}\colon F^e_*(\sO_{D^\nu})^{F_*^e\tau}\to (\sO_{D^\nu})^\tau$ the corresponding (co-)restriction. Then we have a commutative diagram of solid arrows
    $$\begin{tikzcd}
    && \sO_{X^\nu}\arrow[rrr]\arrow[from=dddd, hook] &&& \sO_{D^\nu} \\
    &&&&& \\
    F^e_*\sO_{X^\nu}\arrow[rrr, crossing over]\arrow[uurr, "\varphi"] &&& F^e_*\sO_{D^\nu}\arrow[uurr, "\bar{\varphi}^\nu" below right] && \\
    &&&&& \\
    && \sO_{X}\arrow[rrr] &&& (\sO_{D^\nu})^\tau\arrow[uuuu, hook] \\
    &&&&& \\
    F^e_*\sO_{X}\arrow[rrr]\arrow[uuuu, hook]\arrow[uurr, dotted, "\varphi_0"] &&& F^e_*(\sO_{D^\nu})^{F^e_*\tau}\arrow[uuuu, hook, crossing over]\arrow[uurr, "\bar{\varphi}^{\nu,\tau}" below right]
    \end{tikzcd}$$
and we have to show the existence of $\varphi_0$ making the whole diagram commutative. We claim that the face at the back of the cube is a pullback diagram of sheaves. Since $\sO_X$ is $S_2$ it is sufficient to check it on a big open subset of $X$: then by \autoref{lemma:normal_conductors} we may assume that $D$ and $C$ are normal, and apply the fact that the square \autoref{eqn:pushout_diagram} is co-cartesian. The existence (and unicity) of $\varphi_0$ now follows from the universal property of pullbacks.

Conversely, assume that in the diagram above $\varphi_0$ exists and that we look for $\bar{\varphi}^{\nu,\tau}$. From \cite[Propositions 7.10]{Schwede_Centers_F_purity} we have $\varphi_0F^e_*I_C\subset I_C$, and as in the first paragraph it follows that $\varphi_0$ descends to $\bar{\varphi}_0\colon F^e_*\sO_{C}\to \sO_{C}$. This map in turn extends uniquely to the normalization of $C$. As $C^\nu=D^\nu/\langle \tau\rangle$, we obtain $\bar{\varphi}^{\nu,\tau}\colon F^e_*(\sO_{D^\nu})^{F^e_*\tau}\to (\sO_{D^\nu})^\tau$. By construction this map makes the cube a commutative diagram: hence $\bar{\varphi}^\nu$ preserves the sub-sheaf of $\tau$-invariant elements.
\end{proof}


\subsection{Finite actions and quotients}\label{section:finite_qts}
We collect some basic material on finite group actions. 

\begin{assumption}
We work over a base $B$ which is either $\Spec(k)$ or the spectrum of a DVR with residue field $k$---for example, $W(k)$. We let $X$ be a normal scheme that is flat, separated and of finite type over $B$. Suppose that $X$ is endowed with a good $B$-linear action of a finite flat $B$-group scheme $G$ (\footnote{
    A $G$-action on $S$ is \emph{good} if every point is contained in an affine $G$-stable open sub-scheme. If $X$ is quasi-projective over $B$, any finite group action on $X$ is good.
}).
\end{assumption}

We use the following terminology:

\begin{definition}
The \emph{ramification locus} of the $G$-action is the closed subset where the action is not free. Its image in $Y$ is the \emph{branch locus}.

Let $Z\subset X$ be a closed sub-scheme. We say that $Z$ is \emph{invariant} for the $G$-action if the sheaf of ideals $I_Z$ is a $G$-sub-module of $\sO_X$. We say that $Z$ is \emph{fixed} by the $G$-action if it is invariant and if the induced $G$-action on $\sO_Z$ is trivial.
\end{definition}

The geometric $G$-quotient $q\colon X\to Y=X/G$ always exists as a finite morphism onto a normal separated $B$-scheme of finite type (\footnote{
    Taking in account that the quotient stack $[X/G]$ has these properties and has finite stabilizers, this follows from the Keel--Mori theorem \cite[Theorem 6.12]{Rydh_Geometric_quotients}. 
}).
When $G$ is also linearly reductive (e.g.\ when $p$ does not divide the length of $G$), it turns out that taking quotients commutes with base-change over $Y$:

\begin{proposition}\label{prop:restriction_linearly_red_qt}
Suppose that $G$ is linearly reductive over $B$. Let $W\subset Y=X/G$ be a locally closed sub-scheme and consider the cartesian diagram
    $$\begin{tikzcd}
    X_W\arrow[d, "q_W"]\arrow[r, "f"] & X\arrow[d,"q"] \\
    W\arrow[r, hook] & Y.
    \end{tikzcd}$$
Then $X_W=X\times_YW$ has a good $G$-action which makes $f$ equivariant, and $q_W$ is the geometric quotient.
\end{proposition}
\begin{proof}
By definition $q^{-1}(W)=X_W$ is preserved by the action of $G$, so we obtain a $G$-action on $X_W$ which makes $f$ equivariant. Since $q_W$ is finite and $G$-invariant, it is clear that the $G$-action on $X_W$ is good. The rest of the statement follows from \cite[Theorem 3.2, Corollary 3.3]{AOV_Tame_stacks}.
\end{proof}


Projectivity over $B$ descends finite geometric quotients: 
\begin{proposition}\label{prop:descent_projectivity}
Suppose that $X$ is projective over $B$. Then the geometric quotient $X/G$ is also projective over $B$.
\end{proposition}
\begin{proof}
See \cite[Expos\'{e} 5, Remarque 5.1]{SGA3I} or \cite[Proposition 4.7]{Rydh_Geometric_quotients}.
\end{proof}


Let us now specialize to the case where $X$ is a smooth over $B$ and where $G=\mu_{2^m,B}$. 

\begin{lemma}\label{lemma:fixed_locus_is_smooth}
If $X$ is smooth over $B$ and $G=\mu_{2^m,B}$, then the ramification locus of the $G$-action is closed and smooth over $B$.
\end{lemma}
\begin{proof}
We may assume that the action is faithful, generated by $\sigma\in \Aut_B(X)$. Then the ramification locus of the $G$-action is the fixed locus of $\sigma^{2^{m-1}}$. This automorphism generates an action of $\mu_{2,B}$ on $X$. As $\mu_{2,B}$ is linearly reductive over $B$, the fixed locus of the $\mu_{2,B}$-action on $X$ is closed and smooth over $B$ by \cite[Proposition A.8.10.(2)]{CGP15}.
\end{proof}

For the rest of this sub-section, \emph{we assume that $B=\Spec W(k)$, that $\dim_B X=2$, that $G=\mu_{2m,B}$ and that $2$ is invertible on $B$}. Fix a generator $\sigma\in \Aut_{W(k)}(X)$ of the action; its ramification locus $\Ram(\sigma)$ is smooth over $W(k)$ by \autoref{lemma:fixed_locus_is_smooth}. We study the quotient $Y=X/G$ at the image of a ramification point.

Let $x\in \Ram(\sigma)$ be a closed point, and $y$ is its image in $Y$. If $x=x_1,\dots,x_l$ are the preimages of $y$ in $X$, then their stabilizers $G_i=\Stab_{G}(x_i)$ are conjugate to each other, the complete local rings $\widehat{\sO}_{X,x_i}$ are isomorphic to each other compatibly with the actions of the stabilizers, and
        \begin{equation}\label{eqn:red_to_stabilizer}
        \widehat{\sO}_{Y,y}\cong \widehat{\sO}_{X,x}^{\Stab_{G}(x)}.
        \end{equation}
So, replacing $G$ by $\Stab_G(x)$, we may assume that $x$ is a fixed point of $\sigma$.

If $x$ is a fixed point, then $\sigma$ descends to a $W(k)$-linear automorphism $\sigma_x$ of the local ring $\sO_{X,x}$. In particular, it induces a $k$-linear automorphism $\overline{\sigma}_{x}$ of order $2^m$ of the $k$-vector space $\fm_{X,x}/\fm_{X,x}^2$. By Cohen's structure theorem we have a $W(k)$-isomorphism 
    $$\widehat{\sO}_{X,x}\cong W(k)\llbracket u_1,u_2\rrbracket .$$
So $\fm_{X,x}/\fm_{X,x}^2$ is a three-dimensional $k$-vector space, generated by the images of $p,u_1,u_2$. The image of $p$ is fixed by $\overline{\sigma}_x$. Now $G=\mu_{2^m}$ is linearly reductive so the automorphism $\overline{\sigma}_x$ is diagonalizable with eigenvalues $1,\xi_1,\xi_2\in \mu_{2^m}(k)$. This implies that we can choose the formal coordinates $u_1,u_2$ such that the extension $\widehat{\sigma}_x$ of $\sigma_x$ to $\widehat{\sO}_{X,x}$ satisfies
        $$\widehat{\sigma}_x(u_i)=\xi_iu_i, \quad i=1,2.$$
The following lemma is standard:

\begin{lemma}\label{lemma:sing_qt}
With the reduction step as above:
    \begin{enumerate}
        \item the $\mu_{2^m}$-action on the connected component of $X$ containing $x$ is trivial if and only if $1$ is the only eigenvalue of $\overline{\sigma}_x$;
        \item $\dim_{W(k)}\Fix(\sigma)=0$ at $x$ if and only if both eigenvalues $\xi_1$ and $\xi_2$ are different from $1$;
        \item $y=q(x)$ is a singular point of $Y=X/\mu_{2^m}$ if and only if $\dim_{W(k)}\Fix(\sigma)=0$ at $x$.
    \end{enumerate}
\end{lemma}
\begin{proof}
The natural map $\sO_{X,x}\to \widehat{\sO}_{X,x}$ is equivariant and $\Spec(\sO_{X,x})$ is schematically dense in the connected component of $x\in X$: this implies the first point.

Let $I$ be the ideal of $\sO_{X,x}$ generated by $\im(\sigma_x-1)$, and $\mathcal{I}$ be the ideal of $\widehat{\sO}_{X,x}$ generated by $\im(\widehat{\sigma}_x-1)$: we claim that $I\otimes\widehat{\sO}_{X,x}=\mathcal{I}$.
Let  $\sO_n=\sO_{X,x}/\fm_x^n=\widehat{\sO}_{X,x}/\widehat{\fm}_x^n$, and notice that $\sigma_x-1$ and $\widehat{\sigma}_x-1$ descend to the same linear operator $\delta_n$ on $\sO_n$. So $I\otimes \sO_n$ and $\mathcal{I}\otimes \sO_n$ are generated by the image of $\delta_n$, which shows that $I\otimes \sO_n=\mathcal{I}\otimes \sO_n$. As $\mathcal{I}$ is a complete ideal, it follows that $I\otimes \widehat{\sO}_{X,x}=\mathcal{I}$. 
Since $I$ and $\mathcal{I}$ define respectively $\Fix(\sigma_x)$ and $\Fix(\widehat{\sigma}_x)$, we obtain
        $$\Fix(\sigma)\otimes \widehat{\sO}_{X,x}=\Fix(\widehat{\sigma}).$$
To prove the second point, it is now sufficient to prove that $\dim_{W(k)} V(\mathcal{I})=0$ if and only if $\xi_1$ and $\xi_2$ are different from $1$. If both are different from $1$, then $\mathcal{I}=(u_1,u_2)$. If, say, $\xi_1=1$ and $\xi_2\neq 1$ then $\mathcal{I}=(u_2)$. So the second point follows. 

To prove the third point, as regularity can be checked at completions, by \autoref{eqn:red_to_stabilizer} and the second point it is sufficient to show that $\widehat{\sO}_{Y,y}$ is singular if and only if both $\xi_i$ are different from $1$. If, say, $\xi_1=1$ then $\widehat{\sO}_{Y,y}\cong W(k)\llbracket u_1,u_2^r\rrbracket$ is regular, where $r>1$ is minimal with $\xi_2^r=1$. If both $\xi_1$ and $\xi_2$ are different from $1$, then the maximal ideal of $\widehat{\sO}_{Y,y}$ is generated over $W(k)$ by those $u_1^{r_1}u_2^{r_2}$ with $\xi_1^{r_1}\xi_2^{r_2}=1$. Then it is easily seen that the embedding dimension of $\widehat{\sO}_{Y,y}$ over $W(k)$ is greater than $2$, so it cannot be regular.
\end{proof}

\begin{proposition}\label{prop:sing_qt_pair}
Assume that the quotient $Y=X/\mu_{2^m}$ is smooth over $W(k)$. Let $H=\Ram(\sigma)\subset X$ be the ramification locus, and $K\subset Y$ be the branch locus of $q$.
    \begin{enumerate}
        \item $K$ and $H$ are smooth of pure dimension $1$ over $W(k)$.
        \item Let $E\subset X$ be a $\mu_{2^m}$-invariant smooth relative divisor, such that $H$ is not contained in $\Supp(E)$. Then $(X,E+H)$ and $(Y,q(E)+K)$ are relatively snc.
    \end{enumerate}
\end{proposition}
\begin{proof}
We may assume that the $\mu_{2^m}$-action is nowhere trivial. Since the assertions of the proposition are formal-local on $Y$, by \autoref{eqn:red_to_stabilizer} we may assume that the $\mu_{2^m}$-action on $X$ is free outside $H$, and that $H$ is in fact the fixed locus of the action.

We have $H=q^{-1}(K)$. Since $Y$ is regular, the sub-scheme $H$ has pure dimension $1$ over $W(k)$ by \autoref{lemma:sing_qt}, and as $q$ is finite the same holds for $G$. 

Let $x\in K$ be any closed point. By \autoref{lemma:sing_qt} the eigenvalues of $\overline{\sigma}_x$ are $1$ with multiplicity $2$, and $\xi\in \mu_{2^m}(k)$ with $\xi\neq 1$. Hence an equivariant formal-local model at $x$ is given by
        $$(x\in X)\cong\left(\mathbf{0}\in \widehat{\bA}^2_{u,v}\right), \quad 
        \sigma\colon (u,v)\mapsto (u,\xi v).$$
Then
        $$H=(v=0),\quad (q(x)\in Y)\cong \left(\mathbf{0}\in \widehat{\bA}^2_{u,v^r}\right), \quad 
        K=(v^r=0)$$
where $r>0$ is the smallest integer such that $\xi^r=1$. 

Suppose that $E$ is a $\mu_{2^m}$-invariant curve, smooth over $W(k)$, passing through $x$ and distinct from $H$. The equation of $E$ defines a non-zero element in $\fm_x/\fm_x^2$ which is an eigenvector for the automorphism $\overline{\sigma}_x$: hence, in our equivariant formal-local model, the curve $E$ is given by $(u=0)$ and its image $q(E)$ in $Y$ is given by $(u=0)$. 

It remains to prove that $C=q(E)$ (with its reduced structure) is $W(k)$-smooth away from $K$. As $E$ is $\mu_{2^m}$-invariant we have $E=q^{-1}(C)$ set-theoretically. Since $q$ is \'{e}tale above $Y\setminus K$, the morphism $q^{-1}(C)\to C$ is generically \'{e}tale, and so $q^{-1}(C)$ is generically reduced. Since $Y$ is regular, $C$ is Cartier, and so $q^{-1}(C)$ is $S_1$. Therefore $q^{-1}(C)$ is reduced, and this implies that $E=q^{-1}(C)$. Now smoothness descends \'{e}tale morphisms \cite[02K5]{stacks-project}: as $E$ is $W(k)$-smooth we obtain that $C$ is so away from $K$. This completes the proof.
\end{proof}

In case $G=\mu_{2,B}$, we also record:
\begin{proposition}\label{prop:simplify_mu2_sing}
Assume that $X$ is smooth of relative dimension $2$ over $W(k)$, and that $G=\mu_{2,B}$. If $b\colon X'\to X$ is the blow-up of the components of $\Ram(\sigma)$ of relative dimension $0$, then $b$ is $\mu_{2,B}$-equivariant and the ramification locus of the $\mu_{2,B}$-action on $X'$ has pure relative dimension $1$. Moreover, the $b$-exceptional divisors are fixed by the action.
\end{proposition}
\begin{proof}
In this case the ramification locus coincides with the fixed locus. So by \autoref{lemma:fixed_locus_is_smooth}, the components of relative dimension $0$ of fixed locus are \'{e}tale over $W(k)$. After an \'{e}tale extension of the base, we reduce to the case where these components are sections. This enables a formal-local computation: say $X=\bA^2_{u,v}$ and $\sigma\colon (u,v)\mapsto (-u,-v)$. Then a chart of the blow-up of $(u,v)$ is given by $\bA^2_{\bar{u},\bar{v}}$ with $u=\bar{u}$ and $v=\bar{u}\bar{v}$. The lift of the action is given by $\sigma\colon (\bar{u},\bar{v})\mapsto (-\bar{u},\bar{v})$, and the assertions of the statement are now easily verified.
\end{proof}


\subsection{Cyclic covers}
We briefly recall the construction of cyclic covers associated to sections of line bundles. Let $X$ be a normal separated Noetherian scheme, let $\sL$ be a line bundle on $X$, let $n>0$ be a positive integer, and let $s \colon \sO_X\to \sL^{\otimes n}$ be a global section. Such a global section gives an $\sO_X$-algebra structure to the finite $\sO_X$-module $\bigoplus_{i=0}^{n-1}\sL^{-i}$, see e.g.\ \cite[2.50]{km-book}, and we let
        $$X[\sL,n,s]=\Spec_X \left(\bigoplus_{i=0}^{n-1}\sL^{-i}\right)
        \overset{\pi}{\longrightarrow} X.$$
The $\mathbb{Z}/n\mathbb{Z}$-grading on the direct sum induces an action of the group scheme $\mu_n$ on $X[\sL,n,s]$, and $\pi$ is the geometric quotient.

Let $I$ be the image sheaf of the contraction map $\sL^{-n}\overset{\otimes s}{\longrightarrow}\sO_X$: it is the ideal of an effective $\mathbb{Z}$-divisor $G$ on $X$. In case the map $\otimes s$ is injective (which is equivalent to $s_\xi\in \sL^n_\xi$ not being zero for every generic point $\xi$ of $X$), we get an isomorphism $\sL^{n}\cong \sO(G)$ sending $s$ to $1$. In this case we say that $X[\sL,n,s]$ is the \emph{$n$-th cyclic cover of $X$ branched along $G$}.



\begin{lemma}\label{lemma:restriction_cyclic_cover}
Let $Y\subset X$ be a locally closed sub-scheme. Then we have a cartesian diagram
        $$\begin{tikzcd}
        Y[\sL_Y,n,s_Y] \arrow[d, "\pi_Y"] \arrow[r] & X[\sL,n,s]\arrow[d, "\pi"] \\
        Y\arrow[r, hook] & X
        \end{tikzcd}$$
where $\sL_Y$ is the pullback of $\sL$ to $Y$, and $s_Y\in H^0(Y,\sL_Y)$ is the restriction of $s$.
\end{lemma}
\begin{proof}
This follows immediately from functoriality of the relative spectrum.
\end{proof}

\begin{proposition}\label{prop:quotient_is_cyclic_cover}
Let $X$ be a regular separated scheme with a non-trivial $\mu_2$-action such that the geometric quotient $Y=X/\mu_2$ is regular. Assume that $2$ is invertible on $X$. Then the quotient morphism $q\colon X\to Y$ is a double cover of $Y$ branched along a regular divisor.
\end{proposition}
\begin{proof}
Let $q_*\sO_X=\sO_Y\oplus \sL^{-1}$ be the $\mu_2$-eigensheaves decomposition. The sheaf $\sL^{-1}$ is reflexive of rank one; as $Y$ is regular it follows that $\sL^{-1}$ is invertible. The $\sO_Y$-algebra structure of $q_*\sO_X$ is equivalent to a morphism $\sL^{-2}\to \sO_Y$, dually to a global section $s\colon \sO_Y\to \sL^2$. 
So we obtain an isomorphism $X\cong Y[\sL,2,s]$ over $Y$. The branch locus of this cover is a regular divisor by \cite[2.51]{km-book}.
\end{proof}

\subsection{Logarithmic vector fields and differentials}
Let us quickly review some basic facts about logarithmic derivations and differentials. 

Let $(X,\Upsilon)$ be an snc equi-dimensional scheme of finite type over $k$. The sheaf of logarithmic vector fields $T_X(-\log \Upsilon)$ is the sub-sheaf of the tangent sheaf $T_{X/k}$ consisting of those derivations which send the ideal sheaf of $\Upsilon$ into itself. Its dual is denoted 
        $\Omega_X^1(\log \Upsilon)=\sHom_X(T_X(-\log \Upsilon), \sO_X)$.
Both of these sheaves are locally free, and the stalks of $\Omega_X^1(\log \Upsilon)$ at a point $x$ are described as follows (see e.g. \cite[\S 2]{EV92}): if $x_1,\dots,x_n$ are local parameters at $x$ such that $\Upsilon=(\prod_{i=1}^jx_i=0)$, then $\Omega_X^1(\log \Upsilon)_x$ is generated as $\sO_{X,x}$-module by
        $$\frac{dx_1}{x_1},\dots,\frac{dx_j}{x_j},dx_{j+1},\dots,dx_n.$$
If $\dim X=1$, then obviously $\Omega_X^1(\log \Upsilon)=\omega_X(\Upsilon)$ and $T_X(-\log \Upsilon)=T_X(-\Upsilon)$. In dimension two we will need the following results:

\begin{proposition}\label{prop:restriction_for_log_vector_fields}
Assume that $\dim X=2$ and let $\Upsilon=\bigcup_i\Upsilon_i$ be the decomposition into (regular) irreducible components. Then there is a short exact sequence
    \begin{equation*}
    0\to T_X(-\Upsilon)\to T_X(-\log\Upsilon)\to \bigoplus_i T_{\Upsilon_i}(-\log (\Upsilon-\Upsilon_i)|_{\Upsilon_i})\to 0.
    \end{equation*}
\end{proposition}
\begin{proof}
The inclusion $T_X(-\Upsilon)\hookrightarrow T_X(-\log \Upsilon)$ is the natural one. To identify the cokernel we may work locally at an arbitrary closed point $x\in X$. If $x\notin\Upsilon$ the cokernel is trivial. Otherwise, choose local parameters $x_1,x_2$ such that $\Upsilon=(\prod_{i=1}^jx_i=0)$. 

If $j=1$, then $T_X(-\log \Upsilon)_x$ is generated by the derivations $x_1\partial_{x_1}, \partial_{x_2}$ while $T_X(-\Upsilon)_x$ is generated by $x_1\partial_{x_1}, x_1\partial_{x_2}$. So the quotient is the $\sO_{\Upsilon,x}$-module generated by the image of $\partial_{x_2}$, which is $T_{\Upsilon,x}$.

If $j=2$, let $\Upsilon_i$ be the branch of $\Upsilon$ cut out by $(x_i=0)$. 
Then $T_X(-\log \Upsilon)_x$ is generated by the derivations $x_1\partial_{x_1}$ and $x_2\partial_{x_2}$, while $T_X(-\Upsilon)$ is generated by $x_1x_2\partial_{x_1}$ and  $x_1x_2\partial_{x_2}$. So the quotient is the direct sum of the $\sO_{\Upsilon_1,x}$-module generated by the image of $x_2\partial{x_2}$, and of the $\sO_{\Upsilon_2,x}$-module generated by the image of $x_1\partial_{x_1}$: this direct sum is naturally isomorphic to $T_{\Upsilon_1}(-\log\Upsilon_2)_x\oplus T_{\Upsilon_2}(-\log \Upsilon_1)_x$.
\end{proof}

\begin{proposition}\label{prop:restr_for_twisted_log_vector_fields}
Assume $\dim X=2$ and that $(X, E+D)$ is snc. 
Let $D=\bigcup_i D_i$ be the decomposition into (regular) irreducible components. Then each $E|_{D_i}$ is an snc divisor on $D_i$, and we have a short exact sequence
        \begin{equation}\label{eqn:restr_for_twisted_log_vector_fields}
        0\to T_X(-\log E)(-D) \to T_X(-\log (D+E))
        \to \bigoplus_i T_{D_i}(- \log (E+D-D_i)|_{D_i})\to 0.
        \end{equation}
\end{proposition}   
\begin{proof}
The inclusion $T_X(-\log E)(-D) \hookrightarrow T_X(-\log (D+E))$ is the natural one. To identify the cokernel we may work locally at an arbitrary closed point $x\in X$. If $x\notin D$ the cokernel is trivial, and if $x\notin E$ we are in the situation of \autoref{prop:restriction_for_log_vector_fields}. Thus assume that $x\in E\cap D$, and choose local coordinates $x_1,x_2$ such that $D=(x_1=0)$ and $E=(x_2=0)$. Then $T_X(-\log(E+D))_x$ is generated by $x_1\partial_{x_1}$ and $x_2\partial_{x_2}$, while
$T_X(-\log E)(-D)$ is generated by $x_1\partial_{x_1}$ and $x_1x_2\partial_{x_2}$. So the quotient is the $\sO_{D,x}$-module generated by the image of $x_2\partial_{x_2}$, which is $T_{D}(-\log E|_D)_x$.
\end{proof}


\subsection{Involutions of $\mathbb{P}^1$}
Let $A$ be a local ring with residue field $k$.
We recall some well-known facts about automorphisms and involutions of $\bP^1_A$.
Any rational function of $\mathbb{P}^1_A$ that generates the function field over $\Frac(A)$ will be called a \emph{coordinate} of $\bP^1_A$.

\begin{lemma}\label{lemma:Aut_of_P1}
The Möbius transformations induce $\PGL_2(A)\cong \Aut_A(\mathbb{P}^1_A)$.
\end{lemma}
\begin{proof}
By \cite[Chapter 0, \S 5]{Mumford_GIT} the automorphism functor of $\mathbb{P}^1_A$ over $A$ is represented by the group scheme $\PGL_{2,A}$ via Möbius transformations. So it suffices to show that the $A$-valued points of $\PGL_{2,A}$ is the abstract group $\PGL_2(A)$. In general we have the exact sequence
        $$H^0(\mathbb{G}_{m,A},A)\to H^0(\GL_{2,A},A)\to H^0(\PGL_{2,A},A)\to H^1(\mathbb{G}_{m,A},A)\cong \Pic(A).$$
Since $A$ is local the group $\Pic(A)$ is trivial. Therefore $\PGL_{2,A}(A)$ is the cokernel of the diagonal map $A^\times \to \GL_2(A)$, which by definition is $\PGL_2(A)$.
\end{proof}

\begin{lemma}\label{lemma:PGL2_is_3_transitive}
The group $\Aut_A(\mathbb{P}^1_A)$ acts transitively on the set of ordered triplets $(P,Q,R)\subset \bP^1_A(A)$ with the property that the reductions $P_k,Q_k,R_k\in\bP^1_k(k)$ are pairwise distinct.
\end{lemma}
\begin{proof}
Let $P,Q,R\in \mathbb{P}^1_A(A)$ be three $A$-points. Pick a coordinate $x$ for $\mathbb{P}^1_A$ such that $P=[a:1]$, $Q=[b:1]$ and $R=[c:1]$ where $a,b,c\in A$. We claim that the matrix
        $$M=\begin{pmatrix}
        b-c & -a(b-c) \\
        b-a & -c(b-a)
        \end{pmatrix}$$
belongs to $\GL_2(A)$. This is equivalent to say that its determinant $\det(M)=(a-c)(b-c)(b-a)$ is invertible in $A$. By assumption the reduction of $\det(M)$ modulo the maximal ideal of $A$ is not zero, which shows that $\det(M)\in A^\times$. So $M$ defines an element of $\PGL_2(A)=\Aut_A(\bP^1_A)$, which is easily seen to send the triplet $(P,Q,R)$ to $([0:1],[1:1],[1:0])$. 
\end{proof}

\begin{remark}\label{rmk:PGL2_not_4transitive}
The action of any element of $\PGL_2(A)$ is completely determined by its effect on $3$ distinct $A$-points. Indeed, if $A$ is a field this is well-known and the general case follows by base-changing to $\Frac(A)$ as the homomorphism $\Aut_A(\bP^1_A)\to \Aut_{\Frac(A)}(\bP^1_{\Frac(A)})$ is injective.
\end{remark}

\begin{proposition}\label{prop:involutions_on_P1}
Assume that $2\in A^\times$ and let $P,Q,R\in \mathbb{P}^1_A(A)$ be three $A$-rational points whose reductions in $\bP^1_k(k)$ are pairwise distinct.
    \begin{enumerate}
        \item There exists a unique non-trivial $A$-involution of $\mathbb{P}^1_A$ that fixes $P$ and $Q$.
        \item There is a bijection between the $A$-involutions of $\mathbb{P}^1_A$ exchanging $P$ and $Q$, and the elements of $A^\times$.
        \item There exists a unique $A$-involution of $\mathbb{P}^1_A$ that exchanges $P$ with $Q$, and that fixes $R$.
        \item Assume that $A=k$ is a field. Then a non-trivial involution of $\mathbb{P}^1_k$ has exactly two distinct fixed points.
    \end{enumerate}
\end{proposition}
\begin{proof}
To prove the first three statements, choose a coordinate $x$ for $\mathbb{P}^1_A$ such that $P$ becomes the origin $[0:1]$, $R$ becomes $[1:1]$ and $Q$ becomes $[1:0]$ (which is possible by \autoref{lemma:PGL2_is_3_transitive}). We keep in mind that $A$-automorphisms of $\mathbb{P}^1_A$ are given by Möbius transformations with coefficients in $A$ (\autoref{lemma:Aut_of_P1}).
\begin{enumerate}
    \item  Any $A$-automorphism $\tau$ of $\mathbb{P}^1_A$ that fixes both $P$ and $Q$ is of the form $x\mapsto ax$, where $a\in A$. This is a non-trivial involution if and only if $a=-1$. 
    \item An $A$-automorphism $\tau$ that sends $P$ to $Q$ is of the form
            $$x\mapsto \frac{ax+b}{x},\quad a,b\in A.$$
    This is an involution if and only if $a^2x+ab=ax^2$. As $x$ is not algebraic over $A$, we deduce that $\tau$ is an involution if and only if $a=0$. Thus every $b\in A^\times$ defines an involution $\tau_b\colon x\mapsto b/x$ exchanging $P$ and $Q$, and every involution with this property is of such form.
    \item In the notations of the previous item, we have $\tau_b(R)=R$ if and only if $b=1$.
\end{enumerate}
To see the last statement, notice that an involution on $\bP^1_A$ has some fixed points (e.g.\ by applying the Riemann--Hurwitz formula to the geometric quotient morphism). Choosing a coordinate $x'$ such that $[1:0]$ is a fixed point, the involution is then given by $x'\mapsto -x'+b$, and the second fixed point is $[b/2:1]$.
\end{proof}

Next we consider lifts from $k$ to $W(k)$ of log involutions of the pairs $(\bP^1,\Delta)$ considered in \autoref{corollary:possible_CY_log_curves}. When $\Supp(\Delta)$ contains at most three points this is quite simple to realize; the remaining cases of \autoref{corollary:possible_CY_log_curves} will be studied in \autoref{section:can_lift}.

\begin{corollary}\label{cor:lift_involutions_on_P1}
Assume that $2\in k^\times$ and let $\sE\to \Spec W(k)$ be a proper flat morphism of relative dimension one where $\sE=\bigsqcup_i \sE_i$ is regular. Let $\Gamma=\bigsqcup_i\Gamma_i$ be an effective divisor on $\sE$, whose support is unramified over $\Spec W(k)$. Assume that:
    \begin{enumerate}
        \item $\sE_k=\bigsqcup_i \sE_{i,k}$ is a disjoint union of smooth rational curves, 
        \item the coefficients of $\Gamma$ belong to $\{1\}\cup\{\frac{n-1}{n}\mid n\geq 2\}$ and $|\Supp(\Gamma_k)\cap \sE_{i,k}|\leq 3$ for each $i$,
        \item $K_{\sE_k}+\Gamma_k\sim_\mathbb{Q} 0$, and
        \item there is a log involution $\tau$ of $(\sE_k,\Gamma_k)$ over $\Spec(k)$.
    \end{enumerate}
There there exists an involution $\Phi$ of $(\sE,\Gamma)$ over $\Spec W(k)$ that lifts $\tau$.
\end{corollary}
\begin{proof}
 The pairs $(\sE_i,\Gamma_i)$ are quite constrained. Indeed, since $\sE_{i,k}$ is rational and $\mathbb{P}^1_k$ is rigid, we obtain $\sE_i\cong \mathbb{P}^1_{W(k)}$. Since $\Supp(\Gamma)$ is unramified over $\Spec W(k)$, we obtain that $\Supp(\Gamma_i)$ is the disjoint union of at most three sections of $\sE_i\to\Spec W(k)$. 
 Taking in account the CY condition on the central fiber, we obtain by \autoref{corollary:possible_CY_log_curves} that $(\sE_{i,k},\Gamma_{i,k})$ is log isomorphic to one of the following log CY pairs:
        $$(\mathbb{P}^1_k,P_{k}+P'_{k})\quad\text{or}\quad 
        \left(\mathbb{P}^1_k,aP_{k}+bQ_{k}+cQ_{k}'\right).$$
Therefore each $(\sE_i,\Gamma_i)$ is log isomorphic to one of the following log CY pairs:
        $$(\mathbb{P}^1_{W(k)},P+P')\quad\text{or}\quad 
        \left(\mathbb{P}^1_{W(k)},aP+bQ+cQ'\right)$$
where $P,P'$ (respectively $P,Q,Q'$) are $W(k)$-rational points whose reduction modulo $p$ are pairwise distinct.

Now we lift the involution $\tau$. Since it has order two, we can break it as a collection of non-trivial auto-involutions $\tau_i\colon (\sE_{i,k},\Gamma_{i,k})\to (\sE_{i,k},\Gamma_{i,k})$ and log isomorphisms $\tau_{jl}\colon (\sE_{j,k},\Gamma_{j,k})\rightleftarrows (\sE_{l,k},\Gamma_{l,k})\colon \tau_{jl}^{-1}$ (where $j\neq l$). It is sufficient to lift each one separately.
    \begin{enumerate}
        \item First we lift the auto-involutions $\tau_i$ to involutions $\Phi_i\colon (\sE_{i},\Gamma_{i})\to (\sE_{i},\Gamma_{i})$. We use freely the results of \autoref{prop:involutions_on_P1}. 
            \begin{enumerate}
                \item If $\Gamma_{i,k}=aP_k+bQ_k+cQ_k'$, then at least two of the boundary coefficients have to be equal (otherwise $\tau_i$ would have three fixed points and hence be trivial). So we may assume that $\tau_i(P_k)=P_k$ and $\tau_i(Q_k)=Q_k'$. Let $\Phi_i\colon \sE_i\cong \sE_i$ be the unique $W(k)$-involution that fixes $P$ and exchanges $Q$ with $Q'$. Then necessarily $\Phi_{i,k}=\tau_i$.
                \item If $\Gamma_{i,k}=P_k+P'_k$ and $\tau_i(P_k)=P_k$, then from \autoref{prop:involutions_on_P1} we see as before that there is a unique lift of $\tau_i$ to an involution $\Phi_i$ of $(\sE_i,P+P')$.
                \item Finally, if $\Gamma_{i,k}=P_k+P'_k$ but $\tau_i(P_k)=P'_k$, we can find a $\tau_i$-fixed point $R_k\in \sE_{i,k}$. Let $R\in \sE_i$ be an arbitrary lift of $R_k$ (for example its Teichmüller lift), and take $\Phi_i$ to be the unique involution of $\sE_i$ that fixes $R$ and exchanges $P$ with $P'$. Then $\Phi_i$ necessarily reduces to $\tau_i$.
            \end{enumerate}
        \item Second we lift the log isomorphisms $\tau_{jl}\colon (\sE_{j,k},\Gamma_{j,k})\cong (\sE_{l,k},\Gamma_{l,k})$. 
            \begin{enumerate}
                \item If the boundaries have three components, then after relabelling we may assume that $\Gamma_{j,k}=aP_{j,k}+bQ_{j,k}+cQ_{j,k}'$, that $\Gamma_{l,k}=aP_{l,k}+bQ_{l,k}+cQ_{l,k}'$ and that $\tau_{jl}$ sends the triplet $(P_{j,k},Q_{j,k},Q'_{j,k})$ to $(P_{l,k},Q_{l,k},Q'_{l,k})$. By \autoref{rmk:PGL2_not_4transitive} the automorphism $\tau_i$ is unique with this property. By $3$-transitivity (\autoref{lemma:PGL2_is_3_transitive}) there is a $W(k)$-isomorphism $\Phi_{jl}\colon \sE_j\cong \sE_l$ such that $\Phi_{jl}\colon (P_j,Q_j,Q_j')\mapsto (P_l,Q_l,Q_l')$. Then $\Phi_{jl}$ necessarily reduces to $\tau_{jl}$ modulo $p$.
                \item Assume that $\Gamma_{j,k}=P_{j,k}+P_{j,k}'$ and $\Gamma_{l,k}=P_{l,k}+P'_{l,k}$. We may assume that $\tau_{jl}(P_{j,k})=P_{l,k}$ after relabeling. Take a distinct third point $R_{j,k}\in \sE_{j,k}$ with image $R_{l,k}=\tau_{jl}(R_{j,k})$. Choose arbitrary liftings $R_j$ and $R_l$ over $W(k)$ (for example the Teichmüller lifts). Then as before there is a unique lift $\Phi_{jl}\colon \sE_j\cong \sE_l$ of $\tau_{jl}$, defined by the condition $\Phi_{jl}\colon (P_j,P'_j,R_j)\mapsto (P_l,P'_l,R_l)$.
            \end{enumerate}
    \end{enumerate}
This proof is complete.
\end{proof}

\subsection{Equivariant MMP for surfaces}\label{section:equivariant_MMP_surfaces}
In this section, we let $k$ be an algebraically closed field $k$ of characteristic $p\geq 0$.

\begin{definition}\label{def:action_on_pairs}
Let $(S,\Delta)$ be a surface couple over $k$. We define $\Aut_k(S,\Delta)$ to be the sub-group of those automorphisms $\phi\in \Aut_k(S)$ such that $\phi_*\Delta=\Delta$ as $\bQ$-divisors. 
\end{definition}

Note that $\phi_*\Delta=\Delta$ is equivalent to say that $\phi^*$ sends $\sO_S(m\Delta)$ to itself for any integer $m\neq 0$ such that $m\Delta$ is integral.

Now let $(S,\Delta)$ be a projective irreducible surface pair over $k$, where $S$ and $\Supp(\Delta)$ are smooth and let $H\subset \Aut_k(S,\Delta)$ an (abstract) finite group of automorphisms. 
By \cite[Example 2.18]{km-book} and \cite{Prokhorov_Equivariant_MMP} (see also the proof of \cite[Theorem 2.25]{kk-singbook}), there is a sequence of birational $H$-equivariant contractions
        $$S=S_0\to S_1\to \dots \to S_n=S_{\min}^H$$
such that:
    \begin{enumerate}
        \item each $S_i\to S_{i+1}$ is the contraction of an $H$-extremal ray $C=\sum C_i$ of $\overline{\NE}(S_i)^H$, where $C_i$ are disjoint $(-1)$-curves (so in particular each $S_i$ is regular);
        \item with $\Delta_\text{min}^H$ being the pushforward of $\Delta$, either:
            \begin{enumerate}
                \item $K_{S_{\min}^H}+\Delta_\text{min}^H$ is nef, or
                \item there is an $H$-equivariant fibration $f\colon S_{\min}^H\to B$ such that $-K_{S_{\min}^H}-\Delta_\text{min}^H$ is $f$-ample, $\rk \Pic(S^H_\text{min}/B)^H=1$ and $\dim B<\dim S_{\min}^H$.
            \end{enumerate}
    \end{enumerate}
Let us record:
\begin{proposition}\label{prop:min_equiv_log_resolution}
Let $(T,\Delta)$ be a normal surface pair over $k$, and let $H\subset \Aut_k(T,\Delta)$ be a finite (abstract) group of automorphisms. 
Then there exists a morphism $\varphi\colon T'\to T$ such that:
    \begin{enumerate}
        \item $\varphi$ is a log resolution of $(T,\Delta)$, and it is an isomorphism over the locus where $T$ and $\Delta$ are regular;
        \item if $\Delta'=\varphi^{-1}_*\Delta$, then $H$ acts on $(T',\Delta')$ and $\varphi$ is $H$-equivariant;
        \item we have $K_{T'}+\Delta'+F=\varphi^*(K_T+\Delta)$ for an effective divisor $F$ that is $H$-invariant (that is, $h_*F=F$ for every $h\in H$).
    \end{enumerate}
\end{proposition}
\begin{proof}
Let $\varphi_0\colon T_0\to T$ be the morphism obtained by iterating the following process: blow-up the (reduced) singular locus, and normalize. It is well-known that this process eventually stops with a regular surface \cite[0BGJ]{stacks-project}. The action of $H$ preserves the (reduced) singular locus, so it lifts the blow-ups; it also lifts the normalizations, by the universal property thereof. Hence $\varphi_0$ is $H$-equivariant. Let $\Delta_0=(\varphi_0^{-1})_*\Delta$. Since $\varphi_0$ is an isomorphism above the generic points of $\Delta$, we see that $H$ acts on $(T_0,\Delta_0)$.

Let $Z$ be the finite set of points where $\Delta_0$ is singular. Since $H$ acts on $\Delta_0$, it preserves $Z$. So the blow-up of $Z$ is $H$-equivariant. By repeating this operation, we obtain a birational projective $H$-equivariant $\varphi_1\colon T_1\to T$ where $T_1$ is regular, $\Delta_1=(\varphi_1^{-1})_*\Delta$ is regular and $H$ acts on $(T_1,\Delta_1)$.

Now run an $H$-equivariant $K_{T_1}+\Delta_1$-MMP over $T$:
        $$\begin{tikzcd}
            T_1 \arrow[dr, "\varphi_1"] \arrow[rr, "h"]
            && T' \arrow[dl, "\varphi"] \\
            & T.
        \end{tikzcd}$$
Here all morphisms are $H$-equivariant, $T'$ is regular, $K_{T'}+\Delta'$ is $\varphi$-nef with $\Delta'=h_*\Delta_1=\varphi^{-1}_*\Delta$, and $H$ acts on $(T',\Delta')$. 
We see as in the proof of \cite[Theorem 2.25]{kk-singbook} that $\Delta'+\Exc(\varphi)$ is snc.
By the negativity lemma we can write
        $$K_{T'}+\Delta'+F=\varphi^*(K_T+\Delta)$$
for an effective divisor $F$. It remains to show that $F$ is $H$-invariant, and it suffices to show that the sheaf $\sO(mF)$ is $H$-invariant for some integer $m>0$ as observed after \autoref{def:action_on_pairs}. Take $m$ to be the Cartier index of $K_T+\Delta$. The $H$-action preserves the sheaf of rational $m$-differentials $\omega_T^m(m\Delta)$, its pullback by $\varphi$ (since $\varphi$ is equivariant), the line bundle $\omega_{T'}^{m}(m\Delta')$ and its inverse. Since
        $$\omega_{T'}^{-m}(-m\Delta')\otimes 
        \varphi^*\omega_T^{m}(m\Delta)\cong \sO(mF)$$
we obtain that the $H$-action preserves $\sO(mF)$. This concludes the proof.
\end{proof}

In case $\Delta=0$, we make the following definition:

\begin{definition} \label{invariant-def}
Let $S$ be a regular projective surface over $k$, and $H\subset \Aut_k(S)$ a finite sub-group. An \emph{$H$-conic bundle structure} on $S$ is an $H$-equivariant fibration $f\colon S\to C$ where $C$ is a regular projective curve, whose general fibers are $\bP^1$s and whose singular fibers are the union of two $(-1)$-curves meeting transversally at one point (eq. $-K_S$ is $f$-ample). 

If the $H$-action on $C$ is trivial, we say that $f$ is \emph{$H$-invariant}. 

If there not exist an $H$-equivariant $C$-morphism $S\to S'$ onto a regular projective surface, we say that $f$ is \emph{minimal}.
\end{definition}

\begin{remark}\label{rmk:minimal_conic_bundle}
An $H$-conic bundle $S\to C$ is minimal if and only if the $H$-orbit of any $(-1)$-curve contained in a fiber is not a disjoint union of $(-1)$-curves.
Thus any $H$-conic bundle dominates a minimal one.
\end{remark}

\begin{proposition}\label{prop:min_equiv_model}
Suppose that $\kappa(S)=-\infty$. Then either $S_{\min}^H$ is a del Pezzo surface with $\rk \Pic(S^H_{\min})^H=1$, or it admits a minimal $H$-conic bundle structure.
\end{proposition}
\begin{proof}
Since $\kappa(S)$ is negative, there is no birational model of $S$ on which the canonical divisor is nef. So according to \cite[Example 2.18]{km-book} and \autoref{rmk:minimal_conic_bundle}, one of the two cases of the statement occurs.
\end{proof}

\subsection{Automorphisms of projective bundles} \label{sub-section: aut_proj_bundles}
Consider a connected Noetherian scheme $S$ and a (possibly trivial) automorphism $\sigma\in \Aut(S)$. Let $\sE$ be a locally free coherent $\sO_S$-module, with associated projective bundle $\pi\colon \bP=\bP_S(\sE)\to S$. We let $\Aut_{S,\sigma}(\bP)$ be the set of those automorphisms $\sigma_{\bP}\in \Aut(\bP)$ satisfying $\sigma\circ \pi=\pi\circ \sigma_{\bP}$.

We investigate the structure of $\Aut_{S,\sigma}(\bP)$: the following proposition generalizes \cite[Lemma 3]{Maruyama_Automorphisms_of_ruled_surfaces}.

\begin{proposition}\label{prop:relative_autom_proj_bundles}
The set $\Aut_{S,\sigma}(\bP)$ is in bijection with the set of pairs $(\sL,[\sigma_\sE])$ where $\sL\in \Pic(S)$ and $[\sigma_\sE]$ is a class in the quotient $\Isom_S(\sigma^*\sE, \sE\otimes \sL)/H^0(S,\sO_S^*)$.
\end{proposition}
\begin{proof}
To begin with, suppose that $\sigma_\bP\in \Aut_{S,\sigma}(\bP)$ is given. For any point $s$ of $S$, this automorphism (co-)restricts to an isomorphism between fibers $\bP_{k(s)}\cong \bP_{k(\sigma(s))}$ of $\pi$. So $\sigma_\bP^*\sO_\bP(1)$ restricts on any fiber $\bP_{k(s)}$ to the sheaf $\sO_{\bP_{k(s)}}(1)$. Hence there is an isomorphism
    $$\alpha\colon \sigma_{\bP}^*\sO_\bP(1)\cong\sO_\bP(1)\otimes \pi^*\sL$$
for some uniquely defined line bundle class $\sL\in \Pic(S)$ (this follows from \cite[Ex. III.12.4]{Ha77}, see also \cite[4.2.7]{EGA2}). Taking in account the projection formula and the equalities $\sigma_\bP^*=(\sigma_{\bP}^{-1})_*$, $(\sigma^{-1})_*=\sigma^*$ and $\pi\circ \sigma_\bP^{-1}=\sigma^{-1}\circ \pi$, we obtain
    \begin{eqnarray*}
        \sE\otimes \sL &=& \pi_*(\sO_\bP(1)\otimes \pi^*\sL)\\
        &\cong & \pi_*\sigma_{\bP}^*\sO_\bP(1) \\
        &=& (\sigma^{-1})_*\pi_*\sO_\bP(1) \\
        &=& \sigma^*\sE.
    \end{eqnarray*}
The composition defines an isomorphism $\sigma_\sE\colon \sigma^*\sE\cong \sE\otimes \sL$. It depends on $\alpha$, which is uniquely defined up to multiplication by an element of $H^0(\bP,\sO_\bP^*)=H^0(S,\sO_S^*)$. Thus the class of $\sigma_\sE$ in the quotient $\Isom_S(\sigma^*\sE, \sE\otimes \sL)/H^0(S,\sO_S^*)$ is canonically associated to $\sigma_\bP$.

Conversely, let $(\sL,[\sigma_\sE])$ be as in the statement. Fix a representative $\sigma_\sE\in \Isom_S(\sigma^*\sE,\sE\otimes \sL)$ of $[\sigma_\sE]$. By functoriality, the automorphism $\sigma$ induces a cartesian diagram
        $$\begin{tikzcd}
        \bP_S(\sigma^*\sE)\arrow[d]\arrow[rr, "\phi_\sigma" above, "\sim" below] && \bP_S(\sE) \arrow[d, "\pi"] \\
        S\arrow[rr, "\sigma"] && S
        \end{tikzcd}$$
where the isomorphism $\phi_\sigma$ verifies $\phi^*_{\sigma}\sO_{\bP(\sE)}(1)\cong \sO_{\bP(\sigma^*\sE)}(1)$. Now, $\sigma_\sE$ induces an $S$-isomorphism $\phi_{\sigma_\sE}\colon \bP_S(\sE\otimes \sL)\cong \bP_S(\sigma^*\sE)$ such that $\phi_{\sigma_\sE}^*\sO_{\bP(\sigma^*\sE)}(1)=\sO_{\bP(\sE\otimes \sL)}(1)$. Notice that $\phi_{\sigma_\sE}$ depends only on the class $[\sigma_\sE]$ of $\sigma_\sE$. Furthermore, by \cite[Lemma II.7.9]{Ha77} we have an $S$-isomorphism $\phi_\sL\colon \bP_S(\sE)\cong \bP_S(\sE\otimes \sL)$ satisfying 
$\phi_\sL^*\sO_{\bP(\sE\otimes\sL)}(1)=\sO_{\bP(\sE)}(1)\otimes \pi^*\sL$. Therefore the composition 
    $$\sigma_\bP= \phi_\sigma\circ \phi_{\sigma_\sE}\circ \phi_\sL\colon \bP_S(\sE)\to \bP_S(\sE)$$ 
is an element of $\Aut_{S,\sigma}(\bP)$ satisfying $\sigma_\bP^*\sO_\bP(1)\cong \sO_\bP(1)\otimes \pi^*\sL$, and depending only on the class $[\sigma_\sE]$ of $\sigma_\sE$.

As these two operations are clearly inverse to each other, the proposition follows.
\end{proof}

%% file: 2.5_Canonical_liftings.tex
\section{Canonical liftings}\label{section:can_lift}
In this section we define and study \emph{canonical liftings} over $W(k)$ for some log Calabi--Yau pairs of dimension 1 or 2 over $k$. The paradigmatic case is that of ordinary (equivalently globally $F$-split) irreducible regular curves of genus $1$, which goes back to the work of Serre--Tate on deformations of ordinary abelian varieties, and for which refer to \cite[Appendix]{MS87}. We treat three additional cases that will be needed later: circles of $\mathbb{P}^1$s, split projective bundles over ordinary elliptic curves and the globally sharply $F$-split pair $(\bP^1,\frac{1}{2}\sum_{i=1}^4q_i)$.

\subsection{Canonical liftings of varieties with trivial log tangent bundle}\label{section:can_lifts_trivial_log_tangent_bdl}
We recall the construction and the properties of the canonical liftings of pairs $(X,D)$ with trivial log tangent bundle, generalizing \cite{MS87}*{Appendix} to the logarithmic case as indicated in \cite{AWZ21}*{Variant 3.3.2}. 
For the purposes of this article, the two main examples to keep in mind are ordinary curves of genus 1 and the normalisation of a cycle of $\mathbb{P}^1$s.
We recall the definition of Frobenius log lifting for snc pairs.

\begin{definition}
Let $(X, D)$ be an snc pair over $k$. 
A \emph{log Frobenius lifting} (or log $F$-lifting) of $(X,D)$ over $W_n(k)$ is a triple $(X_n, D_n, F_n)$ where $(X_n, D_n)$ is a log lifting of the pair $(X,D)$ and $F_n$ is a lifting of the Frobenius morphism such that $F_n^*D_n=pD_n$.
We also call $F_n$ a log Frobenius lifting of the Frobenius $F$ on $(X,D)$ over $W_n(k)$. 
If clear from the context, we omit the adjective log.
\end{definition}

Recall for any regular variety $X$ over $k$, we define the sheaf of exact 1-forms as $B_X^1 \coloneqq F_{*}\mathcal{O}_X/\mathcal{O}_X$ which is a vector bundle on $X$ of  rank $\dim X -1$. If $X=\Spec(A)$, we use $B_A^1$ to denote $B_X^1$.

\begin{proposition} \label{prop: def_th_log_Frobenius}
    Let $(X, D)$ be an snc pair over $k$ and let $(X_n, D_n, F_n)$ be a log Frobenius lifting of $(X, D)$ over $W_n(k)$.
    Then the following holds:
    \begin{enumerate}
        \item For every log lifting $(X_{n+1}, D_{n+1})$ over $W_{n+1}(k)$, there is an obstruction class
        $$o^{F}_{(X_{n+1}, D_{n+1})} \in \Ext^{1}(\Omega_X^1(\log D), F_*\mathcal{O}_X),$$
        whose vanishing is sufficient and necessary for the existence of a log lifting $F_{n+1}$ of $F_{n}$ to $(X_{n+1}, D_{n+1})$. 
        If the obstruction vanishes, then the space of such log liftings is a torsor under $\Hom(\Omega_X^1(\log D), F_*\mathcal{O}_X)$.
        \item There exists an obstruction class
    $o_{(X_{n}, D_n, F_n) } \in \Ext^1(\Omega_X^1(\log D), B^1_X)$ whose vanishing is a sufficient and necessary condition for the existence of a log Frobenius lifting of $(X_n, D_n, F_{n})$ over $W_{n+1}(k)$. 
If the obstruction vanishes, then the space
of such liftings is a torsor under $\Hom(\Omega_X^1(\log D), B^1_X)$.
\end{enumerate}
\end{proposition}

\begin{proof}
    The case $D=0$ is proved in \cite[Appendix, Proposition 1]{MS87}. We indicate how to adapt their proof to the logarithmic setting. 
    We start with the description of the local case.
    
    Let $A$ be a regular affine $k$-algebra of finite type, and let $B$ be a flat $W_n(k)$-algebra such that $B/p \simeq A$.
    Let $x_1, \dots, x_r \in B$ be elements such that $D_B=(x_1 \cdots x_r=0)$ is a lift of $D$ in $X=\Spec(A)$ to $X_n=\Spec(B)$, and let $F_B$ a lifting of the Frobenius to $B$ compatible with $D$.

    We prove that, fixed a log lifting $(\Spec(C), D_C)$ over $W_{n+1}(k)$ of $(\Spec(B), D_B)$, the set of log lifting of the Frobenius $F_B$ to $C$ is a torsor under the group $\Hom_{A}(\Omega_A^1(\log D), F_*\mathcal{O}_A).$ 
    Indeed, given two lifts $F_1, F_2$ of $F_B$ to $C$ we have
    $F_1(b)=F_2(b)+p^n\psi({b})$, where $\psi \colon A \to A$ is a function satisfying    
    \begin{enumerate}
        \item[(i)] $\psi(b_1+b_2)=\psi(b_1)+\psi(b_2)$;
        \item[(ii)] $\psi(b_1\cdot b_2)=b_1^p\psi(b_2)+b_2^p\psi(b_1).$
    \end{enumerate}
    This shows that $\psi$ is a derivation with values in $F_*A$, thus
    by the universal property of the K\"ahler differentials, $\psi$ corresponds to an element $\eta \in \Hom_A(\Omega^1_{A},F_*A).$
    By hypothesis we have that $F_{n+1}(x_i)=x_i^p,$ which implies that $\psi(x_i)=0$, equivalently that $\eta \in \Hom_A(\Omega^1_{A}(\log D),F_*A).$ 
    We now prove (a). Fix a lifting $(X_{n+1}, D_{n+1})$ over $W_{n+1}(k)$ and let $(U_{n+1}, D_{n+1}|_{U_{n+1}})$ be an affine covering of $X_{n+1}$ We write $U_n \coloneq U_{n+1} \times X_n$. 
    By the local computations above, the set of log Frobenius liftings of $F_{U_n}$ is a torsor under the action of the group $\Hom_{\mathcal{O}_{X}}(\Omega_X^1(\log D), F_*\mathcal{O}_X).$ 
    Thus the obstruction to a lifting $F_{n+1}$ is a class in $\Ext^1(\Omega_X^1(\log D), F_*\mathcal{O}_X)$. 

    To prove (b), we start with the analysis of the local case.
 Let $\varphi \colon (\Spec(C_1), D_{C_1}) \to (\Spec(C_2), D_{C_2})$  be an isomorphism of pairs such that it restricts to the identity on $B$ and   let $F_1$ (resp. $F_2$) be a Frobenius log lifting of $F_B$ to $C_1$ (resp. $C_2$) over $W_{n+1}(k)$.
    As in the local computation above, there exists  $\eta \in \Hom_A(\Omega_A^1(\log D), F_* \mathcal{O}_A)$ such that 
    $$\varphi \circ F_1 - F_2 \circ \varphi =p^n \eta. $$
    Let $\varphi' \colon (\Spec(C_1), D_{C_1}) \to (\Spec(C_2), D_{C_2})$ be another isomorphism of pairs.
    As $\varphi'=\varphi+p^n\delta$ for some $\delta \in \Hom_A(\Omega_A^1(\log D), \mathcal{O}_A)$, a direct computation shows that 
    $$\varphi' \circ F_1 - F_2 \circ \varphi'=p^n \eta', \text{ where } \eta'=\eta+\delta^p,$$
    where $\delta^p$ is the image of $\delta$ via the natural map $\Hom_A(\Omega_A^1(\log D), A) \to \Hom_A(\Omega_A^1(\log D), F_*A)$.
    This shows that there exists at most a unique $\varphi$ for which $\varphi \circ F_1 = F_2 \circ \varphi$ and that the reside class $\overline{\eta} \in \Hom_A(\Omega_A^1(\log D), B_A^1)$ does not depend on the chosen isomorphism $\varphi$ of log pairs (but only on $F_1$ and $F_2$).
    We can now repeat the exact same argument of \cite[Proposition 2, (vii)]{MS87} to conclude the obstruction class to a lifting $(X_{n+1}, D_{n+1}, F_{n+1})$ exists and it is a class in $\Ext^1(\Omega_X^1(\log D), B_X^1)$. 
\end{proof}

We recall the existence of lifting of line bundles compatible with the Frobenius lifting. 

\begin{proposition} \label{prop: lift_can_linebundles}
    Let $X$ be a reduced connected scheme over $k$ and let $(X_n, F_{n})$ be a Frobenius lifting of $X$ over $W_n$.
    Let $L_n$ be a line bundle on $X_n$ such that $F_n^*L_n=L_n^{\otimes p}$
    Suppose that the
Frobenius action on $H^i
(X, \mathcal{O}_X)$ is bijective for $i = 1, 2$. Then  there exists a unique $L_{n+1} \in \Pic X_{n+1}$ such that
$F_{n+1}^{*}L_{n+1} \simeq L_{n+1}^{\otimes p}$ and $L_{n+1}|_{X_n} \simeq L_n$. 
\end{proposition}

\begin{proof}
    The proof of \cite[Proposition 2, Appendix]{MS87} applies in this generality.
\end{proof}

The following is an analogue of \cite{MS87}*{Appendix, Proposition 3} to the logarithmic case.

\begin{proposition}\label{prop: lift_morphisms}
    Let $(X_i, D_i)$ be smooth snc pairs over $k$ and let $(Z_i, D_{Z_i}, F_i)$ be Frobenius log lifting over $W_{n+1}(k)$.
    Suppose $(Y_i, D_{Y_i})$ be the restriction of $(Z_i, D_{Z_i})$ to $W_{n}(k)$.
    Let $\psi \colon (Y_1, D_{Y_1}) \to (Y_2, D_{Y_2})$ be morphism compatible with Frobenius lifting, i.e. $F_2 \circ \psi=\psi \circ F_1$.
    Then the obstruction of lifting $\psi$ to a morphism $\chi \colon (Z_1, D_{Z_1}) \to (Z_2, D_{Z_2})$ respecting the log structure and compatible with the lift of Frobenius is a class in $H^0(X_1, \varphi^{*}T_{X_2}(-\log D_2) \otimes B^1_X)$, whose vanishing is a sufficient and necessary condition of the existence of a lifting of $f$.  
    Moreover, if $f$ exists it is unique.
\end{proposition}

\begin{proof}
    The proof follows \cite{MS87}*{Appendix, Proposition 3} where, as in the proof of \autoref{prop: def_th_log_Frobenius}, we need to consider the compatibility with the boundaries $D_i$ in the computation of the obstruction classes.
\end{proof}

With this said, we can construct canonical liftings of projective pairs with trivial log tangent bundle, generalising \cite{MS87}*{Appendix, Theorem 1}.

\begin{theorem} \label{thm: canonical-log-liftings}
    Let $(X,D)$ be a projective snc pair such that $\Omega_X^1(\log D)$ is trivial and $X$ is globally $F$-split.
    Then there exists a projective lifting $(\mathbf{X}, \mathbf{D})$ over $W(k)$ together with a lifting of the Frobenius compatible with $D$. 
    Such lifting is called the \emph{canonical (log) lifting} of $(X,D)$.
    Moreover,
    \begin{enumerate}
        \item Let $\varphi \colon (X_1, D_1) \to (X_2, D_2)$ be a morphism of snc pairs  with trivial log tangent bundle such that $(X_i, D_i)$ are globally $F$-split. Then there exists a unique morphism between the canonical lifting $\varphi\colon (\mathbf{X}_1, \mathbf{D}_1) \to (\mathbf{X}_2, \mathbf{D}_2)$ which is a lift of $\varphi$ and that is compatible with the lift of the Frobenii morphisms.
        \item The restriction map $\Pic(\mathbf{X})_{F}=\left\{ L \in \Pic(\mathcal{X}) \mid F_{\mathcal{X}}^* L \simeq L^{\otimes p} \right\} \to \Pic(X)$ is an isomorphism.
    \end{enumerate}
\end{theorem}

\begin{proof}
    We first show the existence and uniqueness of a formal lifting $(\mathfrak{X}, \mathfrak{D})$ over $\text{Spf}(W(k))$.     
    As $X$ is globally $F$-split, we have that $H^i(X, B_X^1)=0$ for every $i \geq 0$. 
    By \autoref{prop: def_th_log_Frobenius}, the obstructions to the existence of Frobenius log lifting from $W_{n}(k)$ to $W_{n+1}(k)$ lie in $H^1(X, T_X(-\log D) \otimes B^1_X) \simeq H^1(X, B_X^1)^{\oplus \dim X}=0$, showing the existence of a formal lifting.
    Moreover the space of such liftings is a torsor over $H^0(X, T_X(-\log D) \otimes B^1_X) \simeq H^0(X, B_X^1)^{\oplus \dim X}=0$, which shows the uniqueness of the formal lifting. 
    We first prove (b) at the formal level, i.e. $\Pic(\mathfrak{X})_F \to \Pic(X)$ is an isomorphism. This follows immediately from \autoref{prop: def_th_log_Frobenius}{(c)}. Therefore there exists a lift of an ample line bundle $L$ and thus by Grothendieck algebraization theorem there exists a projective scheme $(\mathcal{X}, \mathcal{D})$ whose completion is $(\mathfrak{X}, \mathfrak{D})$ and every formal lift of a line bundle is algebraic.
    
    We now prove (a). By \autoref{prop: lift_morphisms}, the obstructions to the existence of a Frobenii compatible lifting of $\varphi$ lie in 
    $H^0(X_1, \varphi^*T_{X_2}(-\log D_2) \otimes B_X^1) = H^0(X_1,  B_X^1)^{\oplus \dim X_2} $. Thus $\varphi$ lifts to $\widehat{\varphi}$ at the formal level, which is algebraizable.
\end{proof}

As an application, we construct a Frobenius lifting of a circle $E$ of $\mathbb{P}^1$s over the Witt vectors. 

\begin{proposition} \label{prop: canonical_lifts_circlesP1}
Let $E=E_0\cup E_1 \cup \dots \cup E_n$ be an  oriented cycle of smooth rational curves $E_i \simeq \mathbb{P}^1_k$ over $k$.
Then there exists a log lifting $\mathbf{E}=\bigcup_i \mathbf{E}_i$ of $E$ over $W(k)$ such that:
\begin{enumerate}
    \item $E$ is globally $F$-split;
    \item there is a unique lift of Frobenius $F_{\mathbf{E}}$ on $\mathbf{E}$, such that it fits in the following commutative diagram 
    \begin{equation*}
       \begin{tikzcd}
    (\mathbf{E}^{\nu}, \mathbf{D}) \arrow[d, "\widetilde{\nu}"]  \arrow[r, "F_{(\mathbf{E}^\nu, \mathbf{D})}"] & 
    (\mathbf{E}^{\nu}, \mathbf{D})\arrow[d, "\widetilde{\nu}"] \\
    \mathbf{E} \arrow[r, "F_\mathbf{E}"]  & \mathbf{E} ,
    \end{tikzcd}   
    \end{equation*}
    where $F_{(\mathbf{E}^\nu, \mathbf{D})}$ is the log lifting of the Frobenius of \autoref{thm: canonical-log-liftings};
    \item for every automorphism $\sigma \colon E \to E$, there exists a unique automorphism $\widetilde{\sigma} \colon \mathbf{E} \to \mathbf{E}$ such that $\widetilde{\sigma} \circ F_{\mathbf{E}}=F_{\mathbf{E}} \circ \widetilde{\sigma}$;
    \item for every line bundle $L$ on $E$, there is a unique lift $\mathbf{L}$ such that $F_{\mathbf{E}}^*\mathbf{L} \simeq \mathbf{L}^{\otimes p}$. 
\end{enumerate}
We call $\mathbf{E}$ the \emph{canonical lifting} of $E$.
\end{proposition}
\begin{proof}
    We fix some notations. We write 
        $$E_i\cap E_{i+1}=p_{i+1}^i\in E_i, \quad 
        E_i\cap E_{i+1}=p^{i+1}_i\in E_{i+1}$$ 
    where the indices are taken modulo $n$. We let $(E^\nu, D)=\bigsqcup_i (E_i, p^i_{i-1}+p^i_{i+1})$ be the normalisation of $E$ with the conductor $D$. 

    We prove (a). As $(E^\nu, D)$ is a disjoint union of $(\mathbb{P}^1, [1:0]+[0:1])$, this pair is globally $F$-split with a unique splitting $\phi^{\nu}_{\bP^1}$. The induced map
            $$\phi_{\bP^1}^\nu|_{[1:0]} \colon F_*k([1:0])\to k([1:0])$$
    is the Frobenius morphism of the residue field $k([1:0])=k$.
    We obtain a global $F$-splitting $\phi^\nu$ of the pair $(E^\nu,D)$, and an analogous description of its restriction $\phi^\nu|_D\colon F_*\sO_D\to \sO_D$. The gluing involution $\tau\in \Aut_k(D)$ is the disjoint sum of the identifications $k(p^i_{i+1})=k(p^{i+1}_i)$ given by the structural $k$-module structure. So the hypothesis of \autoref{prop:descent_1/p_linear_map} are clearly satisfied, and $\phi^\nu$ restricts to a global $F$-splitting $\phi\colon F_*\sO_E\to \sO_E$.
    
    We prove (b). As $K_{E^\nu}+D \sim 0$ and dimension is 1, we can apply \autoref{thm: canonical-log-liftings} to each connected component $(E_i, D_i)$ to construct the canonical log lifting $(\mathbf{E}^\nu, \mathbf{D}) \coloneqq \bigsqcup (\mathbf{E}_i, \mathbf{D}_i)$ with a unique lifting of the Frobenius $F_{(\mathbf{E}^\nu, \mathbf{D})}$. We have
            $$\sD_i=\mathbf{p}^i_{i-1}\sqcup \mathbf{p}^i_{i+1}$$
    where $\mathbf{p}^i_{i-1}$ and $\mathbf{p}^i_{i+1}$ are $W(k)$-points lifting respectively $p^i_{i-1}$ and $p^i_{i+1}$. We define a gluing involution $\widetilde{\tau}$ on $\sD$ as the disjoint sum of the identifications
            $$\mathbf{p}^i_{i+1}= \mathbf{p}^{i+1}_i
            \quad \text{commuting with the projection to } W(k).$$
    It is clear that the quotient $\mathbf{E}=\mathbf{E}^\nu/R(\widetilde{\tau})$ exists and is a lift of $E$.
    
    We now show we can descend the lifting of the Frobenius to $\mathbf{E}$.
    To simplify the notations, say we work in a neighbourhood of $\mathbf{E}_1\cap \mathbf{E}_2$. Let $\varphi_i \colon \mathbf{E}_i \to \mathbf{E}$ be the natural projection maps. 
    We have isomorphisms
        $$\begin{tikzcd}
        \mathbf{p}^2_1\arrow[r, "\varphi_2" above, "\sim" below]
        & \mathbf{E}_1\cap\mathbf{E}_2 &
        \mathbf{p}^1_2\arrow[l, "\varphi_1" above, "\sim" below]
        \end{tikzcd}$$
    and note that the projection $\mathbf{E}\to \Spec(W(k))$ restricts to an isomorphism $\mathbf{E}_1\cap \mathbf{E}_2\cong \Spec(W(k))$.
    By construction (see also \cite[Corollary 2.3.9]{Posva_Gluing_stable_families_surfaces_mix_char}), a section $s\in \sO_\mathbf{E}$ on that neighborhood is the same as two sections $s_i\in \sO_{\mathbf{E}_i}$ such that 
            $$\varphi_2(s_2|_{\mathbf{p}^2_1})=
            \varphi_1(s_1|_{\mathbf{p}^1_2}).$$ 
    This implies that an endomorphism $\alpha$ of $\sO_{\mathbf{E}^\nu}$ descends to an endomorphism of $\sO_\mathbf{E}$ if and only if for every section $s=(s_1,s_2)\in \sO_{\mathbf{E}^\nu}$, we have
            $$\varphi_2(\alpha(s_2)|_{\mathbf{p}^2_1})=
            \varphi_1(\alpha(s_1)|_{\mathbf{p}^1_2}).$$
    Consider $\alpha=F_{(\mathbf{E}^\nu,\mathbf{D})}$: it restricts on both $\mathbf{p}^1_2$ and $\mathbf{p}^2_1$ to a lift of the Frobenius of $k$. As $k$ is perfect, there is a unique lifting of Frobenius on $W(k)$, so we conclude that $F_{(\mathbf{E}^\nu, \mathbf{D})}$ descends to a morphism $F_{\mathbf{E}}$ on $\mathbf{E}$ which lifts the Frobenius of $E$.

    We prove (c). Let $\sigma^{\nu} \colon (E^{\nu}, D) \to (E^\nu, D)$ be the log isomorphism induced by $\sigma$ on the normalisation.
    By \autoref{thm: canonical-log-liftings}, there is a unique extension $\widetilde{\sigma}^{\nu} \colon (\mathbf{E}^\nu, \mathbf{D}) \to (\mathbf{E}^\nu, \mathbf{D})$ of this isomorphism which is compatible with $F_{(\mathbf{E}^\nu,\mathbf{D})}$. 
    As above, one can show that $\widetilde{\sigma}^{\nu}$ descends to an endomorphism $\widetilde{\sigma}$ of $\mathbf{E}$ lifting $\sigma$.

    We now prove (d). By (a), $E$ is globally $F$-split and therefore we can apply \autoref{prop: lift_can_linebundles} to conclude there exists a unique lifting of $L$ which belongs to $\Pic(\mathbf{E})_{F_{\mathbf{E}}}$.
\end{proof}

\subsection{Canonical lifting of split projective bundles over ordinary genus 1 curves}\label{section:can_lift_autom_proj_bundles}

We construct canonical liftings of split projective bundle over ordinary elliptic curves.
Let $E$ be a globally $F$-split genus 1 curve, and let $M, N$ be two line bundles on $E$. 
Let $\mathbb{P}(M \oplus N)$ be the projective bundle associated to $M \oplus N$  together with the two naturally associated sections $C$ and $D$. Note that $(\mathbb{P}_E(M+C), C+D)$ is a log smooth CY surface.
Let $\mathbf{E}$ be the canonical lifting of $E$ to $W(k)$ and let $\mathcal{M}, \mathcal{N}$ be the canonical lifting of $M$ and $N$ on $\mathbf{E}$ given by \autoref{thm: canonical-log-liftings}.

\begin{definition} \label{def: canonical_proj_bundle_split_ell}
    We call $(\mathbb{P}_{\mathbf{E}}(\mathcal{M} \oplus \mathcal{N}), \mathbf{C}+\mathbf{D})$ the \emph{canonical lifting} of $(\mathbb{P}_E(M+C), C+D)$.
\end{definition} 

The adjective canonical is justified by the following proposition. We use the notations of \autoref{sub-section: aut_proj_bundles}.

\begin{proposition} \label{prop: canonical_lifting_split_projective}
    Let $\sigma_{\mathbb{P}}$ be an automorphism of the pair $(\mathbb{P}(M\oplus N), C+D)$. 
    Then
    \begin{enumerate}
        \item there exists $\sigma_E \in \Aut(E)$ such that $\sigma_{\mathbb{P}} \in \Aut_{E, \sigma_E}(\mathbb{P}(M \oplus N))$;
        \item there exists a lifting $\widetilde{\sigma}_{\mathbb{P}}$ of $\sigma_P$ to the canonical lift such that $\widetilde{\sigma}_{\mathbb{P}} \in \Aut_{\mathbf{E}, \sigma_\mathbf{E}}(\mathbb{P}_{\mathbf{E}}(\mathcal{M} \oplus \mathcal{N}))$, where $\sigma_{\mathbf{E}}$ is the canonical lift of $\sigma_E$.
    \end{enumerate}
\end{proposition}

\begin{proof}
    As $E$ is an elliptic curve and the fibres of $\pi$ are rational curves, the automorphism $\sigma$ preserves the fibration. Thus there exists an automorphism $\sigma_E$ of $E$ such that 
    $\sigma_E \circ \pi =\pi \circ \sigma_{\mathbb{P}}$, proving (a).
    
    For (b), by \autoref{prop:relative_autom_proj_bundles} we know that the automorphism $\sigma_{\mathbb{P}}$ is equivalent to the data of a pair $(L, \sigma_{M \oplus N})$ where $L \in \Pic(E)$ and $[\sigma_{M \oplus N}] \in \Isom_E(\sigma_E^*(M \oplus N), (M \oplus N) \otimes L))/k^*$.
    The fact that $\sigma_{\mathbb{P}}$ is an automorphism of the pair implies that 
    \begin{enumerate}
        \item $\sigma(C)=C$: in this case, $[\sigma_{M \oplus N}]$ can be represented by $\sigma_{M} \oplus \sigma_N$, where $\sigma_M \colon \sigma_E^*M \to M \otimes L, \sigma_N \colon \sigma_E^*N \to N \otimes L$ are isomorphisms.
        \item $\sigma(C)=D$: in this case, $[\sigma_{M \oplus N}]$ can be represented by $\sigma'_M \oplus \sigma'_N$, where $\sigma'_M \colon \sigma_E^*M \to N$ and $\sigma'_N \colon \sigma_E^*N \to M$ are isomorphisms.  
    \end{enumerate}
    We suppose case (a) holds as the proof for (b) is analogous. 
    Consider the canonical lifting $\mathcal{L}$ of $L$ and note that the isomorphism $\sigma_{M}, \sigma_{N}$ can be lifted to isomorphism $\sigma_{\mathcal{M}} \colon \sigma_{\mathbf{E}}^* \mathcal{M} \to \mathcal{M}$ and  $\sigma_{\mathcal{N}} \colon \sigma_{\mathbf{E}}^* \mathcal{N} \to \mathcal{N}$.
   We can thus lift the pair $(L, [\sigma_{M \oplus N}])$ to the pair $(\mathcal{L}, [\sigma_{\mathcal{M}} \oplus \sigma_{\mathcal{N}}])$ which, by construction, is an isomorphism of the pair $(\mathbb{P}_{\mathbf{E}}(\mathcal{M} \oplus \mathcal{N}), \mathbf{C}+\mathbf{D})$.
\end{proof}

\subsection{Canonical liftings of $(\bP^1,\frac{1}{2}\sum_{i=1}^4q_i)$}\label{section:can_lift_P1_4pts}
We consider more closely the log Calabi--Yau structure on $\bP^1$ consisting of four points obtained in \autoref{corollary:possible_CY_log_curves}, in case it is globally sharply $F$-split. It turns out that such pairs have lifting properties that are similar to ordinary elliptic curves. We require the following lemma.

\begin{lemma}\label{lemma:GFS_and_crepant_maps}
Let $\pi\colon Y\to X$ be a finite surjective morphism of normal varieties. Suppose that $(X,\Delta)$ is a globally sharply $F$-split pair, and that we can write $K_Y+\Delta_Y=\pi^*(K_X+\Delta)$ where $\Delta_Y$ is effective. Then $(Y,\Delta_Y)$ is globally sharply $F$-split.
\end{lemma}
\begin{proof}
By naturality of the trace maps we have a commutative diagram
        $$\begin{tikzcd}
        H^0(X,\sO_X(\lfloor{(1-p^e)(K_X+\Delta)\rfloor}))\arrow[rr, "\Tr_{(X,\Delta)}"]\arrow[d, "\pi^*"] &&
        H^0(X,\sO_X)\arrow[d, hook] \\
        H^0(Y,\sO_Y(\lfloor{(1-p^e)(K_Y+\Delta_Y)\rfloor}))\arrow[rr, "\Tr_{(Y,\Delta_Y)}"] && H^0(Y,\sO_Y).
        \end{tikzcd}$$
By assumption, for $e$ large enough there exists $\xi\in H^0(X,\sO_X(\lfloor{(1-p^e)(K_X+\Delta)\rfloor}))$ such that $\Tr_{(X,\Delta)}(\xi)=1$. Then $\Tr_{(Y,\Delta_Y)}(\pi^*\xi)=1$, which implies that $(Y,\Delta_Y)$ is globally sharply $F$-split.
\end{proof}

\begin{proposition}\label{prop:canonical_lift_P1_4_pts}
Let $(\bP^1_k,\frac{1}{2}\sum_{i=1}^4q_i)$ be a globally sharply $F$-split log CY pair. Then $p>2$ and there exists a cartesian diagram 
        \begin{equation}\label{eqn:canonical_lift_P1_4_pts}
        \begin{tikzcd}
        E\arrow[d, "f"]\arrow[r, hook] & \mathbf{E}\arrow[d, "\mathbf{f}"] \\
        \bP^1_k\arrow[r, hook] & \bP^1_{W(k)}
        \end{tikzcd}
        \end{equation}
where:
    \begin{enumerate}
        \item E is an ordinary curve of genus $1$ and $\mathbf{E}$ is its canonical lift over $W(k)$;
        \item $f$ is a degree $2$ morphism whose branch locus is $\{q_1,\dots,q_4\}$;
        \item $\mathbf{f}$ is a degree $2$ morphism whose branch locus is four $W(k)$-rational points $\mathbf{q}_1,\dots,\mathbf{q}_4$ which respectively lift $q_1\dots,q_4$, and
        \item if $\iota\in \Aut_k(E)$ is the non-trivial deck transformation induced by $f$, and if $\iota_{W(k)}\in \Aut_{W(k)}(\mathbf{E})$ is its canonical lift, then $\mathbf{f}$ is the geometric quotient by the $\mu_{2,W(k)}$-action generated by $\iota_{W(k)}$.
    \end{enumerate}
\end{proposition}
\begin{proof}
The characteristic $p$ is greater than $2$ by \autoref{rmk:when_is_log_curve_GFS}.
Let $f\colon E\to \bP^1_k$ be the double cover of $\bP^1_k$ branched over the four points $q_1,\dots,q_4$. If $x$ is a coordinate on $\bP^1$ such that $q_i=[a_i:1]$, we can describe $E$ as the normalization of $\bP^1$ in the field $k(E)=k(x)[y]/(y^2-\prod_{i=1}^4(x-a_i))$. The morphism $f$ is the quotient by an involution $\iota$ of $E$ which is given on $k(E)$ by $(x,y)\mapsto (x,-y)$. 

By the Riemann--Hurwitz formula we see that $E$ is a genus one normal curve and that $K_E=f^*(K_{\bP^1}+\sum_{i=1}^4 \frac{1}{2}q_i)$. By \autoref{lemma:GFS_and_crepant_maps} it follows that $E$ is globally sharply $F$-split, or equivalently an ordinary genus one curve. So we are in position to apply \cite[Appendix, Theorem 1]{MS87} and get the canonical lifting $\mathbf{E}$ of $E$ over $W(k)$, equipped with a lifting of the absolute Frobenius. Moreover any endomorphism of $E$ admits a unique lift over $W(k)$ that commutes with the lift of Frobenius. Thus we obtain an involution $\iota_{W(k)}$ of $\mathbf{E}$ that lifts $\iota$.

Let $\mathbf{f}\colon \mathbf{E}\to \mathbf{P}=\mathbf{E}/\langle \iota_{W(k)}\rangle$ be the geometric quotient by the $\mu_{2,W(k)}$-action generated by $\iota_{W(k)}$. By \autoref{prop:restriction_linearly_red_qt} we see that 
    $$\mathbf{P}\otimes_{W(k)} k\cong (\mathbf{E}\otimes k)/\langle \iota_{W(k)}\otimes k\rangle
        \cong E/\langle \iota\rangle 
        \cong \bP^1_k
        \quad \text{and} \quad 
        \mathbf{f}\otimes k=f.$$
As $\bP^1_k$ is rigid, we obtain $\mathbf{P}\cong \bP^1_{W(k)}$.

The morphism $\mathbf{f}$ is separable of degree $2$. By purity of the branch locus \cite[0BMB]{stacks-project}, its ramification locus $D\subset \mathbf{E}$ has pure codimension one. Since $\mathbf{f}\otimes k=f$ is not ramified everywhere, we see that every component of $D$ dominates $W(k)$. If $K=\Frac(W(k))$, a Riemann--Hurwitz computation shows that $D_K$ has degree $4$ over $K$. Since $f=\mathbf{f}\otimes k$ has four ramification $k$-points (namely the points $f^{-1}(q_i)$ for $i=1,\dots,4$), we see that $D_K$ is the sum of four $K$-points $R_1,\dots,R_4\in \mathbf{E}$ which specialize to the ramification points of $f$. The images of these four points on $\bP^1_{W(k)}\otimes K$ define four $W(k)$-rational points $\mathbf{q}_1,\dots,\mathbf{q}_4\in \bP^1_{W(k)}$, which specialize (up to relabeling) to $q_1,\dots,q_4\in \bP^1_k(k)$.
\end{proof}

The construction of the square \autoref{eqn:canonical_lift_P1_4_pts} is functorial with respect to log isomorphisms in the following sense.

\begin{proposition}\label{prop:functoriality_can_lift_P1_4pts}
Let $\tau\colon (\bP^1_k,\frac{1}{2}\sum_{i=1}^4q_i)\cong (\bP^1_k,\frac{1}{2}\sum_{i=1}^4q_i')$ be a log isomorphism of globally sharply $F$-split log CY pairs. Then there is a commutative diagram
        $$\begin{tikzcd}
        &&& E'\arrow[dd, "f'" left] && 
        \mathbf{E}'
        \arrow[from=dlll, crossing over, "\Phi_{W(k)}" below right] \arrow[from=ll, hook]\arrow[dd, "\mathbf{f}'"]\\
        E\arrow[dd, "f"] \arrow[rr, hook] \arrow[urrr, "\Phi"] 
        && \mathbf{E} \\
        &&& \bP^1_k\arrow[rr, hook] && \bP^1_{W(k)} \\
        \bP^1_k\arrow[rr, hook]\arrow[urrr, "\tau"] && 
        \bP^1_{W(k)}
        \arrow[from=uu, crossing over,  "\mathbf{f}"] 
        \arrow[urrr, "\tau_{W(k)}" below right]
        \end{tikzcd}$$
where:
    \begin{enumerate}\setcounter{enumi}{2}
        \item $f,\mathbf{f}, \mathbf{q}_i$ and $f',\mathbf{f}',\mathbf{q}_i'$ are constructed as in \autoref{prop:canonical_lift_P1_4_pts};
        \item $\Phi\colon E\cong E'$ is an isomorphism and $\Phi_{W(k)}\colon \mathbf{E}\cong \mathbf{E}'$ is its canonical lift in the sense of \cite[Appendix]{MS87}, and
        \item $\tau_{W(k)}\colon (\bP^1_{W(k)},\frac{1}{2}\sum_{i=1}^4 \mathbf{q}_i)\cong (\bP^1_{W(k)},\frac{1}{2}\sum_{i=1}^4\mathbf{q}_i')$ lifts $\tau$.
    \end{enumerate}
\end{proposition}
\begin{proof}
We use the notations of the proof of \autoref{prop:canonical_lift_P1_4_pts}, adding a prime if we perform them for the log pair $(\bP^1_k,\frac{1}{2}\sum_{i=1}^4q_i')$. For example, $f'\colon E'\to \bP^1$ is the double cover of $\bP^1$ branched over the four points $q_1',\dots,q_4'$. If $x'$ is a coordinate on $\bP^1$ such that $q_i'=[a_i':1]$, then $E'$ is the normalization of $\bP^1$ in the field $k(E)=k(x')[y']/(y'^2-\prod_{i=1}^4(x'-a'_i))$ and $f'$ is the quotient by an involution $\iota'$ of $E'$ given by $(x',y')\mapsto (x',-y')$. 

We begin by constructing $\Phi\colon E\cong E'$. Since $\tau$ sends the set $\{q_1,\dots,q_4\}$ to $\{q_1',\dots,q_4'\}$ we must have 
    $$\tau^*\left(\prod_{i=1}^4(x'-a'_i)\right)=\lambda \prod_{i=1}^4(x-a_i)\quad
    \text{for some }\lambda\in k^\times .$$ 
Choose a square root $\lambda^{1/2}$ of $\lambda$ and define a $k$-linear morphism $\Phi^*\colon k(E')\to k(E)$ by 
        $$\Phi^*(x')=\tau^*(x'), \quad 
        \Phi^*(y')= \lambda^{1/2}y.
        $$
Then $\Phi^*$ induces an isomorphism $\Phi\colon E\to E'$ such that 
    \begin{equation}\label{eqn:commutativity_of_involutions}
        \iota'\circ\Phi=\Phi\circ \iota.
    \end{equation} 
This implies that the diagram
    \begin{equation}\label{eqn:commutativity_of_involutions_II}
    \begin{tikzcd}
    E\arrow[r, "\Phi"]\arrow[d, "f"]& 
    E'\arrow[d, "f'"]\\
    \bP^1_k \arrow[r, "\tau"] & \bP^1_k
    \end{tikzcd} 
    \end{equation}
is commutative.

Next we lift everything over $W(k)$. By \cite[Appendix, Theorem 1]{MS87}, any morphism with source and target amongst $E$ and $E'$ admits a lift to the canonical lifts $\mathbf{E}$ and $\mathbf{E}'$ over $W(k)$, which is unique with the property that it commutes with the lift of Frobenii. Let $\iota_{W(k)}\colon \mathbf{E}\cong \mathbf{E}$, $\iota'_{W(k)}\colon \mathbf{E}'\cong \mathbf{E}'$ and $\Phi_{W(k)}\colon \mathbf{E}\cong \mathbf{E}'$ be the canonical lifts of the involutions $\iota$ and $\iota'$ and of $\Phi$. It easily follows from \autoref{eqn:commutativity_of_involutions} and the compatibility with the lifts of Frobenii that
        \begin{equation} \label{eqn:commutativity_of_involutions_III}
        \iota'_{W(k)}\circ\Phi_{W(k)}=\Phi_{W(k)}\circ \iota_{W(k)}.
        \end{equation}
As $\mathbf{f}$ and $\mathbf{f}'$ are the geometric quotients by $\langle \iota_{W(k)}\rangle$ and $\langle \iota_{W(k)}'\rangle$ respectively, this equality implies that $\Phi_{W(k)}$ descends to an isomorphism $\tau_{W(k)}\colon \bP^1_{W(k)}\cong \bP^1_{W(k)}$ making the diagram
    $$\begin{tikzcd}
    \mathbf{E}\arrow[d, "\mathbf{f}"] \arrow[r, "\Phi_{W(k)}"] &
    \mathbf{E}' \arrow[d, "\mathbf{f}'"] \\
    \bP^1_{W(k)}\arrow[r, "\tau_{W(k)}"] & \bP^1_{W(k)}
    \end{tikzcd}$$
commutative, and whose specialization over $k$ is the diagram \autoref{eqn:commutativity_of_involutions_II}.

It remains to show that $\tau_{W(k)}$ upgrades to a log isomorphism between $(\bP^1_{W(k)},\frac{1}{2}\sum_{i=1}^4 \mathbf{q}_i)$ and $(\bP^1_{W(k)},\frac{1}{2}\sum_{i=1}^4 \mathbf{q}_i')$. It is equivalent to prove that $\Phi_{W(k)}$ sends the ramification divisor of $\mathbf{f}$ to the ramification divisor of $\mathbf{f}'$. This can be checked on the generic fibers; let $K=\Frac(W(k))$. Notice that the ramification points of $\mathbf{f}\otimes K$ (resp. of $\mathbf{f}'\otimes K$) are exactly the fixed points of $\iota_{W(k)}\otimes K$ (resp. of $\iota'_{W(k)}\otimes K$). If $R\in \mathbf{E}\otimes K$ is such a point, then using \autoref{eqn:commutativity_of_involutions_III} we find
        $$(\Phi_{W(k)}\otimes K)(R)=
        (\Phi_{W(k)}\otimes K)(\iota_{W(k)}\otimes K)(R)=
        (\iota'_{W(k)}\otimes K)(\Phi_{W(k)}\otimes K)(R).$$
So $(\Phi_{W(k)}\otimes K)(R)$ is a ramification point of $\mathbf{f}'\otimes K$, as desired. This completes the proof.
\end{proof}

Let us record the construction of \autoref{prop:canonical_lift_P1_4_pts} with the following definition.

\begin{definition}\label{def:canonical_lift_P1_4pts}
Suppose that $(\bP^1_k,\frac{1}{2}\sum_{i=1}^4q_i)$ is globally sharply $F$-split. Then its \emph{canonical lift over $W(k)$} is the pair $(\bP^1_k,\frac{1}{2}\sum_{i=1}^4\mathbf{q}_i)$ constructed in \autoref{prop:canonical_lift_P1_4_pts}.
\end{definition}

\begin{remark}\label{rmk:can_cyclic_cover}
Consider the diagram \autoref{eqn:canonical_lift_P1_4_pts}. Then it follows from \autoref{prop:quotient_is_cyclic_cover} that the canonical lift $\mathbf{E}$ of $E$ is the double cover of $\bP^1_{W(k)}$ branched along the four $W(k)$-rational points $\mathbf{q}_1,\dots,\mathbf{q}_4$, and that the induced $\mu_{2,W(k)}$-action on $\mathbf{E}$ is the canonical lift of the $\mu_{2,k}$-action on $E$ induced by $f$.
\end{remark}

%% file: 3_Log_liftings.tex
\section{Strong slc liftings: definitions}\label{section:def_log_lifts}

We start by recalling the definition of liftability from \cite[Definitions 2.6 and 2.8]{BBKW24} and by extending it in order to accommodate finite group actions. 

\begin{definition}[Lifting over $W(k)$]\label{def:lifting}
Let $f\colon (Y,D_Y)\to (X,D=\sum_{i=1}^r D_i)$ be a proper birational morphism of normal couples over $k$ with reduced boundaries, where $D_Y=f^{-1}_*D+E$ where $E$ is $f$-exceptional.
    \begin{enumerate}
        \item A \emph{lifting} of $(X,D)$ over $W(k)$ consists of a normal pair $(\mathcal{X},\mathcal{D}=\sum_{i=1}^r\mathcal{D}_i)$ where $\mathcal{X}$ and each $\mathcal{D}_i$ are flat, closed and separated over $W(k)$, together with an isomorphism $ \mathcal{X}\otimes k\cong X$ restricting to $\sD_i\otimes k\cong D_i$.
        \item A \emph{lifting} of $f$ over $W(k)$ is a proper birational morphism $\widetilde{f}\colon (\sY,\sD_\sY)\to (\sX,\sD)$ where:
            \begin{itemize}
                \item $(\sY,\sD_\sY)$ and $(\sX,\sD)$ are liftings of respectively $(Y,D_Y)$ and $(X,D)$ over $W(k)$;
                \item $\widetilde{f}_*\sO_\sY=\sO_\sX$ and
                    $$\Big(\widetilde{f}\colon (\sY,\sD_\sY)\to (\sX,\sD)\Big)\otimes k \ \cong  \ 
                    \Big(f\colon (Y,D_Y)\to (X,D)\Big).$$
            \end{itemize} 
    \end{enumerate}
\end{definition}

We emphasize that if $(\mathcal{X},\mathcal{D})$ is a lifting of $(X,D)$, then specialization over $k$ induces a bijection between the irreducible components of $\mathcal{D}$ and those of $D$. 
We demand $f$ to be closed as $(X, D)$ is not proper in the generality of this subsection (while for the statement of the main theorem we always assume the pair to be projective).

\begin{definition}\label{def:group_action_pair}
Let $(X,D)$ be a couple with reduced boundary over $k$, $(\sX,\sD)$ be a lift over $W(k)$ and $G$ be a finite abstract group. Then an action of $G$ on $(X,D)$ (resp.\ on $(\sX,\sD)$) is a group morphism $\rho\colon G\to \Aut_k(X,D)$ (resp.\ $\rho\colon G\to \Aut_{W(k)}(\sX)$) where the latter sub-groups are defined in \autoref{def:action_on_pairs}.
\end{definition}

\begin{definition}[Equivariant lifting over $W(k)$]\label{def:equiv_lifting}
Let $G$ be a finite abstract group, and let $f\colon (Y,D_Y)\to (X,D)$ be as in \autoref{def:lifting}. Suppose that $G$ acts on $(Y,D_Y)$ and $(X,D)$ and that $f$ is $G$-equivariant. Then a \emph{$G$-equivariant lifting of $f$} over $W(k)$ is a lifting $\widetilde{f}\colon (\sY,\sD_\sY)\to (\sX,\sD)$ of $f$ over $W(k)$ such that:
    \begin{enumerate}
        \item $G$ acts on $(\sY,\sD_\sY)$ and $(\sX,\sD)$,
        \item $\widetilde{f}$ is $G$-equivariant, and 
        \item the isomorphism
            $$\Big(\widetilde{f}\colon (\sY,\sD_\sY)\to (\sX,\sD)\Big)\otimes k \ \cong  \ 
                    \Big(f\colon (Y,D_Y)\to (X,D)\Big)$$
            is $G$-equivariant.
    \end{enumerate}
\end{definition}

Next, we recall the notion of strong log liftability of surfaces from \cite[Definition 2.11]{BBKW24}:

\begin{definition}[Strong log lifting I]\label{def:log_lifting_normal_case}
    Let $(X, D)$ be a normal surface couple over $k$, where $D$ is a reduced Weil divisor. 
    Let $f \colon Y \to (X, D)$ be a log resolution.
    We say that $(X,D)$ is \emph{strongly log liftable} over $W(k)$ if there is a commutative diagram
  \begin{equation*}\label{eqn:}
    \begin{tikzcd}
    (Y, f_*^{-1}D+E) \arrow[r ]\arrow[d, "f"] & 
    (\mathcal{Y}, f_*^{-1}\mathcal{D}+\mathcal{E})\arrow[d, "\widetilde{f}"] \\
    (X, D) \arrow[r] \arrow[d] & (\mathcal{X}, \mathcal{D}) \arrow[d] \\
    \Spec(k) \arrow[r] & \Spec(W(k)),
    \end{tikzcd} 
    \end{equation*}
    such that:
    \begin{enumerate}
        \item the exceptional locus of $f$ (resp.\ of $\widetilde{f}$) is the divisor $E$ (resp.\ $\mathcal{E}$), and
        \item $(\mathcal{Y}, f_*^{-1}\mathcal{D}+\mathcal{E})$ (resp.\ $(\mathcal{X}, \mathcal{D})$) is a lifting of $(Y, f_*^{-1}D+E)$ (resp.\ of $(X, D)$) in the sense of \autoref{def:lifting}.
    \end{enumerate}
    In this case, we say that 
    $(\mathcal{X}, \mathcal{D})$ (resp.\ $\widetilde{f}$) is a \emph{strong log lifting} of $(X,D)$ (resp.\ of $f$). 
\end{definition}

By \cite{BBKW24}*{Lemma 2.12}, the liftability of one log resolution is equivalent to the log liftability of any log resolution.
A strong log lifting should be thought of as an equisingular deformation to characteristic 0, as shown by \cite{BBKW24}*{Proposition 6.2} and by the following lemma (we recall, following \cite{ABP24}*{Definition 4.11}, that a family $f \colon (\mathcal{X}, \mathcal{D}) \to \Spec(W(k))$ is \emph{locally stable} if it is a family of pairs such that $(\mathcal{X}, \mathcal{D}+X_k)$ is slc).

\begin{lemma}\label{lemma:normalization_of_lifted_boundary}
Let $(X,D)$ be a reduced surface pair over $k$, and let $(\mathcal{X},\sD)$ be a strong log lifting over $W(k)$. Assume that $K_\sX+\sD$ is $\bQ$-Cartier. Then:
    \begin{enumerate}
        \item $\sD\otimes k\cong D$ extends to $(\sD^\nu, \Diff_{\sD^\nu}(0))\otimes k\cong (D^\nu, \Diff_{D^\nu}(0))$;
        \item the support of $\Diff_{\sD^\nu}(0)$ is \'{e}tale over $W(k)$;
        \item if $(X,D)$ is log canonical, then $(\sX,\sD)\to \Spec W(k)$ is a locally stable family.
    \end{enumerate}
\end{lemma}
\begin{proof}
By \cite[Proposition 4.5.(4)]{kk-singbook}, as $X$ is Cartier in $\sX$ and $\sD|_X=D$, we have that
        $$(K_\sX+X+\sD)|_X=K_X+D$$
is $\bQ$-Cartier. The different divisors on $\sD^\nu$ and $D^\nu$ are defined by adjunction.

For the rest of the proof, let us use the notations of \autoref{def:log_lifting_normal_case}. By definition, after possibly passing to a higher model, $(Y,f^{-1}_*D+E)$ is snc and $f^{-1}_*D$ is regular. It follows from \cite[Lemma 2.9]{BBKW24} that $(\mathcal{Y},f_*^{-1}\sD+\sE)$ is relatively snc and that $f^{-1}_*\sD$ is regular. Let us write
    $$K_Y+f_*^{-1}D+F=f^*(K_X+D),\quad 
    K_{\mathcal{Y}}+
    \widetilde{f}^{-1}_*\sD+\sF= 
    \widetilde{f}^*(K_{\mathcal{X}}+\sD)$$
where $\Supp(F)$ and $\Supp(\sF)$ are respectively contained in $E$ and $\sE$. Observe that $D^\nu=f_*^{-1}D$ and that under this identification it holds that $\Diff_{D^\nu}(0)=\Diff_{f_*^{-1}D}(F)$. A similar descriptions also holds over $W(k)$, namely
    $$(\sD^\nu,\ \Diff_{\sD^\nu}(0))=
    \left(\widetilde{f}^{-1}_*\sD,\ \Diff_{\widetilde{f}^{-1}_*\sD}(\sF)\right).$$
To prove the first statement, it remains to notice that the relative snc conditions implies
    $$\left(\widetilde{f}^{-1}_*\sD,\  \Diff_{\widetilde{f}^{-1}_*\sD}(\sF)\right)\otimes k \cong 
    \left(f_*^{-1}D, \ \Diff_{f^{-1}_*(D)}(F)\right).$$
These identifications also show that to prove the second statement, we have to guarantee that distinct irreducible components of the intersection of $\sF$ and $\widetilde{f}^{-1}_*\sD$ do not collapse together after specialization over $k$. As $\sF$ is supported on $\sE$, such a collapsing behaviour cannot happen because $(\mathcal{Y},\widetilde{f}^{-1}_*\sD+\sE)$ is relatively snc.

To prove the third statement, observe that since $\widetilde{f}$ is a lifting of $f$, every component of $\Exc(\widetilde{f})$ dominates $W(k)$. Therefore we actually have
    $$ K_{\mathcal{Y}}+Y+
    \widetilde{f}^{-1}_*\sD+\sF= 
    \widetilde{f}^*(K_{\mathcal{X}}+X+\sD).$$
Since
        $$\big(\widetilde{f}^* (K_\sX+X+\sD)\big)|_{Y}=
        f^*\big((K_\sX+X+\sD)|_X\big)
        =f^*(K_X+D)$$
we deduce that $\sF\otimes k=F$ as $\bQ$-divisors. We can now deduce the second statement using \cite{kk-singbook}*{Corollary 2.13}. 
\end{proof}

We generalise the notion of strong log lifting to demi-normal surfaces.

\begin{definition}[Strong log lifting II]\label{def:log_lifting_slc_case}
    Let $X$ be a demi-normal surface over $k$.
    Let $X^{\nu}$ be the normalisation of $X$ with the conductor scheme $D \subset X^{\nu}$.
    We say that $X$ is \emph{strongly log liftable} if there is a commutative diagram
    \begin{equation*}\label{eqn:}
    \begin{tikzcd}
    (X^{\nu}, D) \arrow[r ]\arrow[d, "\nu"] & 
    (\mathcal{X}^{\nu}, \mathcal{D})\arrow[d, "\widetilde{\nu}"] \\
    X \arrow[r]\arrow[d]  & \mathcal{X}\arrow[d] \\
    \Spec(k)\arrow[r] & \Spec(W(k)),
    \end{tikzcd} 
    \end{equation*}
    such that:
        \begin{enumerate}
            \item $(\mathcal{X}^{\nu}, \sD)$ is a strong log lifting of $(X^{\nu}D)$, and
            \item $\sX$ is demi-normal, flat separated closed and of finite type over $W(k)$, with $\sX\otimes k\cong X$, and normalization $\widetilde{\nu}\colon (\sX^\nu,\sD)\to \sX$.
        \end{enumerate}
    We say that $\mathcal{X}$ is a \emph{strong log lifting} of $X$. 

    We say that $X$ is \emph{strongly semi-log canonically liftable} if $X$ has semi-log canonical singularities and there exists a strong log lifting $\mathcal{X}$ of $X$ such that $K_{\mathcal{X}}$ is $\mathbb{Q}$-Cartier. In this case $\sX$ is a \emph{strong slc lifting} of $X$.
\end{definition}

\begin{remark}
We emphasize that 
    $$
    \text{strong slc lifting} \ =
    \ \text{strong log lifting} \ + \ 
    \begin{matrix}
        \bQ\text{-Cartier log canonical} \\
        \text{divisor on the total space}.
    \end{matrix}
    $$
We will only work with the more restrictive notion of strong slc liftings, as it allows us to use \autoref{lem: suff-cond-log-lift-nonnormal} below. Note that if $(X^\nu,D)$ is globally $F$-split and CY, then for any strong log lifting $(\sX^\nu,\sD)$ the log canonical divisor $K_{\sX^\nu}+\sD$ is necessarily $\bQ$-linearly trivial, by the proof of \cite[Theorem 6.8]{BBKW24}.
\end{remark}

We have the analogue of \autoref{lemma:normalization_of_lifted_boundary}:

\begin{lemma}
Notations as in \autoref{def:log_lifting_slc_case}. If $X$ is slc and $\mathcal{X}$ is a strong slc lifting, then $\sX\to \Spec W(k)$ is locally stable.
\end{lemma}
\begin{proof}
By definition $K_\sX$ is $\bQ$-Cartier. Since the fibers of $\sX\to \Spec(W(k))$ are demi-normal, we only have to show that the discrepancies of the pair $(\sX,X)$ are greater or equal to $-1$. We have $\widetilde{\nu}^*(K_\sX)=K_{\sX^\nu}+\sD$ and thus 
        $$\widetilde{\nu}^*(K_\sX+X)=K_{\sX^\nu}+X^\nu+\sD.$$
So we reduce to show that $(\sX^\nu, X^\nu+\sD)$ is log canonical. 
Since $X$ is slc we deduce that $(X^\nu,D)$ is log canonical, and so $(\sX^\nu,X^\nu+\sD)$ is indeed log canonical by \autoref{lemma:normalization_of_lifted_boundary}. 
\end{proof}


The gluing theory for slc varieties is recalled in \autoref{section:gluing_theory}. We show that, when $p>2$, a strong slc lifting is equivalent to a strong log lift of the normalization pair with $\bQ$-Cartier log canonical divisor, together with a lift of the gluing involution. The case $p=2$ will be handled in \autoref{section:p=2}.

\begin{proposition} \label{lem: suff-cond-log-lift-nonnormal}
    Let $X$ be a quasi-projective slc surface, with normalisation triple $(X^{\nu}, D, \tau)$. 
    Assume that $p>2$.
    Then $X$ is strongly slc liftable if and only if there is a quasi-projective strong lifting $(\mathcal{X}^{\nu}, \mathcal{D})$ of $(X^\nu,D)$ such that:
        \begin{enumerate}
            \item $K_{\sX^\nu}+\sD$ is $\bQ$-Cartier, and
            \item there is an involution $\widetilde{\tau}$ of $(\sD^\nu,\Diff_{\sD^\nu}(0))$ lifting $\tau$. 
        \end{enumerate}
\end{proposition}
\begin{proof}
Recall that $(\sD^\nu,\Diff_{\sD^\nu}(0))\otimes k =(D^\nu,\Diff_{D^\nu}(0))$ by \autoref{lemma:normalization_of_lifted_boundary} so the hypothesis about $\widetilde{\tau}$ makes sense. 

\medskip
\textsc{Proof of the ``only if" direction.} Assume that a strong slc lifting $\sX$ of $X$ exists. Then $\sX$ is demi-normal. So if $(\sX^\nu,\sD)$ is the normalization of $\sX$, there is a log involution $\widetilde{\tau}$ of $(\sD^\nu,\Diff_{\sD^\nu}(0))$ such that $\sX=\sX^\nu/R(\widetilde{\tau})$. By definition $(\sX^\nu,\sD)\otimes k= (X^\nu,D)$ and by \autoref{lemma:normalization_of_lifted_boundary} the pair $(\sX^\nu,X^\nu+\sD)$ is log canonical. Thus we may apply \cite[Theorem 2]{Posva_Gluing_stable_families_surfaces_mix_char} and find that $X=X^\nu/R(\widetilde{\tau}\otimes k)$.

It remains to show that $\tau=\widetilde{\tau}\otimes k$. This can be checked at the generic points of $D^\nu$, and so we may assume that $D$ and the conductor sub-scheme $C\subset X$ are normal. Then by \autoref{eqn:pushout_diagram} we have a push-out diagram
    $$\begin{tikzcd}
    D \arrow[r] \arrow[d] & X^\nu \arrow[d, "\nu"] \\
    D/\langle\tau\rangle \arrow[r] & X.       
    \end{tikzcd}$$
On the other hand, since $X=X^\nu/R(\widetilde{\tau}\otimes k)$ we also have a push-out diagram if we replace $D/\langle\tau\rangle$ by $D/\langle\widetilde{\tau}\otimes k\rangle$. Therefore $\tau$ and $\widetilde{\tau}\otimes k$ define the same invariant sub-sheaf in $\sO_{D}$, which implies by Galois theory that $\tau=\widetilde{\tau}\otimes k$.

\medskip
\textsc{Proof of the ``if" direction.} Assume that we are given a strong log lifting $(\sX^\nu,\sD)$ of $(X^\nu,D)$ such that $K_{\sX^\nu}+\sD$ is $\bQ$-Cartier, together with $\widetilde{\tau}$ lifting $\tau$. To construct $\sX$, we adapt the proofs of \cite[Proposition 3.1.8 and 3.2.1]{Posva_Gluing_stable_families_surfaces_mix_char}. Let $R(\widetilde{\tau})\rightrightarrows \sX^\nu$ be the equivalence relation generated by $\widetilde{\tau}$. We claim that it is a finite equivalence relation. This can be checked on each fiber of $\sX^\nu$ over $W(k)$: on the generic one this is shown as in \cite[Proof of Theorem 4.1.1]{Posva_Gluing_for_surfaces_and_threefolds}, and on the special one as in \cite[Proof of Proposition 3.1.8]{Posva_Gluing_stable_families_surfaces_mix_char} (\footnote{
    The assumption of properness in the statement of \cite[Proposition 3.1.8]{Posva_Gluing_stable_families_surfaces_mix_char} and in \cite[\S 3.1]{Posva_Gluing_stable_families_surfaces_mix_char} is not needed to establish this finiteness result.
}). By \cite[Proposition 33]{Kollar_Quotients_finite_equiv_relations} the geometric quotient of the generic fiber by the restricted equivalence relation exists. By \cite[Theorem 1.4]{Wit20}, this implies that the geometric quotient $\sX=\sX^\nu/R(\widetilde{\tau})$ exists as an algebraic space of finite type over $W(k)$. Furthermore, $\sX$ is a scheme by the quasi-projectivity of $\sX^\nu$ and \cite[Corollary 48]{Kollar_Quotients_finite_equiv_relations}.

We claim that $\sX$ is a strong slc lifting of $X$. By \cite[Proposition 3.4.1]{Posva_Gluing_for_surfaces_and_threefolds} (whose proof also works for $W(k)$-schemes), the scheme $\sX$ is demi-normal, the quotient morphism $\sX^\nu\to \sX$ is the normalization and $\sD$ is its conductor. This implies that $\sX\to \Spec(W(k))$ is flat, separated and closed. As $(\sX^\nu, X^\nu+\sD)$ is lc by \autoref{lemma:normalization_of_lifted_boundary}, it follows from \cite[Proposition 3.2.1]{Posva_Gluing_stable_families_surfaces_mix_char} that $K_\sX$ is $\bQ$-Cartier, and therefore $\sX$ is locally stable over $W(k)$. It remains to check that $\sX\otimes k=X$: as $\sD\otimes k=D$ is reduced, this follows from \cite[Theorem 2]{Posva_Gluing_stable_families_surfaces_mix_char}.
\end{proof}

\subsection{Demi-normality deforms}
It turns out that the requirements on $\sX$ in \autoref{def:log_lifting_slc_case} can be weakened, as demi-normality deforms in integral $S_2$-families. Since this is a statement of independent interest that we could not locate in the literature, we record it here.

\begin{lemma}\label{lemma:semi_norm_is_S2}
Let $T$ be a reduced excellent $S_2+G_1$ scheme. Then its semi-normalization $T^{\sn}$ is also $S_2$.
\end{lemma}

Here, $T$ being $S_2$ means that $\sO_T$ satisfies the corresponding Serre condition; $T$ being $G_1$ means that if $\eta\in T$ is a point of codimension $\leq 1$, then $\sO_{T,\eta}$ is a Gorenstein local ring; and we refer to \cite[\S I.7.2]{Kollar_Rational_curves} for the definition of the semi-normalization $T^\text{sn}$.

\begin{proof}
Since $T$ is excellent we regard $\sO_{T^{\sn}}$ as a coherent $\sO_T$-module, and it suffices to show that it is $S_2$ as such. Recall that a section $t\in \sO_{T^\nu}(V)$ belongs to $\sO_{T^{\sn}}(V)$ if and only if $t^2,t^3\in \sO_T(V)$ \cite[Proposition 2.10$^*$]{Hamann_On_the_R_invariance}. 

By \cite[Remark 1.11]{Hartshorne_Generalized_divisors_on_Gor_schemes}, a torsion-free coherent module $M$ on $T$ is $S_2$ if and only if for every open $U$ and closed $Y\subset U$ of codimension $\geq 2$, the restriction $M(U)\to M(U-Y)$ is bijective. Since semi-normalization commute with localization, we reduce to check this property for $M=\sO_{T^{\sn}}$ with $U=T$. Take a section $\overline{s}\in \sO_{T^{\sn}}(T\setminus Y)\subset \sO_{T^\nu}(T\setminus Y)$: as $T^\nu$ is $S_2$ it extends uniquely to a section $s\in \sO_{T^\nu}(T)$. We have $\overline{s}^2,\overline{s}^3\in \sO_T(T\setminus Y)$, and as $T$ is $S_2$ we obtain that $s^2,s^3\in \sO_T(T)$. Therefore $s\in \sO_{T^{\sn}}(T)$ as desired.
\end{proof}

\begin{proposition}
Let $(R,\mathfrak{m},\kappa)$ be an excellent DVR, and let 
$\sX$ be a reduced equidimensional $S_2$ scheme that is flat separated closed of finite type over $\Spec(R)$. Assume that $X=\sX\otimes k$ is demi-normal \emph{(\footnote{
    With this condition on $X$, the $S_2$ condition on $\sX$ is automatic if $\dim \sX=3$. For then its fiber $X$ is Cohen--Macaulay, and thus $\sX$ is also Cohen--Macaulay by \cite[Corollary to Theorem 23.3]{Matsumura_CRT} 
})}. Then $\sX$ is demi-normal.
\end{proposition}
\begin{proof}
First, we show that $\sX$ is Gorenstein in codimension one. Let $\sU$ be the Gorenstein locus of $\sX$: by \cite{Greco_Marinari_Openess_for_Gor_and_compl_int}, it is an open subset of $\sX$. Let $Z=\sX\setminus \sU$: it suffices to show that $\codim_{\sX}Z\geq 2$. Suppose that an irreducible component $Z'$ of $Z$ has codimension one. If $Z'$ does not dominate $\Spec(R)$, then $Z'$ is an irreducible component of $X$. But $X$ is reduced, in particular Gorenstein in codimension $0$: by \cite[0BJL]{stacks-project} this implies that $\sX$ is Gorenstein at the generic points of $X$, and we get a contradiction. Therefore, as $\sX$ is closed over $\Spec(R)$, we see that $Z'\to \Spec(R)$ must be surjective. By flatness this implies that $\codim_X(Z'\cap X)=1$. But $X$ is demi-normal, hence Gorenstein in codimension one, and as before we get a contradiction. Therefore $\codim_{\sX}Z\geq 2$ as desired.


Next we show that $\sX$ is semi-normal. Let $\sA=\sO_{\sX^{\sn}}$ be the semi-normalization of $\sO_{\sX}$: it is a finite $\sO_{\sX}$-algebra, flat over $\Spec(R)$, and the inclusion $\iota\colon \sO_\sX\to \sA$ is an isomorphism on a relatively dense open subset. By \autoref{lemma:semi_norm_is_S2} we see that $\sA$ is $S_2$. Consider $\iota\otimes \kappa\colon \sO_X=\sO_{\sX}\otimes \kappa\to \sA\otimes \kappa$. It is an isomorphism on a dense open subset of $X$. Moreover its domain and codomain are $S_1$ and have full support as $\sO_X$-modules, thus they are torsion-free by \cite[0AUV]{stacks-project}. This shows that $\iota\otimes \kappa$ is injective. So the morphism $\Spec_{X}(\sA\otimes \kappa)\longrightarrow X$ induced by $\iota\otimes \kappa$ factors the normalization morphism of $X$. Since it is obtained by restricting the semi-normalization morphism of $\sX$ over $\kappa$, we also see that it is a homeomorphism inducing isomorphisms on residue field extensions: therefore it factors the semi-normalization of $X$ (see \cite[Definition 2.3.1]{Arvidsson_Posva_Normality_min_lc_centers}). Since $X$ is semi-normal by \cite[Corollary 2.3.8]{Posva_Gluing_stable_families_surfaces_mix_char}, we obtain that $\iota\otimes \kappa$ is bijective. So if $Q$ is the cokernel of $\iota$, we have $Q\otimes \kappa=0$ and by Nakayama's lemma it follows that $Q=0$. In other words, $\sX$ is semi-normal 
(\footnote{
    The fact that semi-normality deforms is proved more generally in \cite{Heitmann_Lifting_seminormality}.
}).

Hence $\sX$ is $S_2+G_1$ and semi-normal. By \cite[Lemma 3.2.4]{Posva_Gluing_for_surfaces_and_threefolds} we deduce that $\sX$ has at worst nodal singularities in codimension one. Since $\sX$ is $S_2$ it follows, by definition, that $\sX$ is demi-normal.
\end{proof}

%% file: 4_Canonical_liftings.tex
\section{Strong slc liftings without $4A_1$-curves in the conductor}\label{section:lift_index_1}
Let $S_0$ be a projective slc globally $F$-split CY surface  with normalization $(S,D)$. In this section, we construct a strong slc lifting of $S_0$, 
except when there is a component in the boundary $D$ passing through four $A_1$ singularities (we call such a component a $4A_1$-curve, see \autoref{def:4A_1}). These special cases will be treated in \autoref{section:4A_1}. 

\subsection{Is there a genus one curve in the conductor?}
Let $(S,D)$ be a log canonical projective surface pair with $K_S +D \sim_{\mathbb{Q}} 0$.
As we will see in \autoref{prop:genus_1_forces_Z_linear_CY_cond}, if the conductor $D$ contains a (possibly reducible) curve of arithmetic genus $1$, then in most cases $K_S+D$ is $\mathbb{Z}$-linearly trivial. We start by introducing what the exceptions will be.

\begin{definition}[Surfaces of Enriques-type]\label{def:BLE-type}
Let $(S,D)$ be a normal projective surface pair with reduced non-empty boundary $D \neq 0$. 
We say that $(S,D)$ is \emph{of Enriques-type} if the following condition is satisfied:
 there exists a crepant birational map $\phi\colon (S,D)\dashrightarrow (Y,D_Y)$ such that: $D_Y=\phi_*D\neq 0$ is irreducible, $Y$ is a $\bP^1$-bundle over an elliptic curve $B$, the restriction $D_Y\to B$ is \'{e}tale of degree $2$, and $K_Y+D_Y\sim_\bQ 0$. 
\end{definition}

\begin{remark}
    We call such surfaces $(S,D)$ of Enriques-type because they appear as normalisation of irreducible components of certain degenerations of Enriques surfaces as explained in \cite{Mor81}.
\end{remark}

\begin{lemma}\label{lemma:BLE_is_canonical}
Suppose that $(S,D)$ is of Enriques-type. Then $(S,D)$ has canonical singularities and $D$ is irreducible.
\end{lemma}
\begin{proof}
Let $\phi\colon (S,D)\dashrightarrow(Y,D_Y)$ be the crepant birational model provided by \autoref{def:BLE-type}. By definition, for every divisor $E$ over $S$, we have $a(E; S, D) = a(E; Y, D_Y)$. Since $(Y,D_Y)$ is log smooth and $D_Y$ is irreducible, the unique divisor over $Y$ with negative discrepancy with respect to $(Y, D_Y)$ is $D_Y$. Assume that $E$ is exceptional over $S$: then the fact that $D_Y=\phi_*D$ implies that $E\neq D_Y$, and so $a(E; S, D) = a(E; Y, D_Y) \geq 0$. This shows that $(S,D)$ is canonical.
If $D$ is not irreducible, there exists an irreducible component $E'$ of $D$ that is contracted by $\phi$. So $E'$ is exceptional over $Y$, and it would follow that
        $$-1=a(E';S,D)=a(E';Y,D_Y)\geq 0,$$
which is impossible. So $D=\phi^{-1}_*D_Y$ is irreducible.
\end{proof}

In what follows we usually assume that $p>2$: the case $p=2$ will be treated separately in \autoref{section:p=2}.

\begin{proposition}\label{prop:boundary_disconnected_over_non_rational_curve}
Let $(T,\Delta)$ be an irreducible projective log canonical surface pair with $0\neq \Delta$ reduced and $K_T+\Delta\sim_{\mathbb{Q}} 0$. Assume that $p>2$, and that there exists a fibration $T\to B$ onto a normal proper curve of genus $1$ such that every irreducible component of $\Delta$ dominates $B$. Then either $\Delta$ is not connected, or $(T,\Delta)$ is of Enriques-type. 
\end{proposition}
\begin{proof}
The strategy of proof is as follows: we will go to a smooth minimal of $T$ over $B$, and observe that in most cases $\Delta $ must be disconnected. The remaining cases will be of Enriques-type.

Notice that any connected component of $\Delta$ is irreducible. For otherwise there are two irreducible components $F_1$ and $F_2$ of $\Delta$ that intersect each other. If $\Upsilon=\Delta-F_1-F_2$ then
        \begin{eqnarray*}
        0 &=& (K_T+\Upsilon +F_1+F_2)\cdot F_1 \\
        &=& \deg (K_{F_1^\nu}+\Diff_{F_1^\nu}(\Upsilon))+(F_2\cdot F_1) \\
        &\geq & 0+\deg \Diff_{F_1^\nu}(\Upsilon)+(F_2\cdot F_1) \\
        &>& 0
        \end{eqnarray*}
where we have used that $\deg K_{F_1^\nu} \geq 0$ since $F_1$ dominates $B$ which has genus 1. This is obviously a contradiction. Moreover, it follows from adjunction and \cite[Proposition 2.35]{kk-singbook} that the components of $\Delta$ are smooth.

Let $f\colon Y\to T$ be a minimal log resolution of $(T,\Delta)$. Then $Y$ is regular and we can write
        $$K_Y+\Delta_Y+\Gamma=f^*(K_T+\Delta)\sim_\mathbb{Q} 0$$
where $\Delta_Y$ is the strict transform of $\Delta$, and $\Gamma$ is effective and $f$-exceptional. By construction of the strict transform, the irreducible components of $\Delta_Y$ are in bijection with the irreducible components of $\Delta$. 

Now we run a $K_Y$-MMP over $B$ to obtain a birational $B$-contraction $h\colon Y\to Y'$. 
Since $K_Y$ is not pseudoeffective and $Y'$ has dimension 2, the structural morphism $g\colon Y'\to B$ must be a Mori fiber space. As $Y'$ is regular, we obtain that $Y'\to B$ is a $\mathbb{P}^1$-bundle. In addition, if $\Delta_{Y'}=h_*\Delta_Y$ and $\Gamma'=h_*\Gamma$ then $K_{Y'}+\Delta_{Y'}+\Gamma'\sim_\mathbb{Q} 0$. Notice that $\Gamma'$ is vertical over $B$, while every component of $\Delta_{Y'}$ dominates $B$. In particular, the irreducible components of $\Delta_{Y'}$ are disjoint regular curves of genus $1$, and are in bijection with the irreducible components of $\Delta_Y$.

We have $\NS(Y')=\mathbb{Z}[F]\oplus \mathbb{Z}[C_n]$ where $F$ is the class of a fiber of $g$ and $C_n$ is a section of $g$ with self-intersection $-n<0$. Write $\Delta_{Y'}\equiv aC_n+bF$. The number $a$ is determined by the intersection of $\Delta_{Y'}$ with a general fiber of $g$. The birational morphism $h\circ f^{-1} \colon T\dashrightarrow Y'$ is an isomorphism over a non-empty open subset of $B$, and as $K_T+\Delta_T\sim_{\mathbb{Q}} 0$ we deduce that $a=2$.

\begin{claim}
If $n>0$ then $\Delta$ is not connected.
\end{claim}
\begin{proof}\renewcommand{\qedsymbol}{$\lozenge$}
It suffices to show that $\Delta_{Y'}$ is not connected. In fact, we claim that $\Delta_{Y'}$ contains two disjoint sections of $g$. We adapt the argument of \cite[Proposition 5.9]{BBKW24} and distinguish two cases:
    \begin{enumerate}
        \item $C_n\not\subset \Supp(\Delta_{Y'})$. Then $0\leq \Delta_{Y'}\cdot C_n=-2n+b$. On the other hand, by adjunction we have
                \begin{eqnarray*}
                    0=\deg K_{\Delta_{Y'}} &=& (K_{Y'}+\Delta_{Y'})\cdot \Delta_{Y'} \\
                    &=& (K_{Y'}+\Delta_{Y'})\cdot 2C_n + (K_{Y'}+\Delta_{Y'})\cdot bF \\
                    &=& [(K_{Y'}+C_n)+(C_n+bF)]\cdot 2C_n + (-\Gamma')\cdot bF \\
                    & = & 2\deg K_{C_n}-2n+2b + 0.
                \end{eqnarray*}
        As $C_n\cong B$ we obtain that $n=b$. This contradicts $0\leq -2n+b$, so this case cannot happen.
        \item $C_n\subset \Supp(\Delta_{Y'})$. Then let us write $\Delta_{Y'}+\Gamma'=C_n+D$, where $D$ is an effective $\mathbb{Q}$-divisor. We have $(\Delta_{Y'}+\Gamma')\cdot F=(-K_{Y'})\cdot F=2$, and therefore $D\cdot F=1$ (which implies in particular that $D\neq 0$). Now we compute
                \begin{eqnarray*}
                    D\cdot C_n &= & (\Delta_{Y'}+\Gamma')\cdot C_n -C_n^2 \\
                    &=& -(K_{Y'}+C_n)\cdot C_n \\
                    &=& -\deg K_{C_n}
                \end{eqnarray*} 
        which is $0$ as $C_n\cong B$. Therefore $D$ is a section of $g$ which is disjoint from $C_n$. 
    \end{enumerate}
We have therefore obtained that $\Delta_{Y'}$ contains two disjoint sections of $g$, as claimed. 
\end{proof}

It remains to consider the case $n=0$.
By construction, $\Gamma'$ is supported on fibers of $Y'\to B$, so we have $\Gamma'\equiv cF$ with $c\geq 0$. Moreover, $\Delta_{Y'}\to B$ has degree $2$ and $\Delta_{Y'} \equiv 2C_0+bF$. Since $C_0^2=0$ and $K_{Y'}\equiv -2C_0$ \cite[III.2.9 and III.2.11]{Ha77}, we obtain
                $$0=(K_{Y'}+\Delta_{Y'}+\Gamma')\cdot C_0=b+c.$$
Thus $\Delta_{Y'}\cdot C_0=b=-c\leq 0$. Assume $\Delta_{Y'}=C_0+G$ with $G$ reduced and not containing $C_0$. Then $-c=C_0\cdot \Delta_{Y'}=C_0\cdot G\geq 0$, so in fact $c=b=0$ and $G$ is disjoint from $C_0$. Hence $\Delta_{Y'}$ is also disconnected in that case.

The only remaining case is $n=0$ and $C_0\not\subset \Supp(\Delta_{Y'})$: then necessarily $c=0$ and so $\Gamma'=0$.
If $\Delta_{Y'}$ is disconnected, then so is $\Delta$. If $\Delta_{Y'}$ is connected, then it must be a single regular curve of genus $1$, and $\Delta_{Y'}\to B$ is \'{e}tale of degree $2$. By construction $h\circ f^{-1}\colon (T,\Delta)\dashrightarrow (Y',\Delta_{Y'})$ is crepant with $(h\circ f^{-1})_*\Delta=\Delta_{Y'}$, so in that case $(T,\Delta)$ is of Enriques-type.
\end{proof} 

We are almost ready to prove the main result of this sub-section, \autoref{prop:genus_1_forces_Z_linear_CY_cond}. For its proof, the following notion is useful:

\begin{definition}
For a sub-pair $(Z,E)$ with $K_Z+E\sim_\bQ 0$, the \emph{torsion index} of $K_Z+E$ is the smallest positive integer $m>0$ such that $m(K_Z+E)\sim 0$. 
\end{definition}

\begin{lemma}\label{lemma:torsion_index_crepant_inv}
Let $(Z,E)\dashrightarrow (Z',E')$ be a crepant birational map of normal sub-pairs (of any dimension, in any characteristic). Then $K_Z+E\sim_\bQ 0$ if and only if $K_{Z'}+E'\sim_\bQ 0$, and in that case the torsion indices of $K_Z+E$ and $K_{Z'}+E'$ are equal.
\end{lemma}
\begin{proof}
By resolving the given crepant map we may assume that we have a crepant morphism $\phi\colon (Z,E)\to (Z',E')$, cf.\ \cite[2.23.2]{kk-singbook}. For any $m>0$, if $m(K_{Z'}+E')\sim 0$ then the pull-back of a nowhere vanishing section together with the crepant condition says that $m(K_{Z}+E)\sim 0$ as well. Hence $K_{Z'}+E'$ being $\bQ$-linearly trivial with torsion index $m$ implies that $K_Z+E$ is $\bQ$-linearly trivial with torsion index dividing $m$. 

Conversely, suppose that $n(K_Z+E)\sim 0$ for some $n>0$. Observe that $nE$ is necessarily a $\mathbb{Z}$-divisor. Since $E-\phi^{-1}_*E'$ is $\phi$-exceptional, we see that $nE'$ is also a $\mathbb{Z}$-divisor. Pick a big open subset $U\subseteq Z'$ where $n(K_{Z'}+E')$ is Cartier. If $\phi_U\colon \phi^{-1}(U)\to U$ is the induced map, then
        $$\phi_U^*\sO(n(K_{Z'}+E')|_U)
        =\sO(n(K_Z+E))|_{\phi^{-1}(U)}\sim 0.$$
Since $\phi_*\sO_{Z}=\sO_{Z'}$, by pushing-forward we obtain  
$\sO(n(K_{Z'}+E'))|_U\sim 0$. Since $U$ is big, it follows that $n(K_{Z'}+E')\sim 0$. So $K_{Z'}+E'\sim_\bQ 0$ with torsion index dividing $n$. This completes the proof.
\end{proof}

\begin{proposition}\label{prop:genus_1_forces_Z_linear_CY_cond}
Let $(S,D)$ be an irreducible projective globally sharply $F$-split log canonical surface pair over $k$ where $D$ is a reduced divisor and $K_S+D\sim_\mathbb{Q} 0$. Assume that $p>2$ and that there exists a connected reduced divisor $0<E\leq D$ with $K_E \sim 0$. Then either:
    \begin{enumerate}
        \item $K_S+D\sim 0$, or
        \item $(S,D)$ is of Enriques-type (cf.\ \autoref{def:BLE-type}).
    \end{enumerate}
\end{proposition}

\begin{proof}
We construct a crepant birational model of $(S,D)$ that witnesses one of the two conditions. Note that these conditions can be checked on a dlt modification: so we can assume that $(S, D)$ is dlt and globally sharply $F$-split.

Let us write $D=D'+E$ where $D'$ and $E$ have no irreducible component in common. Notice that $(S,D')$ is dlt. 
We run a $(K_S+D')$-MMP. Since $K_S+D'\equiv -E$ is not pseudoeffective, we obtain a birational morphism $\psi\colon S\to T$, and a Mori fibre space structure $h\colon (T,\psi_*D')\to B$ \cite{Tan18}. 
We write $E_T=\psi_*E$ and $D'_T=\psi_*D'$. The pair $(T,D'_T)$ is dlt, so in particular $T$ is klt by \autoref{lemma:everything_is_Q_Cartier}. 
Since $K_S+D'+E\sim_\bQ 0$, the morphism $h\colon (S,D'+E)\to (T,E_T+D'_T)$ is crepant. So $K_T+D'_T+E_T\sim_\bQ 0$, which together with the existence of the Mori fiber structure implies that $E_T\neq 0$. 

\begin{claim}\label{claim:image_genus_1_div}
The divisor $E_T$ belongs to the regular locus of $T$. 
\end{claim}

\begin{proof}\renewcommand{\qedsymbol}{$\lozenge$}
We can decompose $\psi$ as 
    $$S=S_0
    \overset{\psi_0}{\longrightarrow}
    S_1\overset{\psi_1}{\longrightarrow}
    \quad \dots \quad
    \overset{\psi_N}{\longrightarrow}
    S_{N+1}=
    T$$
where, upon letting $D'_{i+1}=\psi_{i,*}E_i$ and $E_{i+1}=\psi_{i,*}E_i$, each $\psi_i$ is the contraction of a negative $(K_{S_i}+D_i')$-ray. Notice that $K_{S_i}+D'_i+E_i\sim_\bQ 0$ at every step. Consider the following conditions:
    \begin{itemize}
        \item[$(\flat_i)$] $E_i$ is either a regular genus $1$ curve or a cycle of regular rational curves, and in the latter case $(S_i,E_i)$ is log smooth at the $0$-strata of $E_i$;
        \item[$(\sharp_i)$] $E_i$ belongs to the regular locus of $S_i$.
    \end{itemize}
Note that $(\flat_0)$ holds, as $E$ has only nodal singularities, $p_a(E)=1$, and $(S,D'+E)$ is dlt.
We will prove that $(\flat_i)$ implies $(\sharp_i)$, and that $(\flat_i)+(\sharp_i)$ implies $(\flat_{i+1})$. By induction, this will prove the claim. 

Assume that $(\flat_i)$ holds. By adjunction along $E_i^\nu$, we see that $E_i$ is disjoint from $D'_i$. Let $E_i=\sum_j E_{ij}$ be the decomposition into irreducible components. For every $j$ the divisor $E_i-E_{ij}$ is Cartier in a neighbourhood of $E_{ij}$. Adjunction along $E_{ij}$ entails
        $$0= (K_{S_i}+E_i)\cdot E_{ij}=\deg(K_{E_{ij}}+\Diff_{E_{ij}}(0))
        +(E_i-E_{ij})\cdot E_{ij}$$
and so we must have $\Diff_{E_{ij}}(0)=0$ for every $j$. By \cite[Proposition 2.35.(3)]{kk-singbook} we deduce that $(\sharp_i)$ holds.

Now assume that $(\flat_i)+(\sharp_i)$ holds. Let $\xi$ be the curve class that $\psi_i$ contracts. It satisfies $(K_{S_i}+D_i')\cdot \xi <0$, which is equivalent to $E_i\cdot \xi>0$. There are two cases: 
    \begin{enumerate}
        \item If an irreducible component $E_{ij}$ of $E_i$ belongs to $\mathbb{R}_+\xi$, then $E_{ij}^2$ is a negative integer, as it is contracted and lies on the regular locus. Since $E_i\cdot \xi>0$, this is only possible if $E_i$ is a cycle of rational curves. In that case
        $$0<E_i\cdot E_{ij}=2+E_{ij}^2$$
        forces $E_{ij}^2=-1$. In other words $E_{ij}$ is a $(-1)$-curve, so $(\flat_{i+1})$ holds. 

        \item No component of $E_i$ is contracted by $\psi_i$.
        Then $E_i\to E_{i+1}$ is finite and birational, so $E_{i+1}^\nu$ is the common normalization. Let us write
        $$\Diff_{E_{i+1}^\nu}(D_{i+1}')=N(E_{i+1})+B+(D_{i+1}’\cdot E_{i+1})$$
        where $N(E_{i+1})$ is the preimage of the nodes of $E_{i+1}$ and $B\geq 0$ is induced by the singularities of $S_{i+1}$. We have $K_{E_{i+1}^\nu}+N(E_i)\sim 0$. As $N(E_{i+1})\geq N(E_i)$, we see that $K_{E^\nu_{i+1}}+\Diff_{E^{i+1}_\nu}(D_{i+1}’)$ would have positive degree if $N(E_{i+1})>N(E_i)$ or if $B> 0$. This contradicts the fact that $K_{S_{i+1}}+D_{i+1}’+E_{i+1}$ is $\mathbb{Q}$-linearly trivial. Therefore $E_i\to E_{i+1}$ is an isomorphism and $(S_{i+1},E_{i+1})$ is log smooth in a neighborhood of $E_{i+1}$: this proves both $(\flat_{i+1})$ and $(\sharp_{i+1})$ at once.
    \end{enumerate}
The claim is proved.
\end{proof}

Since $E_T$ belongs to the regular locus of $T$, the adjunction exact sequence along $E_T$ induces the exact sequence of cohomology
	$$H^0(T,\sO(K_T+D’_T+E_T))\to H^0(E_T, \sO_{E_T}(K_{E_T})) \to H^1(T,\sO(K_T+D’_T)).$$
The proof of \autoref{claim:image_genus_1_div} above shows that condition $(\flat_{N+1})$ holds, and therefore $K_{E_T}\sim 0$. So the group in the middle is non-zero. We are going to prove that either the $H^1$ group vanishes, or that $D’_T=0$ and $(T,E_T)$ is of Enriques-type. In the first situation, we get $K_T+D’_T+E_T\sim 0$, and thus $K_S+D\sim 0$ by \autoref{lemma:torsion_index_crepant_inv}. In the second situation, by definition we get that $(S,D)$ is of Enriques-type. To obtain this dichotomy, we distinguish a few cases:

\begin{enumerate}
\item $\dim B=0$. In that case we always have $K_S+D\sim 0$. Indeed:
            \begin{enumerate}
                \item If $D_T'=0$, then we have $H^1(T,\sO_T(K_T))=0$, or dually $H^1(T,\sO_T)=0$, because $T$ is a rational surface.
                \item If $D_T'\neq 0$, then as $T$ has Picard rank one, the effective divisor $D_T'$ is ample. Moreover, by \cite[Lemma 3.5]{SS10} the surface $S$ is globally sharply $F$-split and thus $T$ is globally sharply $F$-split by \cite[Lemma 2.3]{BBKW24}. Hence we apply \cite[Proposition 2.3]{Ber21} to obtain $H^1(T,\mathcal{O}_{T}(-D_T'))=0$ and Serre duality implies $H^1(T,\mathcal{O}_{T'}(K_{T}+D_T'))=0$.
            \end{enumerate}
        
    \item $\dim B=1$. Here we distinguish again two cases:
            \begin{enumerate}
                \item If $D'_T=0$, then $K_T+E_T\sim_\mathbb{Q}0$. As $-K_T$ is $h$-ample, we see that $E_T$ is not contained in a fiber of $h$. 
                
                If $E$ is a cycle of rational curves, then $B$ is rational. As before we obtain that $T$ is rational and so $K_S+D\sim 0$.
                
                Assume $E$ is a regular curve of genus $1$, so $B$ has genus $0$ or $1$. If $B$ has genus $0$ then $T$ is rational and so $K_S+D\sim 0$ as before. If $B$ has genus $1$, then $E_T$ cannot be a cycle of rational curves,so being connected it is must be an irreducible curve of genus $1$. It follows from \autoref{prop:boundary_disconnected_over_non_rational_curve} that $(T,E_T)$ is of Enriques-type. 
                
                \item If $D'_T\neq 0$ we show that $H^1(T,\sO(K_T+D'_T))=0$. By the Leray spectral sequence along $h$ it is equivalent to show the two equalities
                $$H^0(B, R^1h_*\sO_T(K_T+D_T'))=0=H^1(B, h_*\sO_T(K_T+D_T')).$$
                The $H^1$ term is easy to deal with, for $h_*\sO(K_T+D_T')=0$. Indeed, if there was a non-zero section $\varphi\in \sO_T(K_T+D_T')(h^{-1}U)$ for some non-empty open subset $U\subset B$, then we would obtain $0\neq\varphi^{[n]}\in \sO_T(-nE_T)(h^{-1}(U))$ for $n>0$ divisible enough. But $E_T\equiv -(K_T+D_T')$ is $h$-ample, thus $E_T$ dominates $B$, and so $\sO_T(-nE_T)(h^{-1}(U))$ must be trivial.

                To deal with the $H^0$ term, we show that the involved $R^1$ sheaf is zero. Since $T$ is klt, by \cite[Theorem 3.3]{Tan18} we have $R^1h_*\sO(K_T+D'_T)=0$ as soon as $(K_T+D_T')-K_T=D_T'$ is $h$-ample. As $D'_T$ is disjoint from $E_T$ (use adjunction) and the latter is $h$-ample, we see that every component of $D_T'$ dominates $B$. Since the relative Picard rank is 1, we deduce that $D_T'$ is indeed $h$-ample.
            \end{enumerate}
        \end{enumerate}
The proof is complete.
\end{proof}

\subsection{Two classification results}\label{section:lift_BLE_type}

In this sub-section we classify the log smooth minimal models of the surfaces that appear in \autoref{prop:genus_1_forces_Z_linear_CY_cond}. We start with those of Enriques--type.

\begin{proposition} \label{prop: classification-enriques-type}
    Let $(S,D)$ be a globally $F$-split log smooth projective surface pair such that
    \begin{itemize}
        \item $K_S+D \sim_{\mathbb{Q}} 0$;
        \item $f \colon S \to E$ is a $\mathbb{P}^1$-bundle over an elliptic curve $E$;
        \item $D$ is irreducible.
    \end{itemize}
Then $f|_D \colon D \to E$ is an \'etale Galois morphism of degree 2.
If $f_T \colon (T, D_T)\coloneqq (S,D) \times_E D \to D$ denotes the base change of $(S,D) \to E$ via $D \to E$, then either
\begin{enumerate}
    \item $(T, D_T) \cong (\mathbb{P}^1 \times D, D_0 +D_{\infty})$ where $D_0$ and $D_\infty$ are sections; or,
    \item $(T, D_T) \cong (\mathbb{P}_D(\mathcal{O} \oplus \mathcal{L}), D_0 +D_{\infty})$, where $\mathcal{L}$ is a non-trivial $2$-torsion line bundle on $D$ and $D_0$ and $D_\infty$ are sections.
\end{enumerate}
In both cases the isomorphism identifies $f_T$ with the structural morphism onto $D$.
Moreover, if we denote by $\tau$ the involution on $D$ such that $D/\tau=E$, there exists an involution $\sigma$ on $(T,D_T)$ such that $(T,D_T)/\langle\sigma\rangle \cong (S,D)$ and $\tau \circ f_T=f_T \circ \sigma$.
\end{proposition}

\begin{proof}
Note that, by adjunction, $D \to E$ is a 2:1 cover, thus automatically separable in characteristic $p \neq 2$.
As $(S,D)$ is globally $F$-split, we deduce by $F$-adjunction that $D$ is globally $F$-split, and that the generic fibre is geometrically $F$-split by \cite[Proposition 5.7]{Eji19}. 
In particular, the map $D \to E$ is separable also in characteristic 2.
As $K_E$ and $K_D$ are trivial, by Riemann--Hurwitz, we deduce that $f|_D \colon D \to E$ is \'etale of degree 2. Note that in the characteristic  2 case, $f|_D \colon D \to E$ coincides with the Verschiebung morphism of $D$.

By hypothesis, we have $S=\mathbb{P}_E(\mathcal{E})$ for some vector bundle $\mathcal{E}$ of rank 2 on $E$. Note that the base change $T=S \times_E D=\mathbb{P}_D(f|_D^*\mathcal{E})$ is a projective bundle, and by construction, it has two disjoint sections given by $D \times_E D=D_0+D_{\infty}$. By construction, there is an automorphism $\sigma \in \operatorname{Aut}(\mathbb{P}_D(f|_D^*\mathcal{E}))$ such that $\sigma(D_0)=D_\infty$ and $\sigma^2$ is the identity. 
The existence of two disjoint sections shows that we have a splitting of
$f|_D^*\mathcal{E}$ into sum of line bundles. Up to twisting with a line bundle on $D$, we can assume that  $f|_D^*\mathcal{E} \cong \mathcal{O}_D \oplus \mathcal{L}$.
Therefore, up to reordering, the normal bundles to $D_0$ and $D_\infty$ satisfy 
$$\mathcal{N}_{D_0/\mathbb{P}(f|_D^*\mathcal{E})} \cong \mathcal{L} \text{ and } \mathcal{N}_{D_\infty/\mathbb{P}(f|_{D}^*\mathcal{E})} \cong \mathcal{L}^{\vee}.$$ 
As $\sigma$ is an isomorphism mapping $D_0$ to $D_\infty$, we deduce that $\mathcal{L} \cong \mathcal{L}^{\vee}$, which implies that $\mathcal{L}$ is 2-torsion. The two cases (a) and (b) are distinguished on whether $\mathcal{L}$ is trivial or not.
Note that for an ordinary elliptic curve in characteristic 2, there is a unique non-trivial 2-torsion line bundle.
\end{proof}

We now deal with the classification of crepant smooth minimal models of pairs $(X,D)$ with $K_X+D \sim 0$.
This fixes a mistake in \cite[Proposition 5.9]{BBKW24}: surfaces of Enriques--type appear over fields of characteristic $p=2$. 

\begin{proposition} \label{lem: class_logsnc_CY}
Let $k$ be an algebraically closed field of characteristic $p$ and let $(S, D)$ be a log smooth projective surface pair with $K_S+D \sim 0$. 
Suppose $D \neq 0$. 
Then there exists a crepant blow-up $f \colon (Y, D_Y) \to (X,D)$  and a birational contraction $h \colon (Y, D_Y) \to (Z, D_Z)$ such that one of the following holds:
    \begin{enumerate}
        \item[(i)] \label{i} $(Z, D_Z+E_Z) \cong (\mathbb{P}^2, C)$ where $C$ is an elliptic curve;
			\item[(ii)] $(Z, D_Z+E_Z) \cong (\mathbb{P}^2, L_1+L_2+L_3)$ where $L_i$ are three lines in general position; 
			\item[(iii.a)] $ (Z, D_Z+E_Z) \cong (\mathbb{P}_B(M \oplus N), C+D)$, where $B$ is an elliptic curve, $M$ and $N$ are line bundles on $B$, and $C$ (resp.\ $D$) is the section associated to the quotient $M \oplus N \to M$ (resp.\ $M \oplus N \to N$);
            \item [(iii.b)] $p=2$, and $(Z, D_Z) \cong (\mathbb{P}_B(\mathcal{E}), D)$ where $B$ is an elliptic curve, $\mathcal{E}$ a locally free sheaf of rank 2 and $D \to B$ is a $2:1$ \'etale morphism. 
    \end{enumerate}
\end{proposition}

\begin{proof}
The classification of \cite[Proposition 5.9]{BBKW24} is correct except for the case
where $Z$ has a smooth projective surface with a $\mathbb{P}^1$-bundle $g \colon Z \to B$ onto an 
elliptic curve $B$ and $(Z, D_Z)$ satisfies $K_Z+D_Z \sim 0$. Note that $D_Z$ is horizontal and $D_Z \cdot F=2$, where $F$ is a 
fibre of $g$. Here we need to distinguish two cases.

As proven in \cite[Lemma 5.9]{BPRZ}, either $D_Z$ is reducible with 
two components $D_Z=D_1+D_2$ such that each restriction $D_i \to B$ is 
an isomorphism (this correspond to case (iii.a)), or $D_Z \to B$ is a 2-to-1 finite cover (corresponding to (iii.b)). 
As proved in \cite[Lemma 5.10]{BPRZ}, this latter 
case occurs only when $p=2$ if $K_Z+D_Z \sim 0$.
\end{proof}

\begin{remark}
    Surfaces in case (iii.b) exist (see for example \cite[Example 5.11]{BPRZ}) and appear naturally as irreducible components of degenerations of globally $F$-split Enriques surfaces in characteristic 2 (see \cite[Proposition 2.3 and Theorem 2.7]{Sch03}).
\end{remark}

\subsection{Case with a genus 1 curve in the conductor} \label{subsec: genus1_conductor}
In this sub-section we prove \autoref{main_thm} in case each connected component of the normalization has linearly trivial log canonical divisor, or is of Enriques--type.

To begin with, we construct log lifts of the smooth minimal models. Here is the case of surfaces of Enriques--type:

\begin{theorem} \label{thm: liftability_Enriques_type}
    Let $(S,D)$ be a globally $F$-split log smooth projective surface pair such that
    \begin{enumerate}
        \item $K_S+D \sim_{\mathbb{Q}} 0$;
        \item $f \colon S \to E$ is a $\mathbb{P}^1$-bundle over an elliptic curve $E$;
        \item $D$ is irreducible.
    \end{enumerate}
Then there exists a log lifting $\mathbf{f} \colon (\mathbf{S}, \mathbf{D}) \to \mathbf{E}$ over $W(k)$, where $\mathbf{E}$ and $\mathbf{D}$ are the canonical lifting of $E$ and $D$.
\end{theorem}

\begin{proof}
    Note that $D$ and $E$ are globally $F$-split elliptic curves, and denote by $\mathbf{D}$ and $\mathbf{E}$ their canonical liftings.
    The involution $\tau$ induced by $D\to E$ lifts to an involution $\widetilde{\tau}$ on $\mathbf{D}$ by \autoref{thm: canonical-log-liftings}. By canonicity we have $\mathbf{D}/\langle\widetilde{\mathbf{\tau}}\rangle \cong \mathbf{E}$.
    The morphism $f|_D \colon D \to E$ is \'etale of degree 2, and it lifts to a morphism of the canonical liftings $\mathbf{V} \colon \mathbf{D} \to \mathbf{E}$.
    Consider the base change $(T, D_T)=(S,D)\times_E D=(\mathbb{P}_D(M \oplus N), D_0 +D_{\infty})$ constructed in \autoref{prop: classification-enriques-type} and
    let $(\mathbf{T}, \mathbf{D}_0+\mathbf{D}_{\infty})$ be the canonical lifting as in \autoref{def: canonical_proj_bundle_split_ell}.
    Using \autoref{prop: canonical_lifting_split_projective}, we can lift the involution $\sigma$ to an involution $\widetilde{\sigma}$ of $\mathbf{T}$ such that that $\widetilde{\sigma}(\mathbf{D}_0)=\mathbf{D}_\infty$.
    If we define $(\mathbf{S}, \mathbf{D})$ to be the quotients of $(\mathbf{T}, \mathbf{D}_0 + \mathbf{D}_{\infty})$ by $\widetilde{\sigma}$, we have the following commutative diagram:

        $$\begin{tikzcd}
        &&& (S, D) \arrow[dd, left, "f"] && 
        (\mathbf{S}, \mathbf{D})
        \arrow[from=dlll, crossing over, below right] \arrow[from=ll, hook]\arrow[dd, "\mathbf{f}"]\\
        (T, D_0 + D_\infty) \arrow[dd, "f_T"] \arrow[rr, hook] \arrow[urrr] 
        && (\mathbf{T}, \mathbf{D}_0 +\mathbf{D}_{\infty}) \\
        &&& E \arrow[rr, hook] && \mathbf{E} \\
        D \arrow[rr, hook]\arrow[urrr] && 
        \mathbf{D}
        \arrow[from=uu, crossing over] 
        \arrow[urrr, below right]
        \end{tikzcd}$$
Since $D\to E$ is \'{e}tale, the action of $\langle \tau\rangle$ on $D$ is free. Since $f_T\colon T\to D$ is equivariant, the action of $\sigma$ on $T$ is also free. This implies that the actions of $\widetilde{\tau}$ and $\widetilde{\sigma}$ on $\mathbf{T}$ and $\mathbf{D}$ respectively, are free as well; thus taking their geometric quotients commutes with arbitrary base-change over $W(k)$. So
$\mathbf{f}\otimes_{W(k)}k\cong f$, which concludes the proof.
\end{proof}

For surfaces with linearly trivial log canonical divisor, the smooth minimal models given in \autoref{lem: class_logsnc_CY} have log lifts over $W(k)$: see \cite[Proof of Theorem 5.13]{BBKW24} for cases (i), (ii) and (iii.a) (see also \autoref{def: canonical_proj_bundle_split_ell} for (iii.a)), while case (iii.b) is of Enriques--type and thus the previous theorem applies.

Next we prove the lifting statement for the surfaces as in \autoref{prop:genus_1_forces_Z_linear_CY_cond} that are not necessarily smooth minimal. We need a few lemmata on plt and log canonical surface singularities.

\begin{lemma} \label{lem: plt-log-Cartier}
    Let $(S, D)$ be a plt surface singularity such that $K_S+D \sim 0$ (in particular, $D$ is reduced).
    Then one of the following holds:
    \begin{enumerate}
        \item $D=0$ and $S$ has canonical singularities;
        \item $D \neq 0$ and $(S,D)$ is log smooth.
    \end{enumerate}
\end{lemma}

\begin{proof}
    We first show that $S$ has canonical singularities. Suppose by contradiction this is not the case.
    As $(S,D)$ is a plt surface singularity, $S$ is $\mathbb{Q}$-factorial and klt.
    Let $\pi \colon Y \to S$ be the minimal resolution, with $\text{Ex}(\pi) = \bigcup_{i=1}^m E_i$ and write $K_Y+\sum_{i=1}^m b_i E_i = \pi^*K_S$, where $b_i <1$. 
    By the negativity lemma, $b_i \geq 0$ and, as $K_Y$ is not $\pi$-trivial, there exists at least a $b_i>0$.
    As $K_S+D$ is Cartier and $(S, D)$ is log canonical, this implies that  $a(E_i;S,D)=-1$, contradicting the pair is plt.

    Thus $S$ has canonical singularities and, since $S$ has dimension 2, $K_X$ is a Cartier divisor. Since $K_S+D \sim 0$, this implies that also $D$ is Cartier.
    If $D=0$, we are in case (a).
    If $D \neq 0$, since $(X,D)$ is plt, we have $D$ is regular by \cite{kk-singbook}*{Theorem 2.31}, which forces $S$ to be regular and $(S, D)$ to be log smooth, concluding the proof.
\end{proof}

\begin{lemma} \label{lem: nice-log-res}
     Let $(S, D)$ be a log canonical surface singularity such that $K_S+D \sim 0$ (in particular, $D$ is reduced).
     Then there exists a log resolution $f\colon Y \to (S,D)$ such that 
     $K_Y+f_*^{-1}D+E \sim f^*(K_S+D)$, where $E$ is an effective integral divisor contained in $\Ex(f)$.
\end{lemma}

\begin{proof}
    Passing to a dlt modification (see \cite{Tan18}*{Theorem 4.7}), we can suppose that $(S,D)$ is dlt and $K_S+D\sim0$.
    We now distinguish three cases according to the dimension of the minimal dlt stratum that passes through a closed point $x \in X$.
    If $x \notin \Supp(D)$, then $X$ has canonical singularities and we take $Y$ to be the minimal resolution of $x \in X$.
    If $x$ belongs to only one irreducible component of $D$, we have that $(X,D)$ is plt and we conclude that $(X,D)$ is log smooth by \autoref{lem: plt-log-Cartier}.
    Finally, if $x$ is a 0-dimensional stratum, we have $(S,D)$ is log smooth by definition of dlt.
\end{proof}

Note that for a globally $F$-split projective log smooth pair $(Y,D)$ with $K_Y+D\sim 0$, then every connected component of $D$ is either an ordinary elliptic curve or a circles of $\mathbb{P}^1$s, so we can talk about their canonical liftings as defined in \autoref{section:can_lift}.

\begin{proposition} \label{prop: gen_canonical_lift_BBKW}
    Let $(X, D)$ be a globally $F$-split reduced normal surface pair.
    Assume that one of the following holds:
    \begin{enumerate}
        \item $K_X+D \sim 0$,
        \item $(X, D)$ is of Enriques--type.
    \end{enumerate}
    Let $f \colon Y \to (X,D)$ be a log resolution such that 
    $K_Y+D_Y+E \sim f^*(K_X+D)$, where $E$ is reduced and $f$-exceptional, and $D_Y=f^{-1}_*D$.
    
    Then $f$ admits a lifting $\widetilde{f}\colon (\mathcal{Y}, \mathcal{D}_Y+ \mathcal{E}) \to (\mathcal{X}, \mathcal{D})$ where  $\mathcal{D}_Y+\mathcal{E}$ is the canonical lifting of $D_Y+E$.  
    Moreover, $\mathcal{X}$ is projective over $W(k)$ and $K_\sX+\sD$ is $\bQ$-Cartier.
\end{proposition}

\begin{proof} 
    Let $f\colon Y \to (X,D)$ be a log resolution as in the hypothesis, whose existence is guaranteed by \autoref{lem: nice-log-res} in case (a) and by \autoref{lemma:BLE_is_canonical} in case (b).
    Write $\Ex(f)=E+F,$ where $K_Y+D_Y+E=f^*(K_X+D)$ and $F$ are trees of $(-2)$-curves.
    We now follow the proof of \cite{BBKW24}*{Theorem 5.13}, recalling the main steps for the convenience of the reader.
    Possibly after passing to a higher birational model obtained by a blow-up of $Y$ at points on the boundary, by \autoref{lem: class_logsnc_CY} in case $(a)$, or by definition in case $(b)$,
    there exists a proper birational contraction  $h \colon (Y, D_Y+E) \to (Z, D_Z+E_Z)$ such that we have the following diagram  
		\[
		\begin{tikzcd}
		(Y,D_Y+E) \arrow{r}{f} \arrow{d}{h} & (X, D) \\
		(Z,D_Z+E_Z),
		\end{tikzcd}
		\]
		where $\Ex(f) = E+F$, $D_Z=h_*D_Y$, $E_Z = h_*E$, and $\Ex(h)=\sum G_i$, and $(Z, D_Z)$ either belongs to one of the cases (i), (ii), (iii.a) or (iii.b) of \autoref{lem: class_logsnc_CY}, or is smooth minimal of Enriques--type. For simplicity, we subsume case (iii.b) in the Enriques--type.
        As $(Y, E)$ is globally $F$-split, so is $(Z, E_Z)$ and in particular $E_Z$.
We consider the following lifting $(\mathcal{Z}, \mathcal{D}_Z+\mathcal{E}_Z)$ for $(Z, D_Z+E_Z)$:
\begin{enumerate}
    \item[(i)] $(\mathcal{Z}, \mathcal{D}_{Z}+\mathcal{E}_Z) \simeq  (\mathbb{P}^2_{W(k)}, \mathcal{C}_{\can})$, where $\mathcal{C}_{\can}$ is the canonical lifting of the elliptic curve $C$;
    \item[(ii)] $(\mathcal{Z}, \mathcal{D}_{Z}+\mathcal{E}_Z) \simeq (\mathbb{P}^2_{W(k)}, (xyz=0))$, i.e. the standard toric lifting;
    \item[(iii.a)] $(\mathcal{Z}, \mathcal{D}_{Z}+\mathcal{E}_Z)$ is the canonical lifting constructed in \autoref{def: canonical_proj_bundle_split_ell};
    \item [(E-T)] $(\mathcal{Z}, \mathcal{D}_Z)$ is the lifting of surfaces of Enriques-type constructed in \autoref{thm: liftability_Enriques_type}.
\end{enumerate}

Now let $A$ be a very ample Cartier divisor on $X$ and let $A_Y=f^*A$. 
Write \[h^*h_* A_Y = A_Y + \sum_i a_i G_i\] for some integers $a_i$. Since $Z$ is smooth, the Weil divisor $L:=h_*A_Y$ is a Cartier.
We need the following claim, the proof of which is found in \cite{BBKW24}*{Claim 5.14}.

\begin{claim}\label{claim-good-lifting}
		There exists a lifting $\widetilde{h} \colon (\mathcal{Y}, D_{\mathcal{Y}}+\mathcal{E}) \to (\mathcal{Z}, \mathcal{D}_{\mathcal{Z}}+\mathcal{E}_{\mathcal{Z}})$ of $h \colon (Y, D+E) \to (Z, D_Z+E_Z)$ together  with liftings $\mathcal{F}$ of $F$, $\mathcal{G}_{i}$ of $G_i$ and $\mathcal{L}$ of $L$ such that the line bundle
			\[\mathcal{A}_{\mathcal{Y}}:=\widetilde{h}^*\mathcal{L}{\big(}-\sum_{i} a_i \mathcal{G}_{i}\big)\]
			satisfies $\mathcal{A}_{\mathcal{Y}}|_{\mathcal{E}} \sim 0$.
            Moreover, $\mathcal{A}_{\mathcal{Y}}|_{\mathcal{E}}$ is the canonical lifting of $A|_{E}$ in the sense of \cite{MS87} (in the case of elliptic curves) and \autoref{prop: canonical_lifts_circlesP1} (in the case of circles of $\mathbb{P}^1$s).
        \hfill $\lozenge$
\end{claim}

We decompose $f$ as
        $$(Y,D_Y+E+F) \overset{\varphi}{\longrightarrow}
        (T,D_T+E_T) \overset{\psi}{\longrightarrow}
        (X,D)$$
where $\varphi$ is the contraction of $F$, and $\psi$ is the contraction of $E_T=\varphi_*E$. Since $E$ is either an elliptic curve or a circle of $\bP^1$s, by adjunction along $E_T$ we see that $\varphi(F)$ and $E_T$ are disjoint. Let $(\mathcal{Y}, \mathcal{E} + \mathcal{F})$ be the lifting constructed in \autoref{claim-good-lifting}. By \cite{BBKW24}*{Corollary 5.5} we can contract $\mathcal{F}$ to get a lifting $\widetilde{\varphi} \colon (\mathcal{Y}, \mathcal{E} + \mathcal{F}) \to (\mathcal{T}, \mathcal{E}_{T})$ of $\varphi$.
Since $A_Y|_F \sim 0$ and $H^1(F, \mathcal{O}_F)=0$, we deduce that $\mathcal{A}_{\mathcal{Y}}|_{\mathcal{F}} \sim 0$. So $\mathcal{A}_{\mathcal{Y}}$ is $\widetilde{\varphi}$-trivial and it descends to a line bundle $\mathcal{A}_{\mathcal{T}}$ on $\mathcal{T}$. 
As $\mathcal{A}_{\mathcal{T}}|_{\mathcal{E}_{T}} \sim 0$, by \cite{BBKW24}*{Proposition 5.4} we see that $\sA_\sT$ induces a projective morphism $\widetilde{\psi}\colon 
            (\sT,\sD_\sT+\sE_\sT)
            \to
            (\sX,\sD)$
that lifts $\psi$. Thus $\widetilde{f} = \widetilde{\psi} \circ \widetilde{\varphi}$ is the desired lifting of $f$. By construction we have:

    \begin{claim}\label{claim:can_lift_ample_on_cond}
    The line bundle $\sA_\sT$ descends to an ample line bundle $\sA$ on $\sX$ such that $\sA|_{\sD}$ is the canonical lift of $A|_D$.
    \hfill $\lozenge$
    \end{claim}

In particular, $\sX$ is projective. The fact that $K_\sX+\sD$ is $\bQ$-Cartier is proved as in \cite[Theorem 6.8]{BBKW24}.
\end{proof}


\begin{theorem} \label{thm: main_KX-trivial}
    Let $S_0$ be a projective globally $F$-split CY surface with normalization $(S,D)$. Assume that $p>2$ and that each connected component $(S_i, D_i)$ satisfies either
    \begin{enumerate}
        \item $K_{S_i}+D_i \sim 0$; or,
        \item $(S_i,D_i)$ is of Enriques--type.
    \end{enumerate} 
    Then there exists a projective strong slc lifting $\mathcal{S}_0$ of $S_0$ over $W(k)$.
\end{theorem}

\begin{proof}
If $S_0$ is normal, we conclude by \autoref{prop: gen_canonical_lift_BBKW}. Suppose that $S_0$ is not normal and let $(S,D,\tau)$ be its normalization triple. Let $A_0$ be an ample line bundle on $S_0$, and let $A$ be its pullback to $S$. Applying \autoref{prop: gen_canonical_lift_BBKW} and \autoref{claim:can_lift_ample_on_cond} to each connected component of $S$, there exists a projective strong slc lifting $(\sS,\sD)$ of $(S,D)$ where $\sD$ is the canonical lifting of $D$, and an ample line bundle $\sA$ on $\sS$ such that $\sA|_{\sD}$ is the canonical lifting of $A|_D$. 

Combining \cite{MS87}*{Appendix, Theorem 1.2} and \autoref{cor:lift_involutions_on_P1}, there is a canonical involution $\tilde{\tau}$ of $(\sD^\nu,\Diff_{\sD^\nu}(0))$ lifting $\tau$. By \autoref{lem: suff-cond-log-lift-nonnormal}, the quotient $\sS_0=\sS/R(\tilde{\tau})$ is a strong slc lifting of $S_0$. Clearly $\sS_0$ is flat and proper over $W(k)$, and it remains to establish its projectivity.

By canonicity of our liftings, $\sA|_{\sD^\nu}$ is $\widetilde{\tau}$-invariant: this induces an action on the total space of sections $\Tot(\sA|_{\sD^\nu})$. The graph of this action is a finite equivalence relation $R(\widetilde{\tau}) \rightrightarrows \Tot(\sA|_{\sD^\nu})$. Let $R\rightrightarrows \Tot(\sA)$ be the smallest set-theoretic equivalence relation that contains the image of $R(\widetilde{\tau})$ under the finite morphism $\Tot(\sA|_{\sD^\nu})\to \Tot(\sA)$. Notice that its special fiber $R\otimes k\rightrightarrows \Tot(A)$ is a finite equivalence relation, since $A$ descends to $S_0$. Therefore $R\rightrightarrows \Tot(\sA)$ is finite. If the quotient $\Tot(\sA)/R$ exists as a scheme, then one deduces as in \cite[Proposition 9.48]{kk-singbook} that some multiple of $\sA$ descends to $\sS_0$, establishing the desired projectivity of $\sS_0$.

Let $K=\Frac(W)$ and notice that the generic equivalence relation $R_ K\rightrightarrows \Tot(\sA_K)$ is the identity over the complement of 
    $$\Tot(\sA_K)\times_{\sS_K}\sD_K=
    \Tot(\sA_K|_{\sD_K}).$$ 
As $(\sS_K,\sD_K)$ is an lc surface pair the boundary $\sD_K$ is semi-normal. Therefore $\Tot(A_K|_{\sD_K})$ is semi-normal as well, and by \cite[Proposition 33]{Kollar_Quotients_finite_equiv_relations} the quotient $\Tot(\sA_K)/R_K$ exists as an algebraic space. By \cite[Theorem 1.4]{Wit20} it follows that the quotient $\Tot(\sA)/R$ exists as an algebraic space. By \cite[Corollary 48]{Kollar_Quotients_finite_equiv_relations} it is in fact a scheme; so our proof is complete. 
\end{proof}

\subsection{First step towards the main theorem}\label{section:higher_index}
We recall that our goal is to prove:

\begin{theorem}\label{thm:main_with_technical_statement}
Let $S_0$ be a projective semi-log canonical globally $F$-split surface over $k$ with $K_{S_0}\equiv 0$. Then there exists a projective strong slc lifting $\mathcal{S}_0$ of $S_0$ over $W(k)$.
\end{theorem}

Building on the results obtained in this section, we prove this theorem in many cases: the remaining ones will be dealt with in the next section. We introduce the following terminology:

\begin{definition}[$4A_1$-curve]\label{def:4A_1}
Let $(S,D)$ be a log canonical surface pair with reduced boundary: we say that an irreducible component $C$ of $D$ is a $4A_1$-curve of the pair $(S,D)$ if $C \simeq \mathbb{P}^1$ and that $S$ has four $A_1$-singularities along $C$.
\end{definition}

\begin{lemma}\label{lemma:4A1_not_BLE_type}
Suppose that $(S,D)$ has a $4A_1$-curve. Then $(S,D)$ is not of Enriques--type and $K_S+D$ is not linearly trivial.
\end{lemma}
\begin{proof}
If $(S,D)$ is of Enriques--type then it has canonical singularities by \autoref{lemma:BLE_is_canonical} and its discrepancies are non-negative; if $K_S+D\sim 0$ then the discrepancies of $(S,D)$ are integers. Blowing-up an $A_1$-singularity along a $4A_1$-curve produces an exceptional divisor with discrepancy $-1/2$ over $(S,D)$---see \autoref{section:local_analysis} below---, so the statement follows.
\end{proof}

\begin{proof}[Proof of \autoref{thm:main_with_technical_statement} when $p>2$]
By \cite[Theorem 0.1]{Tanaka_Abundance_for_slc_surfaces} we have $K_{S_0}\sim_\bQ 0$. Let $(S,D,\tau)$ be the normalization triple of $S_0$. Then $K_{S}+D\sim_\bQ 0$ and $(S,D)$ is globally sharply $F$-split by \autoref{prop:GFS_ascends_normalization}. 

Let $(S,D)=\bigsqcup_i (S_i,D_i)$ be the decomposition into connected components. Thanks to \autoref{lemma:4A1_not_BLE_type}, we can partition these components into three classes:
    \begin{enumerate}
        \item[(I)] $K_S+D\sim 0$, or $(S,D)$ is of Enriques--type;
        \item[(II)] $K_S+D$ is not linearly trivial, $(S,D)$ is not of Enriques--type, but there is no $4A_1$-curve;
        \item[(III)] $(S,D)$ has some $4A_1$-curve.
    \end{enumerate}
We choose proper strong log liftings over $W(k)$ of these components as follows:
    \begin{enumerate}
        \item[(I)] we take the proper strong log lifting afforded by \autoref{prop: gen_canonical_lift_BBKW};
        \item[(II)] we take an arbitrary proper strong log lifting, afforded by the proof of \cite[Theorem 6.8]{BBKW24};
        \item[(III)] we take a well-chosen proper strong log lifting, in a sense to be made precise shortly.
    \end{enumerate}
For the classes (I) and (II), the given reference shows that the log canonical divisor of the chosen strong log lifting is $\bQ$-Cartier; for class (III), this will be part of our desiderata. So we have a proper strong log lifting
$(\sS,\sD)=\bigsqcup_i (\sS_i,\sD_i)$ over $W(k)$ with $K_\sS+\sD\sim_\bQ 0$. By \autoref{lemma:normalization_of_lifted_boundary} it holds that 
        $$(\sD^\nu,\Diff_{\sD^\nu}(0))\otimes k= (D^\nu,\Diff_{D^\nu}(0))$$
and the support of $\Diff_{\sD^\nu}(0)$ is \'{e}tale over $W(k)$.

Next, thanks to \autoref{corollary:possible_CY_log_curves} we also partition the irreducible components of $D=\bigcup_{i,j}D_{ij}$ into three classes:
    \begin{enumerate}
        \item[$(\alpha)$] those $D_{ij}$ with normalization a genus one curve (and by adjunction we see that $D_{ij}$ is normal already and isolated in $D$);
        \item[$(\beta)$] those $D_{ij}$ that are rational and such that the support of $\Diff_{D_{ij}^\nu}(D-D_{ij})$ contains at most three points;
        \item[$(\gamma)$] the $4A_1$-curves.
    \end{enumerate}
The involution $\tau$ preserves these three classes (with a small abuse of language, since $\tau$ is defined on $D^\nu$). 
We claim that the involution $\tau$ lifts to an involution $\widetilde{\tau}$ of $(\sD^\nu,\Diff_{\sD^\nu}(0))$ over $W(k)$. It suffices to consider each class separately.
    \begin{enumerate}
        \item[$(\alpha)$] If a component $D_{ij}$ of $D_i$ is of type $(\alpha)$, then $(S_i,D_i)$ is necessarily of type (I) by \autoref{prop:genus_1_forces_Z_linear_CY_cond}. By \autoref{thm: main_KX-trivial} the lift $\sD_{ij}\subset \sD_i$ is the canonical lift of $D_{ij}$, so \autoref{thm: canonical-log-liftings} applies.

        \item[$(\beta)$] In that case the lift of $\tau$ is constructed in \autoref{cor:lift_involutions_on_P1}.

        \item[$(\gamma)$] If $D_{ij}$ is a $4A_1$-curve, then the corresponding $(S_i,D_i)$ is of type (III). We will require that the lift $\sD_{ij}\subset \sD_i$ of $D_{ij}$ is the canonical one, in the sense of \autoref{def:canonical_lift_P1_4pts}, so \autoref{prop:functoriality_can_lift_P1_4pts} can be applied.
    \end{enumerate}
We define $\sS_0=\sS/R(\widetilde{\tau})$. By \autoref{lem: suff-cond-log-lift-nonnormal} the surface $\sS_0$ is a proper strong slc lifting of $S_0$ over $W(k)$. 
It remains to show that $\sS_0$ is projective over $W(k)$. Suppose first that $K_{S_0}$ is not linearly trivial. Then $H^2(S_0,\sO_{S_0})=H^0(S_0,\sO(K_{S_0}))^\vee =0$. In that case any ample line bundle on $S_0$ lifts to an ample line bundle on $\sS_0$. So we may assume that $K_{S_0}\sim 0$, which implies that every component of $(S,D)$ belongs to class (I). In this situation \autoref{thm: main_KX-trivial} applies.
Therefore the proof is complete, modulo the following claim:
\begin{claim}
If $(S_i,D_i)$ has a $4A_1$-curve, then there exists a strong log lifting $(\sS_i,\sD_i)$ over $W(k)$ such that: $K_{\sS_i}+\sD_i\sim_\bQ 0$, and the liftings in $\sD_i$ of the $4A_1$-curves are the canonical ones (\autoref{def:canonical_lift_P1_4pts}).
\end{claim}

The proof of this claim requires a careful study of the global geometry. This is the goal of the next section. The overall strategy is explained as the beginning of \autoref{section:4A_1}. In \autoref{section:strategy} (more precisely, in \autoref{prop:existence_of_log_liftings} and \autoref{prop:lift_equiv_MMP}) we reduce the problem to a precise lifting problem, which is solved in \autoref{section:equiv_lift_dP}, \autoref{section:equiv_lift_conic_bd_I}, \autoref{section:equiv_lift_Gamma_non_zero}, \autoref{section:equiv_lift_E_reducible}, \autoref{section:equiv_lift_conic_bd_II} and \autoref{section:equiv_lift_m=4}. 
\end{proof}

\subsubsection{The case $p=2$}\label{section:p=2}
We prove \autoref{main_thm}, whose notation we keep, in characteristic $p=2$. The main observation is the following one:

\begin{lemma}\label{lemma:GFS_when_p=2}
Assume that $p=2$. Then $S_0$ has only separable nodes and $K_{S_0}\sim 0$.
\end{lemma}
\begin{proof}
Since $S_0$ is globally $F$-split, it is weakly normal \cite{Hochster_Roberts_Purity_Frobenius_Local_cohom}. Inseparable nodes are not weakly normal \cite[Lemma 2.3.7]{Posva_Gluing_stable_families_surfaces_mix_char} so $S_0$ has only separable nodes. The Frobenius trace $\Tr_{S_0}\colon H^0(S_0,\sO((1-2^e)K_{S_0}))\to H^0(S_0,\sO_{S_0})$ is split for every $e\geq 1$: taking $e=1$, we obtain a non-zero global section of $\sO(-K_{S_0})$. As $K_{S_0}\equiv 0$ it follows that $K_{S_0}\sim 0$.
\end{proof}

Consequently, the strategy of proof of \autoref{thm: main_KX-trivial} applies to $S_0$ and yields \autoref{main_thm} in the case $p=2$. We spell out the details in this sub-section, as we need to adapt \autoref{lem: suff-cond-log-lift-nonnormal} to this setting. 

\begin{lemma}\label{lemma:involution_in_char_2}
Let $\tau$ be a non-trivial involution of $\bP^1_k$. Then there is a coordinate $x$ and $b\in k^\times$ such that $\tau$ is given by $\tau^*(x)=x+b$.
\end{lemma}
\begin{proof}
The Riemann--Hurwitz formula shows that $\tau$ has at least a fixed point. Let $x$ be a coordinate of $\bP^1_k$ such that $\tau(\infty)=\infty$. Then $\tau^*(x)=ax+b$ for some $a,b\in k$. As $\tau\neq \id$ is an involution, we get $a=1$ and $b\in k^\times$.
\end{proof}

\begin{theorem}\label{thm:main_thm_KX_trivial_p=2}
     Let $S_0$ be a projective globally $F$-split CY surface with normalization $(S,D)$. Assume that $p=2$. Then there exists a projective strong slc lifting $\mathcal{S}_0$ of $S_0$ over $W(k)$.
\end{theorem}
\begin{proof}
By \autoref{lemma:GFS_when_p=2} we have $K_S+D\sim 0$ and $S_0$ has only separable nodes. So we are in position to follow the proof of \autoref{thm: main_KX-trivial}. The hypothesis $p>2$ was used at two places only: to find the lift $\widetilde{\tau}\circlearrowright \sD$ of $\tau\circlearrowright D$ (through the use of \autoref{cor:lift_involutions_on_P1}), and to show that $\sS_0\otimes k=S_0$ (through \autoref{lem: suff-cond-log-lift-nonnormal}). We will show that both steps also work when $p=2$.

First, we consider the problem of lifting the involution $\tau\circlearrowright (D^\nu,\Diff(0))$ to the $\sD$ constructed in \autoref{prop: gen_canonical_lift_BBKW}. The genus one curves in $D^\nu$ are ordinary and permuted by $\tau$, and \cite[Appendix, Theorem 1.2]{MS87} applies to their canonical liftings. So it remains to consider the components of genus zero. Let $D_1$ be such a component, and let $D_2=\tau(D_1)$. If $D_1\neq D_2$ then it is easy to lift the isomorphism $\tau\colon (D_1,\Diff(0))\cong (D_2,\Diff(0))$ to $(\sD_1,\Diff(0))\cong (\sD_2,\Diff(0))$ as $\sD_1\cong \bP^1_{W(k)}$ and each divisor $\Diff$ is the sum of two points. Assume next that $D_1=D_2$. Then by \autoref{lemma:involution_in_char_2} there is a coordinate $x$ such that $\tau|_{D_1}$ is given by $x\mapsto x+b$. This description shows that $\infty=[1:0]\in D_1$ is the unique fixed point. As $\Diff_{D_1}(0)$ is $\tau|_{D_1}$-invariant of degree $2$, we have
$\Diff(0)=[\alpha:1]+[\beta:1]$ (with respect to the coordinate $x$) for some distinct $\alpha,\beta\in k$. Notice that $b=\alpha+\beta$. Now take a coordinate $\mathbf{x}$ of the lift $\sD_1\cong \bP^1_{W(k)}$ of $D_1$ such that $\mathbf{x}$ restricts to $x$. Then
        $$\Diff_{\sD_1}(0)=[\boldsymbol\alpha:1]+[\boldsymbol\beta:1]$$
with respect to $\mathbf{x}$ where $\boldsymbol\alpha$ and $\boldsymbol\beta$ are elements of $W(k)$ reducing respectively to $\alpha$ and $\beta$. Define $\widetilde{\tau}|_{\sD_1}$ by $\mathbf{x}\mapsto -\mathbf{x}+\boldsymbol\alpha+\boldsymbol\beta$. Then clearly $\widetilde{\tau}|_{\sD_1}$ is a log involution of $(\sD_1,\Diff(0))$ lifting $\tau|_{D_1}$.

Let $\sS_0=\sS/R(\widetilde{\tau})$: it remains to show that $\sS_0\otimes k= S_0$. We follow the proof of \cite[Theorem 4.1.9]{Posva_Gluing_stable_families_surfaces_mix_char} with some minor adaptations. To begin with, we have:

\begin{claim}\label{claim:qt_and_fiber}
$(\sD^\nu/\langle \widetilde{\tau}\rangle)\otimes k= D^\nu/\langle \tau\rangle$.
\end{claim}
\begin{proof}\renewcommand{\qedsymbol}{$\lozenge$}
In general, we have a finite universal homeomorphism $\rho\colon D^\nu/\langle \tau\rangle \to (\sD^\nu/\langle \widetilde{\tau}\rangle)\otimes k$ which measures the failure of commutativity of quotienting and taking the closed fiber. Let $\sU\subset \sD^\nu$ be the locus where the action of $\widetilde{\tau}$ is free, and $\sU'$ be its image in the quotient $\sD^\nu/\langle \widetilde{\tau}\rangle$. Then $\sU'$ is relatively dense over $W(k)$, and $\rho$ is an isomorphism over $\sU'\otimes k$. As $D^\nu/\langle \tau\rangle$ is regular, it suffices to show that $(\sD^\nu/\langle \widetilde{\tau}\rangle)\otimes k$ is also regular at the remaining points outside $\sU'\otimes k$.

Let $x\in \sD^\nu$ be a closed point which is fixed by $\widetilde{\tau}$. Then
        $$\widehat{\sO}_{\sD^\nu,x}= W(k)\llbracket t\rrbracket,$$ 
$\widetilde{\tau}$ acts $W(k)$-linearly on this complete local ring and $\widetilde{\tau}(t)\in (t)$. Then (see e.g.\ \cite[Proposition 10.3.48]{Liu02}),
        $$\left(\widehat{\sO}_{\sD^\nu,x}\right)^{\langle \widetilde{\tau}\rangle}
        = W(k)\llbracket t\cdot \widetilde{\tau}(t)\rrbracket,$$
the reduction of which modulo $2$ is clearly regular. The claim follows.
\end{proof}

Returning to the proof of \autoref{thm:main_thm_KX_trivial_p=2}, let $\mathcal{W}=\sD/R_{\sD}(\widetilde{\tau})$ and $W=\mathcal{W}\otimes k$. 
By construction we have a commutative diagram
    $$\begin{tikzcd}
    & D^\nu \arrow[rr, hook] \arrow[d] \arrow[ddl] && 
    \sD^\nu \arrow[d] \arrow[ddl] \\
    & D \arrow[rr, hook, crossing over]  
    && \sD \arrow[dd] \\
    D^\nu/\langle \tau\rangle \arrow[rr, hook] \arrow[dr] && 
    \sD^\nu/\langle \widetilde{\tau}\rangle \arrow[dr]  \\
    & W \arrow[from=uu, crossing over] \arrow[rr, hook] &&
    \mathcal{W}.
    \end{tikzcd}$$
As in \cite[Observation 4.1.15]{Posva_Gluing_stable_families_surfaces_mix_char} we see that to show $\sS_0\otimes k=S_0$ it suffices to prove the following: if $s\in \sO_D$ is such that $s^2\in \sO_W$, then $s\in \sO_W$ already.

Take a section $s\in \sO_D$ such that $s^2\in \sO_W$. As $\sO_W$ embeds into the sub-sheaf of $\tau$-invariants of $\sO_{D^\nu}$, we see that 
    $$\big(\tau(s)-s\big)^2=\tau(s^2)-s^2=0.$$
As $D^\nu$ is reduced, we get that $s\in \sO_{D^\nu}$ is $\tau$-invariant. By \autoref{claim:qt_and_fiber} we can find $v\in \sO_{\sD^\nu}$ which is $\widetilde{\tau}$-invariant and lifts $s$. If we can show that there is a correction term $c\in \sO_{\sD^\nu}$ that is $\widetilde{\tau}$-invariant and such that $v+2c\in \sO_\sD$, then $v+2c\in \sO_{\mathcal{W}}$ and so $s=(v+2c)|_D\in \sO_W$ as desired. This is achieved by the argument in \cite{Posva_Gluing_stable_families_surfaces_mix_char} starting after the proof of Lemma 4.21, as this argument also holds when $p=2$ (\footnote{
    More precisely, the argument starts on p.165 with ``Now let us treat the general case (...)" and ends with the first paragraph on p.168.
    As the reader can verify, that argument does not use \cite[Lemma 4.1.20]{Posva_Gluing_stable_families_surfaces_mix_char}---and fortunately so, as that lemma does not hold without the assumption that $|G|$ is invertible in $A$.
}). This completes our proof.
\end{proof}


%% file: 6_case_of_4A1.tex

\section{Strong slc lifings with $4A_1$-curves in the conductor}\label{section:4A_1}
In this section we study the geometry of projective globally $F$-split slc surface pairs $(S,D)$, where $D$ is a reduced divisor, $K_S+D\sim_\bQ 0$ is not Cartier and $(S,D)$ has some $4A_1$-curve (\autoref{def:4A_1}), with the aim of completing the proof of \autoref{thm:main_with_technical_statement}. 
Throughout this section we assume that $p>2$.

\medskip
The difficulty is the following. In the context of \autoref{thm:main_with_technical_statement}, we are given a gluing log involution $\tau$ of $(D^\nu,\Diff_{D^\nu}(0))$ and we must lift it over $W(k)$ together with a slc lifting of $(S,D)$. Let $C_1,C_2$ be $4A_1$-curves such that $\tau$ restricts to
        $$\tau\colon (C_1,\Diff_{C_1}(0))\cong (C_2,\Diff_{C_2}(0)).$$
Take an arbitrary strong log lift $(\sS,\sD)$ of $(S,D)$, and let $\sC_i$ be the irreducible component of $\sD$ specializing to $C_i$. By definition each $\Supp(\Diff_{C_i}(0))$ contains four points. Whether $\tau$ can be lifted to a log isomorphism $(\sC_1,\Diff_{\sC_1}(0))\cong (\sC_2,\Diff_{\sC_2}(0))$ depends on the boundary divisors $\Diff_{\sC_i}(0)$, essentially because of \autoref{rmk:PGL2_not_4transitive}. So we need to be careful about the strong log lifting of $(S,D)$ that we choose.

We will show that we can find a strong log lifting of $(S,D)$ which induces the canonical lifts of all its $4A_1$-curve, in the sense of \autoref{def:canonical_lift_P1_4pts}. Then we will be able to use \autoref{prop:functoriality_can_lift_P1_4pts} to conclude.

Instead of working directly with $(S,D)$, we study its index one canonical cover (that is, the cyclic cover of $S$ defined by the $\bQ$-Cartier divisor $K_S+D$) and try to lift it together with the induced cyclic group action over the Witt vectors. 
It may seem that we are making our task more difficult by working with these equivariant covers: the advantage is that we have an explicit description of the equivariant minimal models, thanks to \cite{Bayle_Beauville_Birational_involutions_P2}. In each case we are able to construct an equivariant lifting over $W(k)$.

\medskip
This section is divided in several parts, through which we have tried to keep the notations consistent. In \autoref{section:local_analysis} and \autoref{section:rationality} we collect basic properties of the surfaces $(S,D)$. In \autoref{section:index_one_covers} we analyze their index one canonical covers. We classify the equivariant minimal models of these covers in \autoref{section:classification_min_models}. In \autoref{section:strategy} we explain our strategy to construct adequate equivariant log liftings. Finally, in the remaining sub-sections we implement this strategy in each case obtained in our classification.

\begin{notation}\label{notation:4A1_case}
From now on, $(S,D)$ is a projective globally sharply $F$-split surface pair over $k$, with reduced boundary $D$, such that $K_S+D\sim_\bQ 0$ is not Cartier. We allow $S$ to be reducible. We assume that $(S,D)$ has some $4A_1$-curve: in particular $D\neq 0$.
\end{notation}

In \autoref{section:strategy} we will also assume that a log involution $\tau$ of $(D^\nu,\Diff_{D^\nu}(0))$ is given, but it will play no role in the other sub-sections.

\medskip
Let us mention that cases where $(S,D)$ has a $4A_1$-curve indeed occur: the classification results of \autoref{section:classification_min_models} give several ways of constructing examples. Here is one:

\begin{example}\label{example:4A_1}
We start with an ordinary elliptic curve $(E,e)$ embedded in $\bP^2_k$ by $\sO(3e)$. By \cite[Lemma 2.7]{CTW17} the pair $(\bP^2,E)$ is globally $F$-split. Let $\sigma\in \Aut_k(E,e)$ be the standard involution of the elliptic curve: as it acts on $H^0(E,\sO(3e))^\vee$, it extends to an involution $\sigma_{\bP^2}$ of $\bP^2$. We can find a basis of $H^0(E,\sO(3e))^\vee$ for which the $\mu_2$-action induced by $\sigma$ is given by a diagonal matrix $M$ of order $2$. Since $\sigma_{\bP^2}$ is given by the image of $M$ through
        $$\Aut_k \left( H^0(E,\sO(3e))^\vee\right)
        \cong 
        \GL_3(k) \twoheadrightarrow \PGL_3(k)\cong \Aut_k(\bP^2)$$
we see that $M$ is not a scalar matrix. Thus one eigenspace of $M$ has dimension two. This eigenspace defines a line $L$ in $\bP^2$ which belongs to the fixed locus of $\sigma_{\bP^2}$. 

By \autoref{lemma:sing_qt} we see that the codimension one part $C$ of the fixed locus of $\sigma_{\bP^2}$ is regular, and by \autoref{prop:sing_qt_pair} it meets $E$ transversely. So if $d$ is the degree of $C$, the involution $\sigma$ has at least $3d$ distinct fixed points. Now $\sigma$ has exactly four fixed points, namely the sub-group $E[2]$ of $2$-torsion points. So $d=1$ and $C=L$. Write $L\cap E=\{e_1,e_2,e_3\}\subset E[2]$: since the intersection is transversal, these three points are distinct. Since $e_1+e_2+e_3=0$ in the group law of $(E,e)$, it is easily seen that $L\cap E=E[2]\setminus \{e\}$ (and in fact $e$ is the point defined by the one-dimensional eigenspace of the matrix $M$).

Let $b\colon X\to \bP^2$ be the blow-up of $L\cap E$, and let $E'=q^{-1}_*E$. As $b\colon (X,E')\to (\bP^2,E)$ is crepant we obtain $K_X+E'\sim 0$ and $(X,E')$ is globally $F$-split \cite[Lemma 2.4]{BBKW24}. 
The $\mu_2$-action generated by $\sigma_{\bP^2}$ lifts to a $\mu_2$-action on $X$, and a local computation shows that $\Exc(b)\cap E'$ are isolated fixed points. In particular, the codimension one part of the fixed locus of the action on $X$ is $L'=b^{-1}_*L$.

Let $c\colon X\to T$ be the contraction of the $(-2)$-curve $L'$, and write $E_T=c_*E'$. Then $(T,E_T)$ is CY and globally $F$-split \cite[Lemma 2.3]{BBKW24}, and the $\mu_2$-action descends to $(T,E_T)$. The point is that the $\mu_2$-action on $T$ is now free in codimension one, and has exactly four fixed points along $E_T$.

Let $q\colon T\to S$ be the geometric $\mu_2$-quotient, and let $D=q(E_T)$. Then $q$ is unramified in codimension one, so $K_S+D\sim_{\bQ}0$. As $q_*\sO_T((p^e-1)E_T)= \sO_S((p^e-1)D)$ and $\sO_S\to q_*\sO_T$ is split by the Reynolds operator, the argument of \cite[Lemma 1.1.9]{Brion-Kumar} shows that $(S,D)$ is globally $F$-split. The four $\mu_2$-fixed points along $E_T$ create four $A_1$-singularities along $D$: so $(S,D)$ is the desired example.
\end{example}

\subsection{Local analysis}\label{section:local_analysis}
We begin with some local analysis of $S$ along a neighborhood of the $4A_1$-curves.

\medskip
Let $C$ be a $4A_1$-curve of $(S,D)$ and let $s\in S$ be one of the four singular points along $C$. According to \cite[Theorem 3.32]{kk-singbook}, a model for the completion of $(S,C)$ at $s$ is given by
        \begin{equation}\label{eqn:formal_description_A1_sing}
            \widehat{S}= \Spec \left( \frac{k\llbracket x,y,v\rrbracket}{v^2-xy}\right),
            \quad \widehat{C}=(x=0=v).
        \end{equation}
Using the formal coordinates $x$ and $y$, we see that $\omega_{\widehat{S}}$ is invertible and generated by $\frac{dx\wedge dy}{2v}=\frac{dx\wedge dv}{x}$, that 
    \begin{equation}\label{eqn:formal_description_I}
    \text{the non-invertible reflexive sheaf }\omega_{\widehat{S}}(\widehat{C})\text{ is generated by }
    \frac{dx\wedge dv}{x^2} \text{ and }
    \frac{dx\wedge dv}{xv},
    \end{equation}
and that
    \begin{equation}\label{eqn:formal_description_II}
    \omega_{\widehat{S}}^{\otimes 2}(2\widehat{C})\text{ is invertible and generated by }
    \frac{(dx\wedge dv)^{\otimes 2}}{x^3}.
    \end{equation}

\begin{lemma}\label{lemma:4A_1_is_isolated}
The curve $C$ disjoint from $D-C$, and the only singularities of $S$ along $C$ are the four $A_1$-singularities $q_1,\dots,q_4$. Moreover, 
    $$\left(C,\Diff_C(0)=\frac{1}{2}\sum_{i=1}^4q_i\right)$$ 
is globally sharply $F$-split.
\end{lemma}
\begin{proof}
Recall that $C$ is regular by definition. The local analysis above shows that 
        $$\Diff_{C^\nu}(D-C)=\frac{1}{2}\sum_{i=1}^4q_i+G,$$
where $G$ is effective by \cite[Proposition 2.35]{kk-singbook}. Since $K_S+D\sim_\bQ 0$ and $C^\nu\cong \bP^1$, by adjunction along $C^\nu$ we obtain that $G=0$. By \cite[Proposition 2.35]{kk-singbook} again, this implies that $D-C$ is disjoint from $C$
and that the only singularities of $S$ along $C$ are the $q_i$'s.

We show that $(C,\Diff(0))$ is globally sharply $F$-split. Observe that $(1-p^e)(K_C+\Diff(0))$ is a $\mathbb{Z}$-divisor for every $e\geq 1$. Functoriality of the trace morphisms of Frobenius gives the following commutative diagram
		$$\begin{tikzcd}
        H^0(S,F_*\sO_S((1-p^e)(K_S+D))) \arrow[rr, "\Tr_{(S{,}D)}"] \arrow[d]
        && H^0(S,\sO_S) \arrow[d, two heads] \\
        H^0(C, F_{*}\mathcal{O}_{C}((1-p^e)(K_C+\Diff_C(0))) \arrow[rr, "\Tr_{(C{,}\Diff(0))}"] &&
        H^0(C,\sO_C).
		\end{tikzcd}$$
As $(S,D)$ is globally sharply $F$-split, the top arrow $\Tr_{(S,D)}$ is surjective for some $e\geq 1$ and thus $\Tr_{C, \Diff_C(0)}$ is also surjective. This shows that $(C, \Diff_C(0))$ is globally sharply $F$-split.
\end{proof}

\subsection{Rationality of the underlying surface}\label{section:rationality}

\begin{proposition}\label{prop:special_surface_is_rational}
Let $(S,D)$ be as in \autoref{notation:4A1_case}. Then $S$ is a rational surface.
\end{proposition}

\begin{proof}
Suppose, for contradiction, that $S$ is not a rational surface. Let $f\colon S'\to S$ be a minimal log resolution. We can write
        $$K_{S'}+D' =f^*(K_S+D), \quad \text{where }D'\geq 0.$$
We claim that all the connected components of $D'$ are trees of rational curves. Indeed, if a connected component of $\Supp(D')$ has genus one, then all its components appear with coefficient one on $D'$ \cite[3.27,28,31]{kk-singbook}. Replacing $S'$ by an intermediate dlt modification of $S$, we may assume that $D'$ is an integral divisor with a connected component of genus one. By \cite[Lemma 2.4]{BBKW24} the pair $(S',D')$ is globally sharply $F$-split, so it follows from \autoref{prop:genus_1_forces_Z_linear_CY_cond} and \autoref{lemma:4A1_not_BLE_type} that $H^0(S',\sO(K_{S'}+D'))$ is non-trivial. 
Since $f_*\sO(K_{S'}+D')=\sO(K_S+D)$ we find that $H^0(S,\sO(K_S+D))$ is non-trivial either, a contradiction since $K_S+D$ is not Cartier.

Returning to the case where $S'$ is the minimal log resolution, we have that $K_{S'}$ is anti-effective and that $S'$ is not rational. We run a $K_{S'}$-MMP to obtain a birational contraction onto a $\mathbb{P}^1$-bundle over a regular curve $E$ of positive genus:
        $$\varphi\colon S'\longrightarrow Y=\mathbb{P}_E(\mathcal{M}).$$
Let $D_Y=\varphi_*D'$. As $D'$ is a union of trees of rational curves, every irreducible component of $D_Y$ is contained in a fiber of $Y\to E$. As $K_Y+D_Y\sim_\bQ 0$, we have for a closed fiber $F$ of $Y\to E$ that
        $$0=(K_Y+D_Y)\cdot F=K_Y\cdot F=-2,$$
which yields the desired contradiction. 
\end{proof}

\subsection{Index one covers}\label{section:index_one_covers}
Next, we show that the index one covers of $(S,D)$ are quite restrained.

\begin{proposition}\label{prop:log_Gor_cover}
Let $(S,D)$ be as in \autoref{notation:4A1_case}. Assume that $S$ is irreducible. Then:
    \begin{enumerate}
        \item the Cartier index $m$ of $K_S+D$ is $2$ or $4$;
        \item if $m=2$ (resp. $m=4$) then the pair $(S,D)$ has at most two (resp. exactly one) $4A_1$-curves.
    \end{enumerate}
The index one cover of $S$ induced by $K_S+D$ produces a finite crepant morphism
        $$q\colon (T,E+\Gamma)\longrightarrow (S,D)$$
such that:
    \begin{enumerate}\setcounter{enumi}{2}
        \item $(T,E+\Gamma)$ is projective, log canonical and globally sharply $F$-split, with $K_T+E+\Gamma\sim 0$;
        \item there exists a $\mu_m$-action on the pair $(T,E+\Gamma)$ such that $q$ is the geometric $\mu_m$-quotient and $q(E+\Gamma)=D$;
        \item the irreducible components of $E$ are disjoint regular ordinary curves of genus $1$, and they are contained in the regular locus of $T$;
        \item $E$ is preserved by the $\mu_m$-action, and its image $q(E)$ is the set of $4A_1$-curves of $(S,D)$.
    \end{enumerate}
Let $C$ be a $4A_1$-curve of $(S,D)$. Then:
    \begin{enumerate}\setcounter{enumi}{6}
        \item $\pi^*C$ has $m/2$ connected components that are permuted by $\mu_{m/2}(k)$;
        \item if $E_1$ is an irreducible component of $E$ which dominates $C$, then the morphism $E_1\to C$ is the quotient by the restricted action of the subgroup $\mu_{2}\subseteq \mu_{m}$, and it ramifies over the four $A_1$-singularities of $S$ contained in $C$;
        \item in particular, the action of $\mu_2$ on $E_1$ is not free.
    \end{enumerate}
\end{proposition}

\begin{remark}
In \autoref{cor:m=4_does_not_occur} we will see that the case $m=4$ actually does not occur.
\end{remark}

\begin{proof}[Proof of \autoref{prop:log_Gor_cover}]
Let us take $m>1$ minimal such that $m(K_S+D)\sim 0$. Write $m=2^nr$ with $(2,r)=1$. Since $K_S+D$ is $2$-Cartier but not Cartier along the $4A_1$-curves of $(S,D)$, we have $n\geq 1$. 
Let $T$ be the cyclic covering of $S$ given by the $2^n$-Cartier divisor $r(K_S+D)$. More precisely, we let
        $$T=\Spec_S \left(\bigoplus_{i=0}^{2^{n}-1}\sO_S(-ir(K_S+D))\right)
        \overset{\pi}{\longrightarrow} S$$
where the algebra structure on the direct sum is given by the choice of a trivialization $\sO_S(-2^nr(K_S+D))\cong \sO_S$ (\footnote{
    The choice of isomorphism does not matter much, since we are not in characteristic $2$.
}). The direct sum decomposition yields an action of $\mu_{2^{n}}$ on $T$ such that the morphism $\pi$ is the geometric quotient.

The morphism $\pi$ is \'{e}tale over the regular locus of $S$, and the pushforward $\pi_*\sO_T$ is an $S_2$ $\sO_S$-module. Therefore $T$ is normal, and 
        $$K_T+\pi^*D\sim_\bQ \pi^*(K_S+D)\sim_\bQ 0.$$
It follows from \cite[Corollary 2.43]{kk-singbook} that $(T,\pi^*D)$ is lc. It is globally sharply $F$-split since $(S,D)$ is so, by \cite[Lemma 2.5]{BBKW24}. We set $E\leq \pi^*D$ be the pullback of the $4A_1$-curves of $(S,D)$. 

Let us analyze the morphism $\pi$ formal-locally above an $A_1$-singularity $s\in S$ that lies along a $4A_1$-curve $C$. Say $r=2a+1$. Using the coordinates of \autoref{eqn:formal_description_A1_sing} and the generators of \autoref{eqn:formal_description_I} and \autoref{eqn:formal_description_II}, we see that the algebra of $T$ is generated over $\widehat{\sO}_{S,s}\cong k\llbracket x,y,v\rrbracket /(v^2-xy)$ by the dual sections
    $$\mathcal{T}_1=\left(\frac{(dx\wedge dv)^{\otimes r}}{x^{3a+2}}\right)^\vee 
    \quad \text{and}\quad 
    \mathcal{U}_1=\left(\frac{(dx\wedge dv)^{\otimes r}}{x^{3a+1}v}\right)^\vee
    \quad \text{in}\quad \omega_S^{-r}(-rD)\otimes \widehat{\sO}_{S,s}$$
together with
    $$\mathcal{G}_{2}=\left(\frac{(dx\wedge dv)^{\otimes 2r}}{x^{3r}}\right)^\vee 
    \quad \text{in}\quad \omega_S^{-2r}(-2rD)\otimes \widehat{\sO}_{S,s}.$$
The subscripts of these generators indicate the weights for the $\mu_{2^n}$-action. The relations between them are given by
    $$\mathcal{T}_1^2=x\mathcal{G}_{2}, \quad \mathcal{U}_1^2=y\mathcal{G}_{2}, \quad \mathcal{U}_1\mathcal{T}_1=v\mathcal{G}_{2},
    \quad \mathcal{G}_{2}^{2^{n-1}}=u$$
where $u\in \widehat{\sO}_{S,s}^\times$. In other words we have
        $$T\otimes_S \widehat{\sO}_{S,s}\cong 
        \Spec \left( \frac{k\llbracket x,y,v,\mathcal{U}_1,\mathcal{T}_1,\mathcal{G}_2\rrbracket}
        {\left(v^2=xy, \ \mathcal{T}^2_1=x\mathcal{G}_2, \ \mathcal{U}_1^2= y\mathcal{G}_2, \ \mathcal{U}_1\mathcal{T}_1=v\mathcal{G}_2, \ \mathcal{G}_2^{2^{n-1}}=u\right)}\right).$$
Since $u$ has a $2^{n-1}$-th root in $\widehat{\sO}_{S,s}$, we see that $T\otimes_S \widehat{\sO}_{S,s}$ has $2^{n-1}$ connected components, permuted cyclically by the subgroup $\mu_{2^{n-1}}\subset \mu_{2^n}$, and each of those is isomorphic to
    $$\Spec \left( \frac{k\llbracket x,y,v,\mathcal{U},\mathcal{T}\rrbracket}
        {\left(v^2=xy, \ \mathcal{T}^2=x, \ \mathcal{U}^2= y, \ \mathcal{U}\mathcal{T}=v\right)}\right)\cong \widehat{\bA}^2_{\mathcal{U},\mathcal{T}}.$$
On each of these connected components, the divisor $\pi^*C$ is equal to $(\mathcal{T}=0)$ and hence regular.
Moreover:
    $$T\otimes_S \widehat{\sO}_{S,s}/(x,v)\cong \bigsqcup_{i=1}^{2^{n-1}}\Spec\left(
    \frac{k\llbracket \mathcal{U},\mathcal{T}\rrbracket}{(\mathcal{T}^2,\mathcal{U}\mathcal{T})}
    \right).$$
This analysis shows the following. First, the surface $T$ is regular along $\pi^*C$. Moreover, if $C$ is a $4A_1$-curve of $(S,D)$, then:
    \begin{itemize}
        \item each of the irreducible component of $\pi^*C$ is regular, has degree $2$ over $C$ and ramifies over the four $A_1$-singularities. It does not ramify over any other point, since $(S,C)$ is snc along $C\setminus \Sing(S)$ by \autoref{lemma:4A_1_is_isolated};
        \item $\pi^*C$ has $2^{n-1}$ irreducible components that are permuted cyclically by the sub-group $\mu_{2^{n-1}}(k)$; 
    \end{itemize}    
The Riemann--Hurwitz formula shows that each connected component of $\pi^*C$ has genus $1$. Since $(T,\pi^*D)$ is globally sharply $F$-split we deduce that each component of $\pi^*C$ is ordinary. 

By \autoref{prop:genus_1_forces_Z_linear_CY_cond}, either $K_T+\pi^*D\sim 0$ or $(T,\pi^*D)$ is of Enriques--type. In the latter case $(T,\pi^*D)$ has canonical singularities (\autoref{lemma:BLE_is_canonical}), and thus $K_T+\pi^*D$ is Cartier by 
\cite[Theorem 2.29]{kk-singbook}.
In both cases this implies that, after localizing on $S$ if necessary, the sheaf 
        $$\pi_*\sO_T(K_T+\pi^*D)=
        \bigoplus_{i=0}^{2^n-1}\sO_S((-ir+1)(K_S+D))$$
has a nowhere vanishing section. Therefore $(-ir+1)(K_S+D)\sim 0$ locally on $S$ for some $0\leq i\leq 2^n-1$. 
If $ir \neq 1$, then the Cartier index of $K_S+D$ would divide $ir-1$ and be smaller than $2^nr$: this is a contradiction. Hence it is necessary that $r=i=1$, and we obtain that $K_S+D$ is $2^n$-Cartier.

Our local analysis has established that if $l$ is the number of $4A_1$-curves of $(S,D)$, then the divisor $\pi^*D$ contains $2^{n-1}l$ disjoint regular curves of genus $1$. By \autoref{lemma:at_most_two_elliptic_curves} below, we have $2^{n-1}l\leq 2$. Thus either $l=1$ and $n\in \{1,2\}$, or $l=2$ and $n=1$. The proof is complete.
\end{proof}

\begin{lemma}\label{lemma:at_most_two_elliptic_curves}
Let $(T,E+\Gamma)$ be an irreducible projective log CY surface pair. Assume that $E$ is a reduced divisor and its irreducible components are regular curves of genus $1$. Then $E$ has at most two irreducible components.
\end{lemma}
\begin{proof}
This is essentially a consequence of the Koll\'{a}r--Shokurov connectedness principle (see \cite[Proposition 3.3.1]{Prokhorov_Complements_on_log_surfaces} and \cite{FW24}), but let us give a quick proof for latter reference. Replacing $T$ by its minimal resolution and taking the crepant boundary, we may assume that $T$ is regular. Run a $K_T$-MMP: we obtain a birational contraction $h\colon T\to U$ and a smooth Mori fiber space structure on $U$. Notice that no irreducible component of $E$ is contracted by $h$. Let $E_U=h_*E$ so that $K_U+E_U\sim_{\mathbb{Q}} -h_*\Gamma$. We have two cases:
    \begin{enumerate}
        \item $U\cong \bP^2$. --- Then let $E_{U,1}$ be an irreducible component of $E_U$, with degree $d$ and normalization $E_{U,1}^\nu$. Since $E_{U,1}^\nu$ has genus $1$, by considering the long exact sequence of cohomology induced by
                $$0\to \sO_{E_{U,1}}\to \sO_{E^\nu_{U,1}}\to \mathcal{Q}\to 0,$$
        we see that $(d-1)(d-2)\geq 2$ and therefore $d\geq 3$. As $K_{\bP^2}$ has degree $-3$, this implies that $E_U=E_{U,1}$ and that $E_U$ is regular. In this case $E$ is irreducible.
        \item $U$ is a $\bP^1$-bundle over a regular curve $B$. --- Then every irreducible component of $E_U$ dominates $B$. So by intersecting $K_U+E_U$ with a general fiber, we see that $E_U$ has at most two components.
    \end{enumerate}
The lemma is proved.
\end{proof}

\subsection{Classification of equivariant minimal models}\label{section:classification_min_models}
We classify the smooth equivariant minimal models of the index one covers constructed in \autoref{section:index_one_covers}.

\medskip 
Let $(T,E+\Gamma, \mu_m)$ be as in the statement of \autoref{prop:log_Gor_cover}. Consider the $\mu_m$-equivariant log resolution $\varphi\colon T'\to T$ of $(T,E+\Gamma)$ of \autoref{prop:min_equiv_log_resolution}, and let $E'$ be the strict transform of $E$. We write
        $$K_{T'}+E'+\Gamma'=\varphi^*(K_T+E+\Gamma)\sim 0$$
where $\Gamma'$ is the sum of $\varphi^{-1}_*\Gamma$ with an effective $\varphi$-exceptional reduced divisor.
Then $\mu_m$ acts on $(T',E'+\Gamma')$ and both divisors $E'$ and $\Gamma'$ are $\mu_m$-invariant. 
By \cite[Lemma 2.4]{BBKW24} the pair $(T',E'+\Gamma')$ is globally sharply $F$-split. 

Now let $H\subseteq \mu_m$ be any non-trivial subgroup. 
Since $K_{T'}$ is anti-effective, by \autoref{prop:min_equiv_model} an $H$-equivariant $K_{T'}$-MMP produces a birational contraction
    $$h\colon T'\longrightarrow U(H)$$ 
such that: $h$ is $H$-equivariant, and either $U(H)$ is a del Pezzo surface with $\rk \Pic(U(H))^H=1$, or it admits a minimal $H$-conic bundle structure. Write $E_{U(H)}=h_*E'$ and $\Gamma_{U(H)}=h_*\Gamma'$. The morphism $h$ does not contract any irreducible component of $E'$, and the $H$-action on $U(H)$ sends $E_{U(H)}$ to itself. Moreover
        \begin{equation}\label{eqn:equiv_min_model_log_CY}
        K_{U(H)}+E_{U(H)}+\Gamma_{U(H)}\sim 0.
        \end{equation}
By adjunction, it follows that the irreducible components of $E_{U(H)}$ are regular and disjoint from each other and from $\Gamma_{U(H)}$. In addition the restriction
$h|_E\colon E\longrightarrow E_{U(H)}$
is an equivariant isomorphism.

\begin{lemma}\label{lemma:min_equiv_model_rational_case}
Suppose that $T$ is rational. Then $m=2$ and $U=U(\mu_2)$ satisfies one of the following:
    \begin{enumerate}
        \item $U\cong \bP^2$;
        \item $U$ is a del Pezzo of degree $2$ and the $\mu_2$-action is given by the Geiser involution;
        \item $U$ is a del Pezzo of degree $1$ and the $\mu_2$-action is given by the Bertini involution;
        \item $U\cong \bP^1\times\bP^1$ and the $\mu_2$-action is given by $(x,y)\mapsto (y,x)$;
        \item there is a minimal $\mu_2$-equivariant conic bundle $f\colon U\to \bP^1$ such that either:
            \begin{enumerate}  
                \item $f$ is $\mu_2$-invariant in the sense of \autoref{invariant-def}, or
                \item $f$ is smooth and non-$\mu_2$-invariant.
            \end{enumerate}
    \end{enumerate}
In every case, $E_U$ is irreducible and $\Gamma_U=0$.
\end{lemma}
\begin{proof}
Since $U$ is rational and minimal, the action $U\circlearrowleft \mu_2$ is one of those described in \cite[Theorem 1.4]{Bayle_Beauville_Birational_involutions_P2}, which we have listed in the statement. To show that $E_U$ is irreducible and that $\Gamma_U=0$, we distinguish two cases:
    \begin{enumerate}
        \item Suppose that $U$ has a $\mu_2$-conic bundle structure $f\colon U\to \bP^1$. By \autoref{eqn:equiv_min_model_log_CY} and the fact that the components of $E_U$ are genus one curves, we see that $E_U$ must be irreducible of degree $2$ over the base $\bP^1$. This implies that $\Gamma_U$ is contained in fibers. By minimality of $f$ and invariance of $E_U$, we see that $E_U$ intersects every component of any fiber, and so $\Gamma_U=0$ by \autoref{eqn:equiv_min_model_log_CY}.
        
        \item If $U$ is a $\mu_2$-minimal del Pezzo surface,        
        let $g\colon U\to Y$ be the (non-equivariant) minimal model. Then either $Y \simeq \mathbb{P}^2$ or $Y \simeq \mathbb{P}^1 \times \mathbb{P}^1$. First assume that $Y \simeq \mathbb{P}^2$. 
        Then $E_{U}\to g_*E_{U}$ is birational and $K_{\bP^2}+g_*E_{U}\sim -g_*\Gamma_U$. As in the proof of \autoref{lemma:at_most_two_elliptic_curves}, we see that $g_*E_U$ is irreducible and regular of degree $3$. So $g_*\Gamma_U=0$. Write
                $$K_U+E_U\sim g^*(K_{\bP^2}+g_*E_U)+F$$
        where $F$ is $g$-exceptional. Since $g$ is a composition of blow-ups and $g_*E_U$ is regular, we see that $F\geq 0$. From \autoref{eqn:equiv_min_model_log_CY} and $K_{\bP^2}+g_*E_U\sim 0$ it follows that $\Gamma_U+F\sim 0$. As $U$ is projective, this implies that $F=0=\Gamma_U$.

        If $Y\cong \bP^1\times \bP^1$, then any component of $g_*E_U$ must have bidegree $(a,b)$ with $a,b\geq 2$. This implies that $g_*E_U$ is irreducible, regular, and that $g_*\Gamma_U=0$. As above, it follows that $\Gamma_U=0$.
    \end{enumerate}
To conclude we show, by way of contradiction, that $m=2$. If $m=4$, by \autoref{prop:log_Gor_cover} the divisor $E$ would have $m/2=2$ connected components: the same would hold for $E_{U(\mu_2)}$, which is impossible as we have seen. Therefore it must hold that $m=2$. The proof is complete.
\end{proof}



\begin{lemma}\label{lemma:min_equiv_model_elliptic_case}
Assume that $T$ is not rational. Then there is a minimal $\mu_m$-conic bundle $f\colon U=U(\mu_m)\to B$ such that:
    \begin{enumerate}
        \item $f$ is non-$\mu_m$-invariant, and $B$ is an ordinary genus $1$ curve;
        \item if $m=2$, then $f$ is smooth;
        \item if $m=4$ then $\Gamma_U=0$, the divisor $E_U$ has two components which are sections of $f$ and are permuted by $\mu_4$, and the $\mu_4$-action on $B$ is faithful.
    \end{enumerate}
\end{lemma}
We will get more details about the case $m=4$ in \autoref{section:equiv_lift_m=4}, and eventually see that it cannot occur. The argument does not require the following sub-sections, but since it is relatively long we have placed it at the end of the section.
\begin{proof}
Since $T$ is not rational, $U$ is not a del Pezzo surface, and therefore it admits a minimal $\mu_m$-conic bundle structure over a regular curve $B$ of positive genus. The components of $E_U$ are not contained in the fibers, which shows that the genus of $B$ is $1$. Since the components of $E_U$ are ordinary, so is $B$ by \cite{Sil}*{Theorem 3.1}.

The quotient $U/\mu_m$ is a birational model of our original surface $S$. By \autoref{prop:special_surface_is_rational} it is rational, and therefore it is fibered onto $\mathbb{P}^1$. Since $U/\mu_m$ is fibered over the quotient of $B$ by the image of $f_*\colon \mu_m(k)\to \Aut_k(B)$, we deduce that this image is non-trivial, in other words that $f$ is non-invariant.

For the rest of the proof, we let $\sigma_U$ be a generator of the image of $\mu_m(k)\to \Aut_k(U)$.

Suppose that $m=2$ and that $f$ is not smooth. Let $F\cup F'$ be a singular fiber. Since $f$ is minimal, $\sigma_U$ sends $F$ to $F'$, so the intersection $u$ of $F$ and $F'$ is a fixed point of $\sigma_U$. The action on $\fm_u/\fm_u^2$ exchanges the directions given by $F$ and $F'$: so there is a basis which this action is given by the matrix
        $$\begin{pmatrix}
        0&1 \\ 
        1&0
        \end{pmatrix}.$$
The eigenvalues of this matrix are $+1$ and $-1$. So if $\widehat{\sigma}_u$ is the induced involution of $\widehat{\sO}_{U,u}$, we see that $\Fix(\widehat{\sigma}_u)$ is a germ of regular curve. Since $\Fix(\sigma_U)\otimes \widehat{\sO}_{U,u}=\Fix(\widehat{\sigma}_u)$ (as shown in the proof of \autoref{lemma:sing_qt}), it follows that there is a regular $\mu_2$-fixed curve passing through $u$. This curve cannot be contained in the fiber $F\cup F'$, so it must dominate $B$. This implies by equivariance of $f$ that the $\mu_2$-action on $B$ is trivial, contradicting what we have established above. 

Assume that $m=4$. Recall that in this case $E_U=E_{U,1}\sqcup E_{U,2}$ has two disjoint regular components which are exchanged by $\sigma_U$ (cf.\ \autoref{prop:log_Gor_cover}). Since the intersection product of $E_U$ with a general fiber is equal to $2$, we see that both restrictions $f|_{E_{U,i}}\colon E_{U,i}\to B$ are isomorphisms. Since $E_{U,1}$, $E_{U,2}$ and $\Gamma_U$ must be pairwise disjoint, it follows from \autoref{eqn:equiv_min_model_log_CY} and Koll\'{a}r--Shokurov connectedness principle that $\Gamma_U=0$. 

For each $i=1,2$, the restriction $\sigma_U^2|_{E_{U,i}}$ is a non-trivial involution with four fixed points: this follows from \autoref{prop:log_Gor_cover} since $h|_E\colon E\to E_U$ is a $\mu_4(k)$-equivariant isomorphism. If $f$ is $\sigma^2_U$-invariant, then it would follow that both $E_{U,i}\to B$ are degree $2$ coverings, which is not the case by what we have established. Thus the morphism $f_*\colon \mu_4(k)\to \Aut_k(B)$ is injective. 
\end{proof}

In most cases it turns out that $\Gamma_U=0$ or $E_U$ is irreducible. Let us describe the special cases where this does not happen.

\begin{lemma}\label{lemma:min_equiv_model_special_case}
Suppose that $\Gamma_U\neq 0$. Then:
    \begin{enumerate}
        \item  $E_U$ is irreducible and $m=2$;
        \item the divisor $\Gamma'$ on $T'$ is a regular curve of genus $1$, preserved by the action of $\mu_2$, and the restriction $h|_{\Gamma'}\colon \Gamma'\to \Gamma_U$ is an isomorphism;
        \item the equivariant morphism $f\colon U\to B$ is smooth and non-invariant, and $E_U$ and $\Gamma_U$ are two disjoint sections.
    \end{enumerate}
\end{lemma}
\begin{proof}
It follows from \autoref{lemma:min_equiv_model_rational_case} and \autoref{lemma:min_equiv_model_elliptic_case} that $m=2$, that $f$ is smooth and that $B$ is a regular genus $1$ curve. Using \autoref{eqn:equiv_min_model_log_CY} we see that $E_U$ is irreducible and that both $E_U$ and $\Gamma_U$ are sections of $f$. 

This implies that one of the component of $\Gamma'$ is an irreducible regular genus $1$ curve. Thus $T$ has a strictly log canonical elliptic singularity, say at $t$, which is disjoint from the support of $E+\Gamma$ by \cite[Proposition 2.28]{kk-singbook}. By the Koll\'{a}r--Shokurov connectedness principle, it follows that $\Gamma'$ is equal to the exceptional divisor over $t$. Hence $h|_{\Gamma'}\colon \Gamma'\to \Gamma_U$ is an isomorphism.
\end{proof}


\begin{lemma}\label{lemma:min_equiv_model_E_reducible}
Assume that $E_U$ is reducible. Then either:
    \begin{enumerate}
        \item $m=4$ (cf.\ \autoref{lemma:min_equiv_model_elliptic_case}), or
        \item $m=2$ and $\Gamma_U=0$, the morphism $f\colon U\to B$ is smooth non-invariant onto an ordinary genus $1$ curve, and the two irreducible components of $E_U$ are disjoint sections.
    \end{enumerate} 
\end{lemma}
\begin{proof}
Assume that $m=2$. From \autoref{lemma:min_equiv_model_rational_case} and \autoref{lemma:min_equiv_model_elliptic_case} we see that $f\colon U\to B$ is a smooth non-invariant morphism onto an ordinary genus $1$ curve. From \autoref{eqn:equiv_min_model_log_CY} we see that $\Gamma_U=0$ and that $E_U$ has two disjoint irreducible components which are sections of $f$.
\end{proof}

\subsection{Equivariant liftings: strategy}\label{section:strategy}
In the previous sub-sections, we have constructed and studied equivariant models of the index one cover of $(S,D)$. We now present a strategy to lift these models, together with the $\mu_m$-actions, over the Witt vectors, and argue that this will enable us to construct an slc log lifting of $(S,D)$ together with a lift of a given gluing involution. By \autoref{cor:m=4_does_not_occur} (which does not use the results of this sub-section), we have in fact $m=2$, so we restrict to this case.

\medskip
Treating each connected component of $S$ separately, we have constructed the following commutative diagram:
    
    \medskip
    \begin{center}
    \begin{tikzpicture}
    \node at (0,0) {$(T', E'+\Gamma')$} ;
    \node at (5,0) {$(U, E_U+\Gamma_U)$} ;
    \node at (0,-1.5) {$(T,E+\Gamma)$} ;
    \node at (5,-1.5) {$B$} ;
    \node at (0,-3.3) {$(S,D)$} ;

    \node at (7.8,-0.75) {$\circlearrowleft \mu_{2,k}$} ;

    \draw (-2,0.7) -- (7,0.7) ;
    \draw (-2,0.7) -- (-2,-2.2) ;
    \draw (-2,-2.2) -- (7,-2.2) ;
    \draw (7,0.7) -- (7,-2.2) ;

    \draw[->] (1.2,0)--(3.8,0) node[midway, above] {$h$} ;
    \draw[->] (0,-0.3) -- (0,-1.2) node[midway, left] {$\varphi$} ;
    \draw[->] (5,-0.3) -- (5,-1.2) node[midway, right] {$f$} ;
    \draw[-] (0,-1.7) -- (0,-2.1) ;
    \draw[->] (0,-2.3) -- (0,-3) node[midway, left] {$q$} ;
    
    \end{tikzpicture}
    \end{center}
where $q$ is the index one cover of $(S,D)$ (\autoref{prop:log_Gor_cover}), $\varphi$ is the minimal log resolution of $(T,E+\Gamma)$ and $h$ is a composition of birational contractions given by a $\mu_{2,k}$-equivariant $K_{T'}$-MMP (\autoref{section:classification_min_models}).

\begin{remark}\label{rmk:double_cover_of_4A1_curve}
The birational map $h\circ \varphi^{-1}$ restricts to an equivariant isomorphism $E\cong E_U$. So it follows from \autoref{prop:log_Gor_cover} and from \autoref{lemma:4A_1_is_isolated} that if $E_U\to C$ is the geometric quotient by the $\mu_{m}$-action on $E_U$, and if $\mathfrak{b}$ is the branch divisor on $C$, then $(C,\frac{1}{2}\mathfrak{b})\cong (\bP^1,\frac{1}{2}\sum_{i=1}^4 q_i)$ is globally sharply $F$-split pair. Conversely, if $E_{U,1}$ is any irreducible component of $E_U$, then the morphism $E_{U,1}\to C$ is the double cover of $(C,\frac{1}{2}\mathfrak{b})$ constructed in \autoref{prop:canonical_lift_P1_4_pts}. 
\end{remark}

We now give a sufficient criterion on how to construct a log lifting of $(S,D)$, together with a lift of a given log involution of $(D^\nu,\Diff_{D^\nu}(0))$ from a lifting of $(T,E+\Gamma)$ with special properties. 

\begin{proposition}\label{prop:existence_of_log_liftings}
Assume that there exists a $\mu_{2,W(k)}$-equivariant lift 
    $$\varphi_{W(k)}\colon (\mathbf{T}',\mathbf{E}'+\mathbf{\Gamma}')\to (\mathbf{T},\mathbf{E}+\mathbf{\Gamma})$$ 
of $\varphi\colon (T',E'+\Gamma')\to (T,E+\Gamma)$ (in the sense of \autoref{def:equiv_lifting})
such that:
    \begin{enumerate}
        \item $\mathbf{T}$ is projective over $W(k)$,
        \item $\mathbf{E}$ is the canonical lifting of $E$
        (in the sense of \autoref{thm: canonical-log-liftings}), and 
        \item the restriction of the $\mu_{2,W(k)}$-action to $\mathbf{E}$ is the canonical lift of the $\mu_{2,k}$-action on $E$ (in the sense of \autoref{thm: canonical-log-liftings}).  
    \end{enumerate}
Then:
    \begin{enumerate}
        \item there exists a strong log lifting $(\mathbf{S},\mathbf{D})$ of $(S,D)$;
        \item any given log involution $\tau$ of $(D^\nu,\Diff_{D\nu}(0))$ lifts to a log involution of $(\mathbf{D}^\nu,\Diff_{\mathbf{D}^\nu}(0))$.
    \end{enumerate}
\end{proposition}
\begin{proof}
First we construct a lift of $(S,D)$ over $W(k)$. Let $\mathbf{q}\colon \mathbf{T}\to \mathbf{S}$ be the geometric $\mu_{2,W(k)}$-quotient and $\mathbf{D}=\mathbf{q}(\mathbf{E}+\mathbf{\Gamma})$: it exists since $\mathbf{T}$ is projective and $\mu_{2,W(k)}$ is finite, cf.\ the beginning of \autoref{section:finite_qts}. By \autoref{prop:restriction_linearly_red_qt} we have $\mathbf{q}\otimes k=q$ and $\mathbf{D}\otimes k=D$. So $(\mathbf{S},\mathbf{D})$ lifts $(S,D)$ over $W(k)$.

We claim that $(\mathbf{S},\mathbf{D})$ is in fact a strong log lifting of $(S,D)$. We have to produce a log resolution of $(\mathbf{S},\mathbf{D})$ whose specialization over $k$ is a log resolution of $(S,D)$: we construct it as follows. By \autoref{prop:simplify_mu2_sing} we have a $\mu_{2,W(k)}$-equivariant birational projective morphism 
        $$\xi_{W(k)}\colon \widetilde{\mathbf{T}}\longrightarrow \mathbf{T}'$$
such that $(\widetilde{\mathbf{T}}, (\xi_{W(k)})^{-1}_*(\mathbf{E}'+\mathbf{\Gamma}')+\Exc(\psi_{W(k)}))$ is relatively snc, and such that the $\mu_{2,W(k)}$-fixed locus on $\widetilde{\mathbf{T}}$ has pure relative codimension one over $W(k)$. Let 
        $$\widetilde{\mathbf{q}}\colon \widetilde{\mathbf{T}}
        \longrightarrow
        \widetilde{\mathbf{S}}=
        \widetilde{\mathbf{T}}/\mu_{2,W(k)}.$$
By \autoref{prop:descent_projectivity} the quotient $\widetilde{\mathbf{S}}$ is projective. The universal property of geometric quotients gives a birational projective morphism $\psi_{W(k)}\colon \widetilde{\mathbf{S}}\to \mathbf{S}$. Note that 
    $$\widetilde{\mathbf{q}}\left( 
    (\xi_{W(k)})^{-1}_*(\mathbf{E}'+\mathbf{\Gamma}')+\Exc(\psi_{W(k)})
    \right)=
    (\psi_{W(k)})^{-1}_*\mathbf{D}+\Exc(\psi_{W(k)}).
    $$
Thus it follows from \autoref{prop:sing_qt_pair} that $(\widetilde{\mathbf{S}},(\psi_{W(k)})^{-1}_*\mathbf{D}+\Exc(\psi_{W(k)}))$ is relatively snc. This shows $\psi_{W(k)}\otimes k$ is also a log resolution of $(S,D)$, and proves that $(\mathbf{S},\mathbf{D})$ is indeed a strong log lifting of $(S,D)$.

It remains to lift the log involution $\tau$ of $(D^\nu,\Diff_{D\nu}(0))$. Let $\mathfrak{C}_{4A_1}$ be the set of irreducible components of $D^\nu$ which are $4A_1$-curves of $(S,D)$. The involution $\tau$ preserves the set of $4A_1$-curves, since they are characterized by having a boundary with four points. So we can write
        $$\tau=(\tau_{4}, \tau_{<4})\in 
        \Aut_k\left( \bigsqcup_{D_i^\nu\in \mathfrak{C}_{4A_1}} (D^\nu_i,\Diff_{D^\nu_i}(0)) 
        \right)
        \times 
        \Aut_k\left( \bigsqcup_{D_j^\nu\notin \mathfrak{C}_{4A_1}} (D^\nu_j,\Diff_{D^\nu_j}(0)) 
        \right).$$
As we have seen at the end of \autoref{section:higher_index}, we can always lift the log involution $\tau_{<4}$. To lift $\tau_4$, we aim to use \autoref{prop:functoriality_can_lift_P1_4pts}: let us verify that we are in an appropriate set-up.

Let $E_1$ be any irreducible component of $E$, let $\mathbf{E}_1$ be the irreducible component of $\mathbf{E}$ specialising to $E_1$, 
and let $\mathbf{C}=\mathbf{q}(\mathbf{E}_1)$ with specialization $C=q(E_1)$ over $k$. Then $C$ is a $4A_1$-curve of $(S,D)$. Consider the cartesian square
        \begin{equation}\label{eqn:qt_square_can_lift}
        \begin{tikzcd}
        E_1\arrow[d, "q|_{E_1}"]\arrow[r,hook] & \mathbf{E}_1 \arrow[d, "\mathbf{q}|_{\mathbf{E}_1}"] \\
        C\arrow[r, hook] & \mathbf{C}.
        \end{tikzcd}
        \end{equation}
Since $C$ is a $4A_1$-curve, it is isolated in $D$, the $\mu_{2,k}$-action sends $E_1$ to itself and $T$ is regular along $E_1$. So $\mathbf{C}$ is isolated in $\mathbf{D}$, the $\mu_{2,W(k)}$-action sends $\mathbf{E}_1$ to itself and $\mathbf{T}$ is regular along $\mathbf{E}_1$. It follows that the ramification locus of $\mathbf{E}_1\to \mathbf{C}$ is the intersection of $\mathbf{E}_1$ with the ramification locus of the $\mu_{2,W(k)}$-action on $\mathbf{T}$. Taking images in $\mathbf{S}$, we find that the branch locus of $\mathbf{E}_1\to\mathbf{C}$ is the intersection $\Sing(\mathbf{S})\cap \mathbf{C}$; by \cite[Proposition 4.5]{kk-singbook} this intersection is the support of $\Diff_\mathbf{C}(0)$.

By \autoref{prop:restriction_linearly_red_qt}, the morphism $q|_{E_1}$ is the quotient by the restricted action of $\mu_{2,k}$ on $E_1$. By assumption the restriction $\mu_{2,W(k)}|_{\mathbf{E}_1}$ is the canonical lift of the $\mu_{2,k}$-action on $E_1$. Therefore the square \autoref{eqn:qt_square_can_lift} is the square \autoref{eqn:canonical_lift_P1_4_pts} constructed in the proof of \autoref{prop:canonical_lift_P1_4_pts} starting with the globally sharply $F$-split pair $(C,\Diff_C(0))$, and the canonical lift of the latter pair (in the sense of \autoref{def:canonical_lift_P1_4pts}) is $(\mathbf{C},\Diff_\mathbf{C}(0))$. So we may use \autoref{prop:functoriality_can_lift_P1_4pts} to lift $\tau_4$.
\end{proof}

Notice that the $\varphi$-exceptional curves that appear in $\Gamma'$ are the log canonical places of $(T,E+\Gamma)$, while the canonical places do not appear in $\Gamma'$. In particular, it is not sufficient to construct an equivariant lift of $(T',E'+\Gamma')$ in order to produce an equivariant lift of $\varphi$. For latter reference, we give these ``hidden" curves a name:

\begin{definition}\label{def:-2_curves}
An irreducible component $G$ of $\Exc(\varphi)$ is called a \emph{canonical place} for $\varphi$ if the discrepancy $a(G; T, E+\Gamma)=0$.
Notice that the set of canonical places of $\varphi$ is stable under the action of $\mu_{2,k}$. 
\end{definition}

Next we explain how to obtain a lifting of $\varphi$ verifying the hypothesis of \autoref{prop:existence_of_log_liftings}.

\begin{proposition}\label{prop:lift_equiv_MMP}
Let $(\mathbf{U},\mathbf{E}_\mathbf{U}+
\mathbf{\Gamma}_\mathbf{U})$ be a $\mu_{m,W(k)}$-equivariant projective lift of $(U,E_U+\Gamma_U)$ over $W(k)$. Suppose in addition that:
    \begin{enumerate}
        \item $\mathbf{E}_\mathbf{U}$ is the canonical lift of $E_U$;
        \item the $\mu_{2,W(k)}$-action on $\mathbf{E}_\mathbf{U}$ is the canonical lift of the $\mu_{2,k}$-action on $E_U$;
        \item if $G$ is sum of the canonical places of $\varphi$, then $h_*G$ lifts equivariantly to a divisor on $\mathbf{U}$ \emph{(\footnote{
            In accordance with \autoref{def:equiv_lifting}, we require that specialization over $k$ induces a bijection between the irreducible components of the lift and those of $h_*G$.
        })}.
    \end{enumerate}
Then there exists a $\mu_{2,W(k)}$-equivariant lift $\varphi_{W(k)}\colon (\mathbf{T}',\mathbf{E}'+\mathbf{\Gamma}')\to (\mathbf{T},\mathbf{E}+\mathbf{\Gamma})$ of $\varphi$ 
satisfying the hypothesis of \autoref{prop:existence_of_log_liftings}.
\end{proposition}
\begin{proof}
We divide the proof in three steps.

\medskip
\textsc{Step 1: Equivariant lifting of $h$.} 
Let us decompose $h$ as:
        $$T'=U_N\overset{g_N}{\longrightarrow} 
        \dots 
        \overset{g_2}{\longrightarrow}
        U_1\overset{g_1}{\longrightarrow}
        U_0=U,$$
where each $g_i\colon U_i\to U_{i-1}$ is the blow-up of the $\mu_{2}(k)$-orbit of a closed point $u_{i-1}$. Since $h^*(K_U+E_U+\Gamma_U)=K_{T'}+E'+\Gamma'$, we see that if we write
        $$K_{U_i}+E_{U_i}+\Gamma_{U_i}=(g_1\circ \dots\circ g_i)^*
        (K_U+E_U+\Gamma_U)$$
where $E_{U_i}$ is the strict transform of $E_U$ (equivalently the push-forward of $E'$), then:
    \begin{itemize}
        \item we have $\Gamma_{U_i}\geq 0$,
        \item the orbit of $u_{i}$ belongs to $\Supp(E_{U_i}+\Gamma_{U_i})$ for every $i$, 
        \item we have $K_{U_i}+E_{U_i}+\Gamma_{U_i}=(g_{i+1}\circ \dots\circ g_N)_*(K_{T'}+E'+\Gamma')\sim 0$,
        \item as $h$ is a sequence of blow-down of $(-1)$-curves, every such curve intersects the pushforward of $E'+\Gamma'$ with multiplicity $1$, and so each pair $(U_i,E_{U_i}+\Gamma_{U_i})$ is snc.
    \end{itemize}
We also let $G_{U_i}=(g_{i+1}\circ\dots\circ g_N)_*G$: it is a $\mu_{2,k}$-invariant divisor on $U_i$.

We will construct, by induction on $i$, an equivariant crepant morphism 
    \begin{equation}\label{eqn:lifting_equiv_MMP}
    (\mathbf{U}_i,\mathbf{E}_{\mathbf{U}_i}
    +\mathbf{\Gamma}_{\mathbf{U}_i})
    \longrightarrow (\mathbf{U},\mathbf{E}_\mathbf{U}+\mathbf{\Gamma}_\mathbf{U})
    \end{equation}
as a composition of blow-ups of $W(k)$-points, which is an equivariant lift of 
    $$g_1\circ\dots\circ g_i\colon (U_i,E_{U_i}+\Gamma_{U_i})\longrightarrow (U,E_U+\Gamma_U),$$ 
and such that we have a $\mu_{2,W(k)}$-invariant relative divisor $\mathbf{G}_{\mathbf{U}_i}$ lifting $G_{U_i}$. Notice that in this situation the restricted morphism $\mathbf{E}_{\mathbf{U}_i}\to \mathbf{E}_\mathbf{U}$ is a $\mu_{2,W(k)}$-equivariant morphism. Moreover, since $\mathbf{U}$ is smooth over $W(k)$, every $\mathbf{U}_i$ will also be smooth over $W(k)$.

In case $i=0$, the morphism \autoref{eqn:lifting_equiv_MMP} exists by hypothesis. Assume that the morphism \autoref{eqn:lifting_equiv_MMP} is constructed for some $i\geq 0$. We lift the point $u_i\in U_i$ to $\mathbf{U}_i$ as follows:
    \begin{enumerate}
        \item Suppose that the action of $\mu_{2,k}$ at $u_i$ is not free. Then by \autoref{lemma:fixed_locus_is_smooth} there exists a unique $W(k)$-point $\mathbf{u}_i\in \mathbf{U}_i$ specializing to $u_i$ at which the $\mu_{m,W(k)}$-action is not free. 
        
        Notice that since $\mathbf{E}_{\mathbf{U}_i}+\mathbf{\Gamma}_{\mathbf{U}_i}$ is relatively snc and $\mu_{2,W(k)}$-invariant. So by applying \autoref{lemma:fixed_locus_is_smooth} to $\mathbf{U}_i$ and to the components of $\mathbf{E}_{\mathbf{U}_i}+\mathbf{\Gamma}_{\mathbf{U}_i}$ whose specializations contain $u_i$, we see that $\bold{u}_i$ belongs to $\mathbf{E}_{\mathbf{U}_i}+\mathbf{\Gamma}_{\mathbf{U}_i}$. 
        
        \item If $u_i$ is a $0$-stratum of $E_{U_i}+\Gamma_{U_i}$ and the action at $u_i$ is free, we let $\mathbf{u}_i$ be the unique relative $0$-stratum of $\mathbf{E}_{\mathbf{U}_i}+\mathbf{\Gamma}_{\mathbf{U}_i}$ which specializes to $u_i$. Since $\mathbf{E}_{\mathbf{U}_i}+\mathbf{\Gamma}_{\mathbf{U}_i}$ is relatively snc, we see that $\mathbf{u}_i$ is indeed a $W(k)$-point.

        \item If $u_i$ is not a $0$-stratum of $E_{U_i}+\Gamma_{U_i}$ and if the $\mu_{2,k}$-action is free at $u_i$, then: we let $\mathbf{u}_i$ be an arbitrary $W(k)$-lift on $\mathbf{E}_{\mathbf{U}_i}+\mathbf{\Gamma}_{\mathbf{U}_i}$. 
        
    \end{enumerate}
Then we let $\mathbf{g}_i\colon \mathbf{U}_{i+1}\to \mathbf{U}_i$ be the blow-up of the $\mu_{2,W(k)}$-orbit of $\mathbf{u}_i$: this is clearly an equivariant lift of $g_i$. We define $\mathbf{E}_{\mathbf{U}_{i+1}}$ to be the strict transform of $\mathbf{E}_{\mathbf{U}_i}$, and $\mathbf{\Gamma}_{\mathbf{U}_{i+1}}$ by the formula
        $$K_{\mathbf{U}_{i+1}}+\mathbf{E}_{\mathbf{U}_{i+1}}
        +\mathbf{\Gamma}_{\mathbf{U}_{i+1}}
        = \mathbf{g}_{i+1}^*(
        K_{\mathbf{U}_{i}}+\mathbf{E}_{\mathbf{U}_{i}}
        +\mathbf{\Gamma}_{\mathbf{U}_{i}}).$$
Since the $W(k)$-points we blow-up is contained in $\mathbf{E}_{\mathbf{U}_{i}}+\mathbf{\Gamma}_{\mathbf{U}_{i}}$, we have $\mathbf{\Gamma}_{\mathbf{U}_{i+1}}\geq 0$. Finally, we let $\mathbf{G}_{\mathbf{U}_{i+1}}$ be the strict transform of $\mathbf{G}_{\mathbf{U}_{i}}$ plus the $\mathbf{g}_{i+1}$-exceptional relative curves which specialize to exceptional $g_{i+1}$-curves.




After finitely many steps, we obtain an equivariant lift $\mathbf{h}\colon (\mathbf{T}',\mathbf{E}'+\mathbf{\Gamma}')\to (\mathbf{U},\mathbf{E}_\mathbf{U}+\mathbf{\Gamma}_\mathbf{U})$ of $h$, together with an equivariant lift $\mathbf{G}$ of the canonical places of $\varphi$ on $\mathbf{T}'$. By construction $\mathbf{T}'$ is projective over $W(k)$.

\medskip
\textsc{Step 2: Equivariant lifting of $\varphi$.} We can factor $\varphi$ into $\mu_{m,k}$-equivariant morphisms
        $$T'\overset{\eta}{\longrightarrow} T'' \overset{\theta}{\longrightarrow} T$$
where:
    \begin{itemize}
        \item $\theta$ is a dlt modification of $(T,E+\Gamma)$ obtained as the minimal resolution of the strictly log canonical singularities of $(T,E+\Gamma)$; we let $E''=\theta^{-1}_*E$ and $\Gamma''\geq 0$ be defined by $K_{T''}+E''+\Gamma''=\theta^*(K_T+E+\Gamma)$;
        \item $\eta$ is  the the minimal resolution of the canonical singularities of $(T'',E''+\Gamma'')$.
    \end{itemize}
Notice that $\Exc(\eta)$ is the sum $G$ of the canonical places of $\varphi$, and that $K_{T'}\sim -E'-\Gamma'$ forces 
    $$H^2(T',\sO_{T'})=H^0(T',\sO(K_{T'}))^\vee=0.$$
Since $\varphi_*^{-1}\Gamma\leq \Gamma'$ and since $(\mathbf{T}',\mathbf{E}'+\mathbf{\Gamma}'+\mathbf{G})$ is a $\mu_{2,W(k)}$-equivariant lift of $(T',E'+\Gamma'+G)$, we see that we have an equivariant lift of $(T',E'+\varphi_*^{-1}\Gamma+G)$ over $W(k)$. Therefore we apply \cite[Proposition 5.2]{BBKW24} and obtain a projective $\mathbf{T}''$ and a birational morphism $\eta_{W(k)}\colon \mathbf{T}'\to \mathbf{T}''$ lifting $\eta$ and contracting exactly $\mathbf{G}$. Clearly the $\mu_{2,W(k)}$-action descends to $\mathbf{T}''$ and makes $\eta_{W(k)}$ equivariant. We let $\mathbf{E}''=\eta_{W(k),*}\mathbf{E}'$ and $\mathbf{\Gamma}''=\eta_{W(k),*}\mathbf{\Gamma}'$. By construction, the divisors $\mathbf{E}''$ and $\mathbf{\Gamma}''$ are $\mu_{2,W(k)}$-invariant.

If $(T,E+\Gamma)$ is already dlt then 
    $$\eta_{W(k)}\colon 
        (\mathbf{T}',\mathbf{E}'+\mathbf{\Gamma}')\longrightarrow
        (\mathbf{T}'',\mathbf{E}''+\mathbf{\Gamma}'')$$
is the desired equivariant lifting of $\varphi\colon (T',E'+\Gamma')\to (T,E+\Gamma)$.
 
Assume that $(T,E+\Gamma)$ is not dlt, so $\theta$ is not an isomorphism. Let $\Sigma$ be the sub-divisor of $\Gamma''$ consisting of the $\theta$-exceptional curves, and le $\mathbf{\Sigma}$ be the sub-divisor of $\mathbf{\Gamma}''$ lifting $\Sigma$. Both divisors are invariant for the group action. 

Let $A$ be a very ample divisor on $T$, and let $A_{T''}=\theta^*A$. Then $A_{T''}$ is semi-ample and a suitable multiple defines the morphism $\theta$. Since $K_{T''}+\Sigma\sim -E''-(\Gamma''-\Sigma)$ we have
        $$H^2(T'', \sO(-\Sigma))=H^0(T'',\sO(K_{T''}+\Sigma))^\vee =0.$$
It follows from \cite[Corollary 2.18]{BBKW24} that we can lift $\sA_{T''}=\sO_{T''}(A_{T''})$ to a line bundle $\sA_{\mathbf{T}''}$ on $\mathbf{T}''$ with the property that $\sA_{\mathbf{T}''}|_{\mathbf{\Sigma}}\cong \sO_{\mathbf{\Sigma}}$. 
By \cite[Proposition 5.4]{BBKW24} the line bundle $\sA_{\mathbf{T}''}$ is semi-ample and induces a projective birational morphism $\theta_{W(k)}\colon \mathbf{T}''\to \mathbf{T}$ lifting $\theta$ and contracting exactly $\mathbf{\Sigma}$. Clearly the $\mu_{2,W(k)}$-action descends to $\mathbf{T}$. We let $\mathbf{E}=\theta_{W(k),*}\mathbf{E}''$ and $\mathbf{\Gamma}=\theta_{W(k),*}\mathbf{\Gamma}''$. Then 
        $$\theta_{W(k)}\circ \eta_{W(k)}\colon 
        (\mathbf{T}',\mathbf{E}'+\mathbf{\Gamma}')\longrightarrow
        (\mathbf{T},\mathbf{E}+\mathbf{\Gamma})$$
is the desired equivariant lifting of $\varphi\colon (T',E'+\Gamma')\to (T,E+\Gamma)$.

\medskip
\textsc{Step 3: Conclusion of the proof.} We have constructed an equivariant lifting $\varphi_{W(k)}$ of $\varphi$. By construction $\mathbf{T}'$ is projective over $W(k)$, and the birational map $\mathbf{h}\circ \varphi_{W(k)}^{-1}$ restricts to a $\mu_{2,W(k)}$-equivariant isomorphism $\mathbf{E}\cong \mathbf{E}_\mathbf{U}$. Therefore the hypothesis of \autoref{prop:existence_of_log_liftings} are satisfied.
\end{proof}

We are reduced to construct projective liftings of $(U,E_U+\Gamma_U)$ satisfying the hypothesis of \autoref{prop:lift_equiv_MMP} one connected component at the time. Clearly we can proceed one connected component at a time. We will now consider all the cases obtained in \autoref{section:classification_min_models} and show that we can satisfy the hypothesis of \autoref{prop:lift_equiv_MMP} in each one. 

\begin{remark}
In the next sections, we will take $\mu_2$-quotients at several places: we will tacitly use \autoref{prop:descent_projectivity}, which guarantees that when doing so we stay in the category of projective varieties. 
\end{remark}


\subsection{Equivariant liftings: case of del Pezzo surfaces}\label{section:equiv_lift_dP}
We start with the cases of \autoref{lemma:min_equiv_model_rational_case} where $U$ is a del Pezzo surface. In this context, recall that $m=2$ and that $E_U$ is an irreducible regular ordinary curve of genus $1$. We denote by $\sigma_U$ the non-trivial involution of $U$ given by the $\mu_2$-action, and by $\sigma_{E_U}$ its restriction to $E_U$.

\subsubsection{Case of $\bP^1\times\bP^1$}\label{section:equiv_lift_P1xP1}
Suppose that $U=\bP^1\times\bP^1$ and that the involution $\sigma_U$ is given by the switch of the two factors. We record the following easy lemma:

\begin{lemma}\label{lemma:inv_map_to_P1xP1}
Let $A$ be a Noetherian ring and $X$ be an $A$-scheme. Then a proper $A$-embedding $X\hookrightarrow \bP^1_A\times \bP^1_A$ whose image is invariant under the switch of factors $\varrho$ on $\bP^1_A\times\bP^1_A$, is equivalent to the data of:
    \begin{itemize}
        \item two line bundles $L_i\in \Pic(X)$ equipped with two-dimensional sub-linear systems $\Lambda_i\subset H^0(X,L_i)$ ($i=1,2$), 
        \item an involution $\sigma_X\in \Aut_A(X)$, and
        \item an isomorphism $\sigma_X^*(L_2)\cong L_1$ which takes $\sigma_X^*(\Lambda_2)$ to $\Lambda_1$.
    \end{itemize}
\end{lemma}
\begin{proof}
The bijection sends $\iota\colon X\hookrightarrow \bP^1\times \bP^1$ to the two-dimensional $A$-modules associated to $\pr_i\circ \iota\colon X\to \bP^1_A$ and to the restricted involution $\sigma_X=\varrho|_X$.
\end{proof}

In particular, the embedding $g\colon E_U\hookrightarrow \bP^1\times\bP^1$ is given two line bundles $L_i\in \Pic(E_U)$, two sub-linear systems $\Lambda_i\subset H^0(E_U,L_i)$ ($i=1,2$) and by an isomorphism $\alpha\colon L_1\cong \sigma_{E_U}^*(L_2)$ such that $\alpha(\Lambda_1)=\sigma_{E_U}^*(\Lambda_2)$.

Let $\mathbf{E}_\mathbf{U}$ be the canonical lifting of $E_U$ over $W(k)$. Let $\sL_i$be the canonical lift of $L_i$ (for $i=1,2$), and $\sigma_{\mathbf{E}_\mathbf{U}}$ be the canonical lift of $\sigma_{E_U}$. The isomorphism $\alpha$ lifts to $\alpha_{W(k)}\colon \sL_1\cong \sigma_{\mathbf{E}_\mathbf{U}}^*(\sL_2)$ by \cite[Appendix, Proposition 2]{MS87}.

\begin{lemma}\label{lemma:lift_switch_inv_map}
For each $i\in\{1,2\}$ there is a two-dimensional free sub-$W(k)$-module $\mathbf{\Lambda}_i\subset H^0(\mathbf{E}_\mathbf{U},\sL_i)$ lifting $\Lambda_i$, such that $\alpha_{W(k)}(\mathbf{\Lambda}_1)=\sigma_{\mathbf{E}_\mathbf{U}}^*(\mathbf{\Lambda}_2)$.
\end{lemma}
\begin{proof}
Let $\mathbf{V}_i=H^0(\mathbf{E}_\mathbf{U},\sL_i)$. Define a $W(k)$-linear involution on the free $W(k)$-module $\mathbf{V}=\mathbf{V}_1\oplus \mathbf{V}_2$ as follows:
        $$(\mathbf{t}_1,\mathbf{t}_2)\mapsto 
        \left(
        \alpha_{W(k)}^{-1}(\sigma_{\mathbf{E}_\mathbf{U}}^*(\mathbf{t}_2)), \ 
        \sigma_{\mathbf{E}_\mathbf{U}}^{-1,*}(\alpha_{W(k)}(\mathbf{t}_1)
        \right).$$
Let $V=\mathbf{V}\otimes_{W(k)}k=\bigoplus_{i=1}^2H^0(E_U,L_i)$. Since $\mathbf{V}$ is a torsion-free $W(k)$-module, we have an exact sequence
        $$0\to \mathbf{V}\overset{p}{\longrightarrow} \mathbf{V}\overset{\res}{\longrightarrow} V\to 0.$$
This sequence is $\mu_{2,W(k)}$-equivariant: the corresponding $\mu_{2,k}$-action on $V$ is given by
        $$(t_1,t_2)\mapsto \left(
        \alpha^{-1}(\sigma_{E_U}^*(t_2)), \ 
        \sigma^{-1,*}_{E_U}(\alpha(t_1))
        \right).$$
Now let $x_2,y_2\in \Lambda_2$ be a basis and define $x_1=\alpha^{-1}(\sigma_{E_U}^*(x_2))$, $y_1=\alpha^{-1}(\sigma_{E_U}^*(y_2))$: by assumption $x_1,y_1$ give a basis of $\Lambda_1$. Both pairs $(x_1,x_2)$ and $(y_1,y_2)$ belong to the fixed sub-space $V^{\mu_2}$. 

Take the $\mu_{2,W(k)}$-invariants of the above exact sequence: this gives a long exact sequence of group cohomology with respect to $\mu_{2,W(k)}$. Since $\mu_{2,W(k)}$ is linearly reductive, all the higher cohomology groups vanish \cite[I, Lemma 4.3]{Jantzen_Representations_alg_grps}. So the restriction modulo $p$ map $\mathbf{V}^{\mu_{2, W(k)}}\to V^{\mu_{2,k}}$ is surjective. Hence we can find $\mathbf{x}_1,\mathbf{y}_1\in \mathbf{V}_1$ and $\mathbf{x}_2,\mathbf{y}_2\in\mathbf{V}_2$ such that $(\mathbf{x}_1,\mathbf{x}_2)$ and $(\mathbf{x}_2,\mathbf{y}_2)$ are $\mu_{2,W(k)}$-fixed and lift respectively $(x_1,x_2)$ and $(y_1,y_2)$.
We let $\mathbf{\Lambda}_i$ be the free $W(k)$-module generated by $\mathbf{x}_i$ and $\mathbf{y}_i$: by construction $\mathbf{\Lambda}_1$ and $\mathbf{\Lambda}_2$ satisfy the desired properties.
\end{proof}

The sub-modules $\mathbf{\Lambda}_i\subset H^0(\mathbf{E}_\mathbf{U},\sL_i)$ define a projective morphism $\mathbf{g}\colon \mathbf{E}_\mathbf{U}\to \bP^1_{W(k)}\times\bP^1_{W(k)}$ whose reduction modulo $p$ is the closed embedding $g\colon E_U\hookrightarrow \bP^1\times \bP^1$: this forces $\mathbf{g}$ to be a closed embedding. By \autoref{lemma:inv_map_to_P1xP1} the image of $\mathbf{g}$ is invariant under the switch of factors, which restricts to $\sigma_{\mathbf{E}_\mathbf{U}}$ on $\mathbf{E}_\mathbf{U}$.

To satisfy the hypothesis of \autoref{prop:lift_equiv_MMP}, it remains to lift equivariantly some $\sigma_U$-invariant curves on $\bP^1_k\times \bP^1_k$ to $\bP^1_{W(k)}\times \bP^1_{W(k)}$. We have a cartesian diagram
        $$\begin{tikzcd}
        \bP^1_k\times \bP^1_k \arrow[r, hook] \arrow[d, "q"] &
        \bP^1_{W(k)}\times \bP^1_{W(k)} \arrow[d, "\mathbf{q}"] \\
        \bP^2_k \arrow[r, hook] & \bP^2_{W(k)}
        \end{tikzcd}$$
where $q$ (resp.\ $\mathbf{q}$) is the geometric quotient by the $\mu_{2}$-action induced by the switch of factors. Let $G$ be an $\sigma_U$-invariant curve, with image $C=q(G)$. We have $q^{-1}(C)=G$, see the proof of \autoref{prop:sing_qt_pair}. We can lift the couple $(\bP^2_k,C)$ to a couple $(\bP^2_{W(k)},\mathbf{C})$, in the sense of \autoref{def:lifting}, and we let $\mathbf{G}=\mathbf{q}^{-1}(\mathbf{C})$. Then $\mathbf{G}$ lifts $G$ (such that specialization over $k$ induces a bijection of irreducible components). So we are done.

\subsubsection{Case of $\bP^2$}\label{section:case_of_dP:P2}
Suppose that $U=\bP^2$. Then $E_U$ is a regular curve of degree $3$. Recall that an \emph{inflection point} of $E_U\subset \bP^2$ is a point $e\in E_U$ such that there exists a line in $|\sO_{\bP^2}(1)|$ which meets $E_U$ at $e$ with multiplicity $3$.

\begin{lemma}
There exists an inflection point $e$ of $E_U\subset \bP^2$ which is fixed by $\sigma_{E_U}$.
\end{lemma}
\begin{proof}
The action of $\sigma_U$ preserves $E_U$ and sends lines to lines, hence its restriction $\sigma_{E_U}$ sends inflection points to inflection points. There is an odd number of inflection points (\footnote{
    Fix one inflection point $e'$ of $E_U\subset \bP^2$. Then the set of inflection points is equal to the set of $3$-torsion points of the elliptic curve $(E_U,e')$. If the characteristic of $k$ is different from $3$ we have $|E[3]|=9$. If the characteristic of $k$ is $3$, as $E_U$ is ordinary we have $|E_U[3]|=3$.
}): since $\sigma_{E_U}$ is an involution, it must fix at least one of them.
\end{proof}

Let $e\in E_U$ be the inflection point given by the above lemma. We think of the given embedding $\iota\colon E_U\hookrightarrow \bP^2$ as the one given by the very ample line bundle $\sO(3e)$ on $E_U$.
The action of $\sigma_{E_U}$ on $\sO(3e)$ induces an involution $\sigma_{U}'$ of $\bP^2=\bP\left( H^0(E_U,\sO(3e))^\vee \right)$, with the property that $\sigma_{U}|_{E_U}=\sigma_{E_U}$. We claim that $\sigma_{U}=\sigma_{U}'$: this follows from the following lemma applied to $\sigma_{U}'\circ\sigma_{U}^{-1}$.

\begin{lemma}
Let $\varphi$ be an automorphism of $\bP^2$ which restricts to the identity on a reduced curve $C\in |\sO_{\bP^2}(3)|$. Then $\varphi$ is trivial.
\end{lemma}
\begin{proof}
Pick a general line $L\subset \bP^2$: it intersects $C$ in three distinct points. Since $\varphi|_C$ is trivial, each of these three points is fixed by $\varphi$. Since $\varphi$ sends lines to lines, we deduce that $\varphi(L)=L$. An automorphism of $L\cong \bP^1$ with three fixed points is the identity, so we find that $\varphi$ is the identity along $L$. As $L$ is general, we deduce that $\varphi$ must be trivial.
\end{proof}

Let $\mathbf{E}_\mathbf{U}$ be the canonical lifting of the ordinary curve $E_U$. Let $\sigma_{\mathbf{E}_\mathbf{U}}$ be the canonical lift of $\sigma_{E_U}$: it defines a $\mu_{2,W(k)}$-action on $\mathbf{E}_\mathbf{U}$. By \autoref{lemma:fixed_locus_is_smooth} the fixed locus of $\sigma_{\mathbf{E}_\mathbf{U}}$ is smooth over $W(k)$, and therefore the point $e$ lifts uniquely to a $\mu_{2,W(k)}$-fixed $W(k)$-point $\mathbf{e}$ of $\sigma_{\mathbf{E}_\mathbf{U}}$. Let 
    $$\mathbf{U}=
    \bP_{W(k)}\left(H^0(\mathbf{E}_\mathbf{U},\sO(3\mathbf{e}))^\vee\right)
    \cong \bP^2_{W(k)}.$$
The action of $\sigma_{\mathbf{E}_\mathbf{U}}$ on $\sO(3\mathbf{e})$ induces an involution $\sigma_\mathbf{U}$ of $\mathbf{U}$ lifting $\sigma_U$. Moreover, we have a closed embedding $\iota_{W(k)}\colon \mathbf{E}_\mathbf{U}\hookrightarrow \mathbf{U}$ which is a $\mu_{2,W(k)}$-equivariant lift of $\iota$. 

\begin{lemma}
Let $G$ be a $\sigma_U$-invariant divisor on $\bP^2_k$. Then there exists a $\sigma_\mathbf{U}$-invariant divisor on $\bP^2_{W(k)}$ lifting $G$.
\end{lemma}
\begin{proof}
By assumption we have $G\in H^0(\bP^2_k,\sO(n))^{\mu_{2,k}}$ for some $n>0$. We have an exact sequence
        $$0\to H^0(\bP^2_{W(k)},\sO(n))\overset{p}{\longrightarrow}
        H^0(\bP^2_{W(k)},\sO(n))
        \overset{\res}{\longrightarrow}
        H^0(\bP^2_k,\sO(n))\to 0.$$
Taking $\mu_{2,W(k)}$-invariants, we see as in the proof of \autoref{lemma:lift_switch_inv_map} that there exists a section $\mathbf{G}$ of $H^0(\bP^2_{W(k)},\sO(n))^{\mu_{2,W(k)}}$ that lifts $G$. By taking $\mathbf{G}$ general, we can assume that specialization over $k$ induces a bijection between the irreducible components of $\mathbf{G}$ and those of $G$.
\end{proof} 

Therefore the hypothesis of \autoref{prop:lift_equiv_MMP} are satisfied, and we move on to the next case.

\subsubsection{Case of the Geiser involution}\label{section:Geiser_inv}
Suppose that $U$ is a del Pezzo surface of degree $2$ and that the $\mu_2$-action is the Geiser involution, whose construction we recall now. We have an embedding
        $$U=V(w^2-f_4(x,y,z))\hookrightarrow \bP(1,1,1,2)=\bP$$
where $(x,y,z,w)$ are the coordinates of the weighted projective space $\bP$ with respective degrees $(1,1,1,2)$ and $f_4$ an homogeneous polynomial of degree $4$. The projection
        $$\bP\dashrightarrow \bP^2, \quad 
        [x:y:z:w] \mapsto [x:y:z]$$
restricts to a double covering $U\to \bP^2$ whose non-trivial deck transformation is the Geiser involution of $U$. In other words, it is given by
        $$U\to U, \quad [x:y:z:w]\mapsto [x:y:z:-w].$$
The fixed locus of this involution is $U \cap V_{\bP}(w)$, 
and that the branch locus of $U\to \bP^2$ is the divisor $G=V_{\bP^2}(f_4)$. The image $L$ of $E_U$ in $\bP^2$ is the geometric $\mu_2$-quotient of $E_U$ by \autoref{prop:restriction_linearly_red_qt}, so $L\cong \bP^1$. It which meets $G$ at four points $q_1,\dots,q_4$, which are the images of the four fixed points of the $\mu_2$-action on $E_U$. By Bézout theorem we deduce that $L\in |\sO_{\bP^2}(1)|$. By \autoref{rmk:double_cover_of_4A1_curve} the pair $(L,\frac{1}{2}\sum_{i=1}^4q_i)$ is globally sharply $F$-split.

We lift $L$ to a line $\mathbf{L}$ on $\bP^2_{W(k)}$. Let $\mathbf{q}_1,\dots,\mathbf{q}_4\in \mathbf{L}$ be the canonical lifts over $W(k)$ of the four intersection points of $L$ and $G$ (cf.\ \autoref{def:canonical_lift_P1_4pts}). By \autoref{lemma:lifts_divisor_with_correct_intersections} below we can find $\mathbf{G}\in H^0(\bP^2_{W(k)},\sO(4))$ which lifts $G$, and such that $\mathbf{G}$ meets $\mathbf{L}$ properly at $\mathbf{q}_1,\dots,\mathbf{q}_4$.

\begin{lemma}\label{lemma:lifts_divisor_with_correct_intersections}
Let $\mathbf{r}_1,\dots,\mathbf{r}_4$ be $W(k)$-points of $\bP^2_{W(k)}$ whose reduction $r_1,\dots,r_4$ in $\bP^2_k$ are pairwise distinct. Let $F\in H^0(\bP^2_k,\sO(n))$ be such that $F(r_i)=0$ for each $i=1,\dots,4$. Then there exists $\mathbf{F}\in H^0(\bP^2_{W(k)},\sO(n))$ lifting $F$ such that $\mathbf{F}(\mathbf{r}_i)=0$ for every $i=1,\dots,4$.
\end{lemma}
\begin{proof}
Pick any $\mathbf{F}\in H^0(\bP^2_{W(k)},\sO(n))$ lifting $F$. It is sufficient to find an $\mathbf{H}\in H^0(\bP^2_{W(k)},\sO(n))$ such that
        $$\mathbf{F}(\mathbf{r}_i)+p\mathbf{H}(\mathbf{r}_i)=0,
        \quad i=1,\dots,4.$$
By hypothesis we have $\mathbf{F}(\mathbf{r}_i)=p\mathbf{a}_i$ where $\mathbf{a}_i\in W(k)$, for every $i$. Thus we reduce to solve the system of equations
        $$\mathbf{H}(\mathbf{r}_i)=-\mathbf{a}_i,\quad i=1,\dots,4.$$
The variables are the coefficients of the homogeneous polynomial $\mathbf{H}$. In particular, this is a system of linear equations. Thus by Hensel's lemma, it is enough to show that a solution exists modulo $p$. So we reduce to find $H\in H^0(\bP^2_k,\sO(n))$ such that 
        $$H(r_i)=-a_i, \quad i=1,\dots,4$$
where $a_i$ is the reduction of $\mathbf{a}_i$ in $k$. 

Choose coordinates $x,y,z$ on $\bP^2_k$ such that $r_1=[1:0:0]$, $r_2=[0:1:0]$, $r_3=[0:0:1]$. Let $r_4=[\alpha:\beta:\gamma]$ be the coordinates of the fourth point. Writing $H=\sum_{i+j+l=n}u_{ijl}x^iy^jz^l$, we want to solve
        $$\left\{
        \begin{array}{lcllllcl}
             u_{n00} &=& -a_1, & & &
             u_{0n0} &=& -a_2, \\
             u_{00n} &=& -a_3, &&&
             \sum_{i+j+l=n}u_{ijl}\cdot
             \alpha^i\beta^j\gamma^l &=& -a_4.
        \end{array}
        \right.
        $$
By hypothesis at least two coordinates of $r_4$ are non-zero, and so the last equation does not involve only $u_{n00}$, $u_{0n0}$ and $u_{00n}$. Therefore the system has a solution.
\end{proof}

Now let $\mathbf{q}\colon \mathbf{U}\to \bP^2_{W(k)}$ be the cyclic cover of degree $2$ branched above the divisor $\mathbf{G}$. By \autoref{lemma:restriction_cyclic_cover} the reduction modulo $p$ of $\mathbf{U}\to \bP^2_{W(k)}$ is the double cover $U\to \bP^2$. By \autoref{lemma:restriction_cyclic_cover} and \autoref{rmk:can_cyclic_cover}, the preimage $\mathbf{E}_\mathbf{U}$ of $\mathbf{L}$ in $\mathbf{U}$ is the canonical lift of the ordinary elliptic curve $E_U$ over $W(k)$, and the induced $\mu_{2,W(k)}$-action on $\mathbf{E}_\mathbf{U}$ is the canonical lift of the $\mu_{2,k}$-action on $E_U$.

Let $G$ be an effective $\mu_2$-invariant divisor on $U$, with image $q(G)=M$ on $\bP^2$. By decomposing $G$ into minimal $\mu_2$-invariant divisors, we may assume that $M$ is irreducible. Let $\mathbf{M}$ be a divisor on $\bP^2_{W(k)}$ lifting $M$, and let $\mathbf{G}=\mathbf{q}^{-1}(\mathbf{M})$. Then $\mathbf{G}$ is a $\mu_{2,W(k)}$-invariant divisor lifting $G$. Thus the hypothesis of \autoref{prop:lift_equiv_MMP} are all satisfied.

\subsubsection{Case of the Bertini involution}\label{section:Bertini}
Suppose that $U$ is a del Pezzo surface of degree $1$ and that the $\mu_2$-action is the Bertini involution. We have an embedding
        $$U=V(w^2-f_6(x,y,z))\hookrightarrow \bP(1,1,2,3)=\bP$$
where the homogeneous coordinates $(x,y,z,w)$ have degrees $(1,1,2,3)$, and where $f_6$ is an homogeneous polynomial of degree $6$. The projection
        $$\bP\dashrightarrow \bP(1,1,2)=Q,\quad 
        [x:y:z:w]\mapsto [x:y:z]$$
restricts to a double covering $q\colon U\to Q$, whose non-trivial deck transformation is the Bertini involution of $U$. The divisor $G=V_Q(f_6)$ belongs to the branch locus of $q$. 

The surface $Q$ is the quadric cone, with a single singular point that we call its vertex $\upsilon$. Since $Q$ is not smooth over $k$ while its cover $U$ is is, we deduce with the help of \autoref{lemma:sing_qt} that the branch locus of $q$ is the disjoint union of $G$ and $\upsilon$ (so $\upsilon\notin G$).

The preimage of the vertex $\upsilon$ in $U$ is the base locus of $|-K_U|$. In particular, it belongs to $E_U$. So $q(E_U)$ is an irreducible curve which contains the vertex and which is regular away from it (by \autoref{prop:restriction_linearly_red_qt}). The Riemann--Hurwitz formula gives $\sO_U(-K_U)= q^*\sO_{\mathbb{P}(1,1,2)}(1)$, and so we see that $L=q(E_U)$ is a standard line going through $\upsilon$.

We think of the quadric cone $Q$ as embedded into $\bP^3$ by the complete linear system $|\sO_Q(2)|$. Then $G\in |\sO_Q(6)|$ is the trace on $Q$ of a cubic surface $C\in |\sO_{\bP^3}(3)|$. The surface $C$ intersects $L$ in $3$ points, which are different from the vertex (since otherwise $\upsilon \in C\cap Q=G$).

Let $q_1,q_2,q_3$ be the intersection points of $L$ with $C$, and write $q_4=\upsilon$ as a point of $L$. Riemann--Hurwitz shows that the quotient $E_U\to L$ has four branching points: by what we have said, they must be $q_1,\dots,q_4$. Moreover, the pair $(L,\frac{1}{2}\sum_{i=1}^4q_i)$ is globally sharply $F$-split by \autoref{rmk:double_cover_of_4A1_curve}.

Let $\mathbf{Q}$ be the quadric cone over $W(k)$ with vertex $\mathbf{v}$, embedded via $\sO_{\mathbb{P}_{W(k)}(1,1,2)}(2)$ over $W(k)$. Let $(\bP^1_{W(k)},\mathbf{q}_1+\dots +\mathbf{q}_4)$ be the canonical lift of $(L,q_1+\dots+q_4)$ (cf.\ \autoref{def:canonical_lift_P1_4pts}). We can embed it as a standard line $\mathbf{L}\subset \mathbf{Q}$ such that $\mathbf{q}_4$ is sent to $\mathbf{v}$. By an argument similar to that of \autoref{lemma:lifts_divisor_with_correct_intersections}, we can lift $C$ to a cubic surface $\mathbf{C}\in |\sO_{\bP^3_{W(k)}}(3)|$ which intersect $\mathbf{L}$ at the points $\mathbf{q}_1,\mathbf{q}_2$ and $\mathbf{q}_3$.

Consider the reflexive sheaf $\sL=\sO_{\mathbf{Q}}(3)$: it is invertible away from $\mathbf{v}$, $2$-Cartier, and the divisor $\mathbf{G}=\mathbf{Q}\cap \mathbf{C}$ defines a global section of the reflexive power $\sL^{[2]}$. So we let
        $$\mathbf{U}=\Spec_{\mathbf{Q}} \left(
        \sO_{\mathbf{Q}} \oplus \sL^{[-1]}
        \right)
        \overset{\mathbf{q}}{\longrightarrow} 
        \mathbf{Q}$$
where the algebra structure is given by $\mathbf{G}$. We claim that $\mathbf{q}$ is a lift of $q\colon U\to Q$. By \autoref{lemma:restriction_cyclic_cover}, it is clear that 
    $$\left(\mathbf{q}\times_{\mathbf{Q}} \mathbf{Q}^*\right) \otimes k
    \cong q\times_Q Q^*,$$ 
where $\mathbf{Q}^*$ and $Q^*$ are respectively the complements of $\mathbf{v}$ and $\upsilon$. Since $\mathbf{Q}^*$ and $Q^*$ are two big open subsets, and since $\mathbf{q}$ and $q$ are finite, this isomorphism extends to $\mathbf{q}\otimes k\cong q$.

So $\mathbf{U}$ lifts $U$, and the deck involution of $\mathbf{q}$ lifts the Bertini involution. Moreover $\mathbf{E}_\mathbf{U}=\mathbf{q}^{-1}(\mathbf{L})$ is the canonical lift of $E_U$ and the restriction of the deck involution to $\mathbf{E}_\mathbf{U}$ is the canonical lift of the restriction of the Bertini involution to $E_U$ (this can be checked away from $\mathbf{v}$, and so it follows from \autoref{lemma:restriction_cyclic_cover} and \autoref{rmk:can_cyclic_cover}).

\medskip
To lift $\mu_2$-invariant effective divisors from $U$ to $\mathbf{U}$, it is sufficient to lift effective divisors from $Q$ to $\mathbf{Q}$ (cf.\ the argument at the end of \autoref{section:equiv_lift_P1xP1} and of \autoref{section:Geiser_inv}). We have a cartesian diagram
        $$\begin{tikzcd}
        \Sigma_2=\bP_{\bP^1_k}(\sO\oplus \sO(-2)) \arrow[r, hook] \arrow[d, "r"] 
        & \bP_{\bP^1_{W(k)}}(\sO\oplus \sO(-2)) \arrow[d, "\mathbf{r}"]=\mathbf{\Sigma}_2 \\
        Q\arrow[r, hook] & \mathbf{Q}
        \end{tikzcd}$$
where $r$ (resp.\ $\mathbf{r}$) is the contraction of the $(-2)$-section (resp.\ of the relative $(-2)$-section). Then it is sufficient to lift effective divisors from $\Sigma_2$ to $\mathbf{\Sigma}_2$. 

For latter use we show that it can be done on any Hirzebruch surface. Let $\Sigma_n=\bP_{\bP^1_k}(\sO\oplus \sO(-n))$ and fix a lift $\mathbf{\Sigma}_n=\bP_{\bP^1_{W(k)}}(\sO\oplus \sO(-n))$.

\begin{proposition}\label{prop:lift_div_on_Hirzebruch}
Any effective divisor of $\Sigma_n$ lifts to $\mathbf{\Sigma}_n$.
\end{proposition}
\begin{proof}
Fix a fiber $F$ of the projection $\Sigma_n\to \bP^1_k$, and let $s$ be the unique $(-n)$-section (or an arbitrary section of the projection if $n=0$). Clearly both divisors lift to $\mathbf{\Sigma}_n$, say to $\mathbf{F}$ and $\mathbf{s}$ respectively; we may take $\mathbf{F}$ to be a fiber above a $W(k)$-point of $\bP^1_{W(k)}$. Then 
    $$\mathbf{\Sigma}_n\setminus (\mathbf{F}\cup \mathbf{s})
    \cong \bA^2_{W(k)}.$$
Let $G$ be an irreducible divisor of $\Sigma_n$ which is not equal to $F$ nor $s$. Then $G$ intersects $\Sigma_n\setminus (F\cup s)\cong \bA^2_k$ non-trivially. On this open chart it is an hypersurface, and so it is easy to lift it to an irreducible divisor on $\mathbf{\Sigma}_n\setminus (\mathbf{F}\cup \mathbf{s})$. Its closure in $\mathbf{\Sigma}_n$ is an irreducible divisor which lifts $G$.
\end{proof}

We have thus proved that the hypothesis of \autoref{prop:lift_equiv_MMP} are all satisfied in case $U$ is a del Pezzo surface of degree $1$.

\subsection{Equivariant liftings: case of invariant conic bundles}\label{section:equiv_lift_conic_bd_I}
Next, we consider the cases where the minimal equivariant model $f\colon U\to B$ is an invariant conic bundle. According to \autoref{lemma:min_equiv_model_rational_case} and \autoref{lemma:min_equiv_model_elliptic_case} we have: $m=2$, the base $B\cong \bP^1$, the curve $E_U$ is regular and irreducible, and $\Gamma_U=0$.

We let $q\colon U\to V=U/\mu_2$ be the geometric quotient morphism. Let $C=q(E)$, and let $G$ be the branch divisor of $q$ on $V$. 

\begin{lemma}\label{lemma:qt_invar_conic_bdl}
We have a commutative diagram
        $$\begin{tikzcd}
        U \arrow[rr, "q"]\arrow[dr, "f" below left] && V \arrow[dl, "g"] \\
        & \bP^1 &
        \end{tikzcd}$$
and $g\colon V\to \bP^1$ is a $\bP^1$-bundle. Moreover:
    \begin{enumerate}
        \item $(V,C+G)$ is snc,
        \item $C$ is a section of $g$, 
        \item every component of $G$ dominates $\bP^1$ and $G\to \bP^1$ has degree $2$,
        \item $K_V+C+\frac{1}{2}G\sim_\bQ 0$.
    \end{enumerate}
\end{lemma}
\begin{proof}
Since $f$ is $\mu_2$-invariant, it factors through $V$. We have $\sO_{\bP^1}\subseteq g_*\sO_V\subseteq f_*\sO_U=\sO_{\bP^1}$ so $g$ is a fibration. We show that $g$ is a $\bP^1$-bundle. If $b\in \bP^1$ is a closed point, by \autoref{prop:restriction_linearly_red_qt} the fiber $V_b$ is given by the $\mu_2$-quotient of the fiber $U_b$. If $U_b\cong \bP^1$ then the quotient is isomorphic to $\bP^1$. 
If $U_b$ is a singular fiber, then: either it is the union of two $\bP^1$s meeting properly at one point and the $\mu_2$-action exchanges the two components as $\rho(V/\mathbb{P}^1)^G=1$, or it is a double line. In the former case, the quotient of $U_b$ is also isomorphic to $\bP^1$. We exclude the latter case as follows:
Suppose by contradiction $U_b=2L$ for a line $L$. If $S \subset U$ is a section of $f$, we have $1=S \cdot \mathfrak{f}=S\cdot 2L>1,$ contradiction.

So $g$ is a smooth morphism, and the surface $V$ is regular. \autoref{prop:sing_qt_pair} shows that $(V,C+G)$ is snc and that $G$ is of pure dimension one over $\bP^1$. The morphism $f|_{E_U}\colon E_U\to \bP^1$ is invariant of degree $2$, hence the induced morphism $g|_C\colon C\to \bP^1$ is of degree $1$: this shows that $C$ is a section of $g$. The Riemann--Hurwitz formula, together with the facts that $K_U+E_U\sim 0$ and $q^*C=E_U$ (cf.\ the proof of \autoref{prop:sing_qt_pair}), shows that $2(K_V+C+\frac{1}{2}G)\sim 0$. 

It remains to show that $g|_G\colon G\to \bP^1$ has degree $2$. The $\mu_2$-action on a general fiber of $f$ is non-trivial, and thus it has exactly two distinct fixed points, say $u$ and $u'$. Then $q(u)$ and $q(u')$ are the intersection points of $G$ with the corresponding fiber of $g$: this shows that $G$ has degree $2$ over $\bP^1$.
\end{proof}

The morphism $g$ exhibits $V$ as an Hirzebruch surface, say $V\cong \bP(\sO\oplus\sO(-n))$ over $\bP^1$. The N\'{e}ron--Severi group of $V$ is given by 
        $$\NS(V)=\mathbb{Z} \mathfrak{c}\oplus \mathbb{Z} \mathfrak{f}$$
where $\mathfrak{f}$ is the class of a fiber of $g$, and $\mathfrak{c}$ is the class of the section with self-intersection $-n$. Moreover, the effective cone is closed polyhedral, given by
        $$\overline{\Eff}(V)=\bR_{\geq 0}\mathfrak{c}\oplus \bR_{\geq 0}\mathfrak{f},$$
see e.g.\ \cite[\S 1.35]{Debarre_Higher_dimensional_algebraic_geometry}.

Let $q_1,\dots,q_4\in C$ be the four branching points of $q|_{E_U}\colon E_U\to C$. By \autoref{rmk:double_cover_of_4A1_curve}, the pair $(C,\frac{1}{2}\sum_{i=1}^4q_i)$ is globally sharply $F$-split. These four points are also the intersection points of $C$ and $G$, as $G$ is the image of the fixed locus of the $\mu_2$-action on $U$.

\begin{proposition}\label{prop:lifting_log_Hirzebruch_canonically}
There is a lift $(\mathbf{V},\mathbf{C}+\mathbf{G})$ of $(V,C+G)$ over $W(k)$ such that $(\mathbf{C},\mathbf{G}|_{\mathbf{C}})$ is the canonical lift of the pair $(C,q_1+\dots+q_4)$ (cf.\ \autoref{def:canonical_lift_P1_4pts}).
\end{proposition}
\begin{proof}
By \autoref{prop:lift_div_on_Hirzebruch} there is a lift of the snc pair $(V,C+G)$ over the Witt vectors: we have to prove that there exists one satisfying the statement of the proposition.

For each given deformation $(V_n,C_n+G_n)$ of that pair over $W_n(k)$, the set of further deformations over $W_{n+1}(k)$ is a torsor under $H^1(V,T_V(-\log(C+G))$ \cite[Proposition 8.4]{Kato_Log_smth_deformations_and_mod_of_curves}. At the same time, $(V_n,C_n+G_n)$ induces a deformation of $(C,C\cap G)$ over $W_n(k)$, and the relation between the further deformations over $W_{n+1}(k)$ is encoded by the restriction map
        $$H^1(V,T_V(-\log(C+G)))\longrightarrow 
        H^1(C,T_C(-\log(C\cap G)))$$
induced by the exact sequence \autoref{eqn:restr_for_twisted_log_vector_fields}. Since we have a canonical lifting of $(C,C\cap G)$---the one from \autoref{def:canonical_lift_P1_4pts}---, and since formal deformations of $V$ are algebraizable as $H^2(V,\sO_V)=0$ \cite[Corollary 8.5.6]{FAG}, the statement of the proposition follows if we can show that the above map of cohomology groups is surjective. Its cokernel embeds in the cohomology group
        $$H^2\left(V,T_V(-\log G)\otimes \sO(-C)\right)$$
and by Serre duality, it is sufficient to prove that
        $$H^0\left( V, \Omega^1_V(\log G)\otimes \sO(K_V+C)\right)=0.$$
Notice that:
    \begin{enumerate}
        \item We have $(-K_V-C)\cdot \mathfrak{f}=\frac{1}{2}G\cdot \mathfrak{f}>0$ as $G$ dominates $\bP^1$;
        \item Since $G$ has degree $2$ over $\bP^1$ and $C$ is a section of $g$, we have $G\equiv 2\mathfrak{c}+\alpha \mathfrak{f}$ and $C\equiv \mathfrak{c}+\beta\mathfrak{f}$ for some $\alpha,\beta\geq 0$.
        Thus
                $$-K_V-G\equiv C-\frac{1}{2}G \equiv \left(\beta-\frac{\alpha}{2}\right)\mathfrak{f}$$
        is $g$-trivial, in particular $g$-nef.
    \end{enumerate}
So by \cite[Lemma 4.11]{Kaw21} we obtain that $H^0\left( V, \Omega^1_V(\log G)\otimes \sO(K_V+C)\right)=0$. In \emph{op.cit.} it is assumed that $p>3$: but this is only used there to ensure that $(V,G)$ is relatively snc over a dense open subset of $\bP^1$, which is a given in our situation thanks to \autoref{lemma:qt_invar_conic_bdl}. The proof is therefore complete.
\end{proof}

Fix a lift $(\mathbf{V},\mathbf{C}+\mathbf{G})$ of $(V,C+G)$ as in \autoref{prop:lifting_log_Hirzebruch_canonically}. We let $\mathbf{q}\colon \mathbf{U}\to \mathbf{V}$ be the double cover branched along $\mathbf{G}$. It follows from \autoref{lemma:restriction_cyclic_cover}, \autoref{prop:quotient_is_cyclic_cover} and \autoref{rmk:can_cyclic_cover} that:
$\mathbf{q}$ is a lift of $q$, and 
$\mathbf{E}_\mathbf{U}=\mathbf{q}^{-1}(\mathbf{C})$ is the canonical $\mu_{2,W(k)}$-equivariant lift of $E_U\circlearrowleft \mu_{2,k}$.

To satisfy the hypothesis of \autoref{prop:lift_equiv_MMP}, it remains to show that $\mu_2$-invariant effective divisors can be equivariantly lifted from $U$ to $\mathbf{U}$. As at the end of \autoref{section:equiv_lift_P1xP1} or of \autoref{section:Geiser_inv}, it is sufficient to lift effective divisors from $V$ to $\mathbf{V}$. Since these are Hirzebruch surfaces, we apply \autoref{prop:lift_div_on_Hirzebruch} to conclude.

\subsection{Equivariant liftings: case where $\Gamma_U\neq 0$}\label{section:equiv_lift_Gamma_non_zero}
We treat the special case of \autoref{lemma:min_equiv_model_special_case}. Recall that in this case $m=2$, $U=\bP_B(L_1\oplus L_2)$ where $B$ is an ordinary regular curve of genus $1$, and the divisors $E_U$ and $\Gamma_U$ are the sections respectively given by the projections $L_1\oplus L_2 \twoheadrightarrow L_1$ and $L_1\oplus L_2\twoheadrightarrow L_2$.
The morphism $U\to B$ is non-invariant in this case, so let $\sigma_U$ and $\sigma_B$ denote the non-trivial involutions on $U$ and $B$ respectively. Recall that $E_U$ and $\Gamma_U$ are both $\sigma_U$-invariant. With this data, we let $\mathbf{U}\to\mathbf{B}$ be the canonical lifting constructed in \autoref{section:can_lift_autom_proj_bundles}.
Notice that if $G$ is the sum of the canonical places of $\varphi$ (cf.\ \autoref{def:-2_curves}), then as each irreducible component of $G$ is rational, the components of $h_*G$ on $U$ are fibers of $U\to B$. These fibers lift to $W(k)$-fibers of $\mathbf{U}\to \mathbf{B}$ by construction. So the hypothesis of \autoref{prop:lift_equiv_MMP} are verified.

\subsection{Equivariant liftings: case where $ m=2$ and $E_U$ is reducible}\label{section:equiv_lift_E_reducible} 
The structure of $f\colon (U,E_U)\to B$ is given by \autoref{lemma:min_equiv_model_E_reducible}; since $m=2$, each component of $E_U$ is $\mu_2$-invariant. So up to a change of notations, the situation is the same as in \autoref{section:equiv_lift_Gamma_non_zero}. In particular, we can construct a lift satisfying the hypothesis of \autoref{prop:lift_equiv_MMP}.

\subsection{Equivariant liftings: $m=2$, remaining cases}\label{section:equiv_lift_conic_bd_II}
Assuming $m=2$, the only cases that are left are the following ones: the morphism $f\colon U\to B$ is smooth and non-invariant, the divisor $E_U$ is irreducible and $\Gamma_U=0$.

By \autoref{prop:simplify_mu2_sing}, by blowing-up once the isolated fixed points of the $\mu_2$-action on $U$, we obtain an equivariant birational morphism $v\colon U'\to U$ such that $U'$ is regular and the ramification locus of the $\mu_2$-action on $U'$ has pure codimension one. Let $q'\colon U'\to V'=U'/\mu_2$ be the geometric quotient.
As $V'$ is rational by \autoref{prop:special_surface_is_rational}, $B/\mu_2\cong \bP^1_k$ is the geometric quotient of the induced $\mu_2$-action and we have an induced fibration $g'\colon V'\to \bP^1$ whose general fiber is a regular rational curve. It particular, the minimal model $\widetilde{V}$ of $V'$ over $\bP^1$ is an Hirzebruch surface. The following commutative diagram summarizes the situation:
        \begin{equation}\label{eqn:blow_up_simplify_action}
        \begin{tikzcd}
        U' \arrow[r, "q'"] \arrow[d, "v"] & V' \arrow[dd, "g'" left] \arrow[r, "\phi"]
        & \widetilde{V} \arrow[ddl, "\widetilde{g}"] \\
        U \arrow[d, "f"] & \\
        B\arrow[r] & \bP^1.
        \end{tikzcd}
        \end{equation}
We denote by $\Fix^i(U)$ the union of the $i$-dimensional irreducible components of the fixed locus of the $\mu_2$-action on $U$. We let $E_{U'}=v_*^{-1}E_U$ and $G'$ be the branch locus of $q'$ on $V'$ (it is a divisor by \autoref{prop:sing_qt_pair}). We write $C'=q'(E')$ and $\widetilde{C}=\phi_*C'$. Notice that $\phi$ restricts to an isomorphism $C'\cong \widetilde{C}$ over $\bP^1$.

\begin{lemma}\label{lemma:properties_simplification_action}
In the situation above:
    \begin{enumerate}
        \item $G'$ is $g'$-vertical;
        \item $C'\cong \bP^1$ and $g'|_{C'}\colon C'\to \bP^1$ has degree $2$;
        \item there exists a divisor $\widetilde{F}$ on $\widetilde{V}$ such that:
            \begin{itemize}
                \item $\widetilde{F}$ is the sum of two fibers of $\widetilde{g}$;
                \item $K_{\widetilde{V}}+\widetilde{C}+\frac{1}{2}\widetilde{F}\sim_\bQ 0$, and
                \item $\widetilde{F}\cap \widetilde{C}$ are the four branching points $q_1,\dots,q_4$ of the quotient $E_{U'}\to C'\cong \widetilde{C}$.
            \end{itemize}
    \end{enumerate}
\end{lemma}
\begin{proof}
Since $f\circ v\colon U'\to B$ is equivariant and the $\mu_2$-action on $B$ is non-trivial, it follows that the ramification divisor of $q'$ is vertical over $B$. Therefore $G'$ is vertical over $\bP^1$. 

By \autoref{prop:restriction_linearly_red_qt} and \autoref{prop:log_Gor_cover} we have $C'\cong E_{U'}/\mu_2\cong \bP^1$. The degree of $C'$ over $B/\mu_2$ is equal to the degree of $E_{U'}$ over $B$, since the $\mu_2$-action on $B$ is non-trivial. Therefore $g'|_{C'}\colon C'\to \bP^1$ has degree $2$.

By definition $v$ is the blow-up of $\Fix^0(U)$. As $K_U+E_U\sim 0$, it follows that
        $$K_{U'}+E_{U'}-\sum_{F'\to \Fix^0(U)\setminus E_U}F' \sim 0.$$
By the Riemann--Hurwitz formula we have
        $$K_{U'}= (q')^*( 
        K_{V'})
        +\sum_{F'\to \Fix^0(U)}F' 
        +\sum_{H'\to \Fix^1(U)}H'.$$
As $q^*C'=E_{U'}$, the combination of these two relations yields
        $$0\sim_\bQ K_{V'}+C'+
        \frac{1}{2}\left[ \sum_{F'\to \Fix^0(U)\cap E_U}q'(F')
        +\sum_{H'\to \Fix^1(U)}q'(H')\right].$$
Here $q'(F')+q'(H')\leq G'$ is $g'$-vertical. So pushing forward onto $\widetilde{V}$, we find that 
        $$0\sim_\bQ K_{\widetilde{V}}+\widetilde{C}+\frac{1}{2}\widetilde{F}$$
where $\widetilde{F}$ is a reduced union of fibers of $\widetilde{g}$. By adjunction along $\widetilde{C}$, it follows that $\widetilde{F}$ is the sum of two fibers.

To prove the last statement, observe that $C'\cap G'=\{q_1,\dots,q_4\}$ and that $\widetilde{F}\leq \phi_*G'$. Since $\phi|_{C'}\colon C'\to \widetilde{C}$ is an isomorphism and since $\widetilde{F}\cdot \widetilde{C}=4$, we obtain that $\widetilde{C}\cap \widetilde{F}=\{q_1,\dots,q_4\}$. 
\end{proof}

We use the notations introduced in \autoref{section:equiv_lift_conic_bd_I} for Hirzebruch surfaces: in particular we let $\mathfrak{c}$ be the numerical class of the most negative section of $\widetilde{g}$, and $\mathfrak{f}$ be the class of a fiber. The following easy lemma describes the ample cone.

\begin{lemma}\label{lemma:big_nef_on_Hirzebruch}
Let $H$ be an effective $\bQ$-divisor on $\widetilde{V}$. Then $H$ is big and nef if and only if $H\cdot \mathfrak{f}> 0$ and $H\cdot \mathfrak{c}\geq 0$. If moreover $H\cdot \mathfrak{c}>0$, then $H$ is ample. 
\end{lemma}
\begin{proof}
It is clear that $H$ is nef. To prove that $H$ is big, it suffices to show that $H^2>0$. Write $H=a\mathfrak{c}+b\mathfrak{f}$. 
We compute
        $$0< H\cdot \mathfrak{f}=a, \quad
        0\leq H\cdot \mathfrak{c}=-na+b.$$
Hence $b\geq na>0$ and so both $a$ and $b$ are positive. Moreover:
        $$H^2=-na^2+2ab=a(-na+2b)\geq ab>0$$
as desired. If moreover $H\cdot \mathfrak{c}>0$, then ampleness of $H$ follows from the Nakai--Moishezon criterion.
\end{proof}

\begin{lemma}
We have:
    \begin{enumerate}
        \item $\widetilde{C}^2=4$, and
        \item $\widetilde{V}$ is a del Pezzo surface.
    \end{enumerate}
\end{lemma}
\begin{proof}
Since $\widetilde{V}$ is an Hirzebruch surface, we have $K_{\widetilde{V}}^2=8$. Therefore, as $\widetilde{C}\cdot \widetilde{F}=4$, we obtain:
    $$8=(-K_{\widetilde{V}})^2=
    \left( \widetilde{C}+\frac{1}{2}\widetilde{F}\right)^2
    = \widetilde{C}^2+4$$
and so $\widetilde{C}^2=4$. Since $\widetilde{C}$ has degree $2$ over the base $\bP^1$, its class in the N\'{e}ron--Severi group is not a scaling of the class $\mathfrak{c}$ of the negative section. So $\widetilde{C}+\frac{1}{2}\widetilde{F}$ intersects $\mathfrak{c}$ and a fiber of $\widetilde{g}$ positively. It follows from \autoref{lemma:big_nef_on_Hirzebruch} that $-K_{\widetilde{V}}=\widetilde{C}+\frac{1}{2}\widetilde{F}$ is ample, so $\widetilde{V}$ is a del Pezzo surface.
\end{proof} 

There are only two Hirzebruch surfaces which are also del Pezzo surfaces: either $g$ is the trivial $\bP^1$-bundle, or $\widetilde{V}$ is the blow-up of $\bP^2$ at one point $q$ and $\widetilde{g}$ is the resolution of the map given by the pencil of elements of $|\sO_{\bP^2}(1)|$ through $q$. We distinguish these two cases.

\subsubsection{Case of $\Bl_q\bP^2$}\label{section:Bl_qP2}
Let $\varphi\colon \widetilde{V}\cong \Bl_q\bP^2\to \bP^2$ be the blow-down of the unique $(-1)$-curve to a point $q\in \bP^2$. The pencil of lines through $q$ is resolved by the blow-up $\varphi$, and gives the $\bP^1$-bundle structure $\widetilde{g}$ of $\widetilde{V}$. We let $\overline{C}=\varphi_*\widetilde{C}$ and $\overline{F}=\varphi_*\widetilde{F}$. 

\begin{lemma}
$\overline{F}$ is the sum of two lines passing through $q$, and $\overline{C}\cong \widetilde{C}$ is a regular curve of degree $2$ that does not contain $q$.
\end{lemma}
\begin{proof}
The statement about $\widetilde{F}$ is clear. Since $K_{\bP^2}+\overline{C}+\frac{1}{2}\overline{F}\sim_\bQ 0$, we deduce that $\overline{C}\in |\sO_{\bP^2}(2)|$. As $\overline{C}$ is an irreducible of degree $2$, it must be regular. So the restriction $\varphi|_{\widetilde{C}}\colon \widetilde{C}\to \overline{C}$ is an isomorphism.

Write $\widetilde{C}\equiv a\mathfrak{c}+b\mathfrak{f}$. As $\widetilde{C}\cdot \mathfrak{f}=2$ we have $a=2$. As
        $$4=\widetilde{C}^2=-a^2+2ab=-4+4b$$
we have $b=2$. Therefore $\widetilde{C}\cdot\mathfrak{c}=-2+2=0$, which shows that $C$ is disjoint from the exceptional divisor of $\varphi$. Hence $\widetilde{C}$ does not pass through $q$. 
\end{proof}

By abuse of notation, the images of the points $q_1,\dots,q_4\in \widetilde{C}$ in $\overline{C}$ will also be denoted by $q_1,\dots,q_4$.

\begin{proposition}
There is a lift $(\bP^2_{W(k)},\overline{\mathbf{F}}+\overline{\mathbf{C}})$ of $(\bP^2,\overline{F}+\overline{C})$ such that $\overline{\mathbf{F}}$ intersects $\overline{\mathbf{C}}$ at the canonical lifts $\mathbf{q}_1,\dots,\mathbf{q}_4$ of $q_1,\dots,q_4$ (in the sense of \autoref{def:canonical_lift_P1_4pts}).
\end{proposition}
\begin{proof}
It is clear that the couple $(\bP^2,\overline{F}+\overline{C})$ can be lifted over $W(k)$, so as in the proof of \autoref{prop:lifting_log_Hirzebruch_canonically} it is sufficient to prove that the restriction map
        $$H^1\left(\bP^2, T_{\bP^2}(-\log \overline{F})(-\overline{C})\right)\longrightarrow 
        H^1\left( \overline{C}, T_{\overline{C}}(-\log (\overline{C}\cap \overline{F}))\right)$$
induced by the exact sequence \autoref{eqn:restr_for_twisted_log_vector_fields}, is surjective. In turn, it is sufficient to prove that the cohomology groups 
    $$H^2\left(\bP^2, T_{\bP^2}(-\log \overline{F})(-\overline{C})\right)
    \cong 
    H^0\left( \bP^2, \Omega_{\bP^2}^1(\log \overline{F})\otimes \sO(K_{\bP^2}+\overline{C})\right)^\vee$$
vanish. Since $(\bP^2,\overline{F})$ is snc and $K_{\bP^2}+\overline{C}$ is anti-ample, this follows from Akizuki--Nakano vanishing theorem \cite[Theorem 3.2.(b)]{BBKW24}. The proof is complete.
\end{proof} 

Let $(\bP^2_{W(k)}, \overline{\mathbf{F}}+\overline{\mathbf{C}})$ be as in the above proposition. Starting from it, we are going to lift the whole diagram \autoref{eqn:blow_up_simplify_action} over $W(k)$.
    \begin{enumerate}
        \item We blow-up the $0$-stratum of $\overline{\mathbf{F}}$ (which is a $W(k)$-point), and take strict transforms of divisors, to obtain a lift 
            $$\widetilde{\mathbf{g}}\colon (\widetilde{\mathbf{V}}, \widetilde{\mathbf{C}}+\widetilde{\mathbf{F}})
            \longrightarrow\bP^1_{W(k)}$$ 
        of $\widetilde{g}\colon (\widetilde{V},\widetilde{C}+\widetilde{F})\to \bP^1$. Through the isomorphism $\widetilde{\mathbf{C}}\cong \overline{\mathbf{C}}$, we get $W(k)$-points $\mathbf{q}_1,\dots,\mathbf{q}_4$ on $\widetilde{\mathbf{C}}$.
        \item The morphism $\phi$ is a composition of blow-ups. Each center can be lifted to a $W(k)$-point above $\widetilde{\mathbf{V}}$, which we blow-up. If the center is one some $q_i$, we take $\mathbf{q}_i$ as lift. We obtain a lift
                $$\mathbf{g}'\colon (\mathbf{V}',\mathbf{C}'+\mathbf{G}')\longrightarrow \bP^1_{W(k)}$$
        of $g'\colon (V',C'+G')\to \bP^1$, with the property that $\mathbf{G}'\cap \mathbf{C}'=\{\mathbf{q}_1,\dots,\mathbf{q}_4\}$. 
        \item We let $\mathbf{q}'\colon \mathbf{U}'\to \mathbf{V}'$ be the double cover branched over $\mathbf{G}'$, and let $\mathbf{E}_{\mathbf{U}'}=(\mathbf{q}')^{-1}(\mathbf{C}')$. 
        Then by \autoref{lemma:restriction_cyclic_cover} and \autoref{rmk:can_cyclic_cover}: $\mathbf{q}'$ lifts $q'$, the $\mu_{2,W(k)}$-pair $(\mathbf{U}',\mathbf{E}_{\mathbf{U}'})$ is an equivariant lift of $(U',E_{U'})$,
        and the $\mu_{2,W(k)}$-action on $\mathbf{E}_{\mathbf{U}'}$ it is the canonical lift of the $\mu_{2,k}$-action on $E_{U'}$. 
        The Stein factorization of $\mathbf{q}'\circ\mathbf{g}'$ gives a morphism $\mathbf{U}'\to \mathbf{B}$, where $\mathbf{B}$ is a lift of $B$.
        \item The irreducible components of $\mathbf{G}'$ in $\mathbf{U}'$ are $\mu_{2,W(k)}$-invariant relative $(-1)$-curves that are vertical over $\mathbf{B}$. So we can blow them down to obtain a flat morphism $\mathbf{f}\colon \mathbf{U}\to \mathbf{B}$ whose fibers are regular. The morphism $\mathbf{U}'\to \mathbf{U}$ is naturally $\mu_{2,W(k)}$-equivariant, and restricts to an equivariant isomorphism $\mathbf{E}_{\mathbf{U}'}\cong \mathbf{E}_\mathbf{U}\subset \mathbf{U}$.
    \end{enumerate}
We have therefore obtained a $\mu_{2,W(k)}$-equivariant lift $(\mathbf{U},\mathbf{E}_\mathbf{U})$ restricting to the canonical equivariant lift of $E_U$.

\subsubsection{Case of $\bP^1\times\bP^1$}
Suppose that $\widetilde{V}\cong \bP^1\times\bP^1$ and that $\widetilde{g}$ is the first projection. 
Then $\widetilde{F}$ has bidegree $(2,0)$, and it follows easily from $K_{\widetilde{V}}+\widetilde{C}+\frac{1}{2}\widetilde{F}\sim_\bQ 0$ that $\widetilde{C}$ has bidegree $(1,2)$. So $\widetilde{g}|_{\widetilde{C}}\colon \widetilde{C}\to \bP^1$ has degree $2$ (as we already knew) and $\pr_2|_{\widetilde{C}}\colon \widetilde{C}\to \bP^1$ is an isomorphism. Up to relabelling, we may assume that
        $$\widetilde{g}|_{\widetilde{C}}^{-1}(0)=\{q_1,q_2\}, \quad 
        \widetilde{g}|_{\widetilde{C}}^{-1}(\infty)=\{q_3,q_4\}.$$
Let $(\widetilde{\mathbf{C}},\mathbf{q}_1+\dots+\mathbf{q}_4)$ be the canonical lift of $(\widetilde{C}, q_1+\dots+q_4)$ (cf.\ \autoref{def:canonical_lift_P1_4pts}).

\begin{lemma}\label{lemma:lift_deg_2_map_to_P1}
There exists a lift $\widetilde{\mathbf{C}}\to \bP^1_{W(k)}$ of $\widetilde{g}|_{\widetilde{C}}$ such that the preimage of $0_{W(k)}$ is $\{\mathbf{q}_1,\mathbf{q}_2\}$, and the preimage of $\infty_{W(k)}$ is $\{\mathbf{q}_3,\mathbf{q}_4\}$.
\end{lemma}
\begin{proof}
The morphism $\widetilde{g}|_{\widetilde{C}}$ is equivalent to a non-trivial involution $\iota$ of the pair $(\widetilde{C},\sum_{i=1}^4q_i)$ such that 
$\iota(q_1)=q_2$ and $\iota(q_3)=q_4$. 
By \autoref{prop:functoriality_can_lift_P1_4pts} there is an automorphism $\iota_{W(k)}$ of $\widetilde{\mathbf{C}}$ that lifts $\iota$ and such that
        $$\iota_{W(k)}(\mathbf{q}_1)=\mathbf{q}_2, \quad
        \iota_{W(k)}(\mathbf{q}_2)=\mathbf{q}_1, \quad 
        \iota_{W(k)}(\mathbf{q}_3)=\mathbf{q}_4, \quad
        \iota_{W(k)}(\mathbf{q}_4)=\mathbf{q}_3.$$
As $\iota_{W(k)}^2$ fixes each $\mathbf{q}_i$, by \autoref{rmk:PGL2_not_4transitive} it must be the identity. The quotient $\widetilde{\mathbf{C}}\to \widetilde{\mathbf{C}}/\langle \iota_{W(k)}\rangle$ is the desired lift of $\widetilde{g}|_{\widetilde{C}}$.
\end{proof}

It follows from this lemma that there is a lift $(\widetilde{\mathbf{V}},\widetilde{\mathbf{C}}+\widetilde{\mathbf{F}})$ of $(\widetilde{V}, \widetilde{C}+\widetilde{F})$ over $W(k)$, where $\widetilde{\mathbf{C}}$ is as above and $\widetilde{\mathbf{C}}\cap \widetilde{\mathbf{F}}=\{\mathbf{q}_1,\dots,\mathbf{q}_4\}$.

Starting from this pair, we lift the entire diagram \autoref{eqn:blow_up_simplify_action} over $W(k)$. The method is the same as in \autoref{section:Bl_qP2}, and so we omit the details. We obtain a $\mu_{2,W(k)}$-equivariant lift $(\mathbf{U},\mathbf{E}_\mathbf{U})$ of $(U,E_U)$ such that the induced $\mu_{2,W(k)}$-action on $\mathbf{E}_\mathbf{U}$ is the canonical lift of the $\mu_{2,k}$-action on $E_U$.

\subsubsection{Conclusion in both cases}
Irrespectively of the isomorphism class of $\widetilde{V}$, we have constructed a suitable lift of $(U,E_U)$ over $B$. To satisfy the hypothesis of \autoref{prop:lift_equiv_MMP}, it remains to lift the divisor $h_*G$ equivariantly. If $B$ has genus $1$, then the components of $h_*G$ are fibers of $U\to B$ and this is clear. Otherwise, observe that $v^{-1}_*(h_*G)$ is a $\mu_{2}$-invariant divisor on $U'$, and so it is sufficient to lift its image $G_{V'}$ in $V'$. In turn it is sufficient to lift $\phi_*G_{V'}$ along $\widetilde{V}\hookrightarrow \widetilde{\mathbf{V}}$, which follows from \autoref{prop:lift_div_on_Hirzebruch}.

\subsection{The case $m=4$}\label{section:equiv_lift_m=4}
We conclude with a further analysis of this case, building on \autoref{lemma:min_equiv_model_elliptic_case}. We follow the method used in \autoref{section:equiv_lift_conic_bd_II}, and we will eventually see that $m=4$ cannot occur. 

Recall that $\Gamma_U=0$, that $E_U$ has two disjoint irreducible components which are sections of $f\colon U\to B$, and that the $\mu_{4}$-action on the ordinary genus $1$ curve $B$ is faithful.  We require the following lemma:

\begin{lemma}\label{lemma:order_4_on_P1}
Suppose that $\mu_4$ acts faithfully on $\bP^1$. Then there are exactly $2$ fixed points, and no strictly $2$-periodic points.
\end{lemma}
\begin{proof}
Let $\phi\in \Aut_k(\bP^1)$ be a generator of the action. It has a fixed point (intersect the graph of $\phi$ in $\bP^1\times\bP^1$ with the diagonal), so let us choose a coordinate $x$ such that $[1:0]$ is a fixed point. Then $\phi(x)=ax+b$ for some $a,b\in k$. The condition $\phi^{4}=\id$ is then equivalent to 
        $$a^4=1, \quad b(a^3+a^2+a+1)=0.$$
Together with the condition that $\phi^2\neq \id$, this implies that $a=\zeta_4$ is a primitive fourth root of unity and $b\in k$. Solving $\zeta_4 x_0+b=x_0$, we find that the only other fixed point besides $[1:0]$ is $[b/(1-\zeta_4):1]$.

Suppose that $[\alpha:1]$ is $2$-periodic, so $\zeta_4^2\alpha+\zeta_4b+b=\alpha$. Then $\alpha=b(1+\zeta_4)/2$. But
        $$\frac{b(1+\zeta_4)}{2}=\frac{b}{1-\zeta_4}$$
so in fact $[\alpha:1]$ is a fixed point.
\end{proof} 

Now let $\sigma_U\in \Aut_k(U)$ be a generator of the $\mu_4$-action on $U$, and let $\sigma_B\in \Aut_k(B)$ be such that $f\circ \sigma_U=\sigma_B\circ f$ holds. We write $E_U=E_{U,1}\sqcup E_{U,2}$.

\begin{lemma}\label{lemma:m=4_case}
In the above situation:
    \begin{enumerate}
        \item $\sigma_B$ has exactly $2$ fixed points $\beta_1,\beta_2$ and $2$ strictly $2$-periodic points $\gamma_1,\gamma_2$;
        \item any singular fiber of $f$ must dominate a $\sigma_B$-fixed point;
        \item For each $\beta_i$, either:
            \begin{itemize}
                \item $U_{\beta_i}$ is smooth, $\sigma^2_U|_{U_{\beta_i}}$ is trivial and the $\sigma_U$-fixed points on $U_{\beta_i}$ are disjoint from $E_U\cap U_{\beta_i}$; or,
                \item $U_{\beta_i}$ is singular, $\sigma^2_U|_{U_{\beta_i}}$ is not trivial on any component and the unique $\sigma_U$-fixed point on $U_{\beta_i}$ is its singular point.
            \end{itemize}
        \item For each $\gamma_i$, either:
            \begin{itemize}
                \item $\sigma_U^2|_{U_{\gamma_i}}$ is trivial; or,
                \item $\sigma_U^2|_{U_{\gamma_i}}$ is not trivial on $U_{\gamma_i}$ and $E_U\cap U_{\gamma_i}$ is the set of $\sigma_U^2|_{U_{\gamma_i}}$-fixed points.
            \end{itemize}
    \end{enumerate}
\end{lemma}
\begin{proof}
Let $a_1,\dots,a_4\in E_{U,1}$ and $b_1,\dots,b_4\in E_{U,2}$ be the fixed points of $\sigma^2_U|_{E_{U,1}}$ and $\sigma^2_U|_{E_{U,2}}$ respectively. We let $\mathfrak{R}=\{a_1,\dots,a_4,b_1,\dots,b_4\}$. Up to relabeling we have 
        $$\sigma_U\colon a_i\mapsto b_i\mapsto a_i.$$
By equivariance of $f$, this implies that 
        $$\sigma_B\colon f(a_i)\mapsto f(b_i)\mapsto f(a_i).$$
Since $E_U\to B$ has degree $2$, for any index $i$ there is at most one index $j\in\{1,\dots,4\}$ such that $f(a_i)=f(b_j)$.

Since the action of $\sigma_B$ admits $2$-periodic points (namely, the image of $\mathfrak{R}$), the quotient morphism $q_B\colon B\to B/\mu_4$ is not \'{e}tale. In particular $B/\mu_4\cong \bP^1$. Let $R$ be the ramification divisor of $q_B$: by Riemann--Hurwitz we see that $\deg R=8$. 

We know that every point of $f(\mathfrak{R})$ is fixed by $\sigma^2_B$, so we get $f(\mathfrak{R})\subseteq \Supp(R)$. If $b\in B$ is fixed by $\sigma_B$, then by equivariance $\sigma_U$ exchanges the two points of $E_U$ above $b$: by definition these two points belong to $\mathfrak{R}$. If $c_1,c_2\in B$ are such that $\sigma_B\colon c_1\mapsto c_2\mapsto c_1$, then again by equivariance we see that the points of $E_U$ over $c_1$ and $c_2$ must belong to $\mathfrak{R}$. This proves that $\Supp(R)= f(\mathfrak{R})$, and that for every $i$ there is a unique $j$ such that $f(a_i)=f(b_j)$.

The divisor $R$ is the sum of the $\sigma_B$-fixed points with coefficient $3$, and of the strictly $2$-periodic points with coefficient $1$. As $f(\mathfrak{R})$ contains four points and $\deg R=8$, we see that $\sigma_B$ has exactly $2$ fixed points $\beta_1,\beta_2$ and $2$ strictly $2$-periodic points $\gamma_1,\gamma_2$. 

Suppose that $U_{\beta_i}$ is a regular fiber. As the points of $E_U\cap U_{\beta_i}$ are strictly $2$-periodic, by \autoref{lemma:order_4_on_P1} the action of $\mu_4$ on $U_{\beta_i}$ cannot be faithful (nor trivial). So $\sigma^2_U|_{U_{\beta_i}}$ must be trivial.

Now let $F_1\cup F_2$ be a singular fiber. Both $F_0$ and $F_1$ are $(-1)$-curves, so as $K_U+E_U\sim 0$ we get that $F_i\cdot E_U=1$ for both $i$. We choose indices so that $F_i\cdot E_{U,i}=1$ for $i=1,2$. Since $f$ is minimal, we have $\sigma_U(F_1)\neq F_1$. Suppose that $\sigma_U(F_1)$ is distinct from $F_2$. Then there is another singular fiber $F_1'\cup F_2'$, where the indices are chosen with the same convention, such that $\sigma_U(F_1\cup F_2)=F_1'\cup F_2'$. As
        $$\sigma_U(F_1)\cdot E_{U,2}
        =F_1\cdot \sigma_U^{-1}(E_{U,2})
        =F_1\cdot E_{U,1}=1,$$
we must have
        $$\sigma_U(F_1)=F_2', \quad \sigma_U(F_2)=F_1'.$$
The same reasoning shows that $\sigma_U^2(F_1)$ must meet $E_{U,1}$, so it cannot be equal to $F_2$. It follows that if $\sigma_U(F_1)\neq F_2$ then the orbit of $F_1$ consists of disjoint irreducible curves: this contradicts the minimality of $f$. Thus $\sigma(F_1)=F_2$, which by equivariance implies that $f(F_1\cup F_2)$ is a fixed point of $\sigma_B$. 

In case the fiber over $\beta_i$ is singular, let us record the following claim:
\begin{claim}\label{lemma:eigenvalues_at_fixed_pt}
Suppose that $U_{\beta_i}$ is singular, with singular point $u$. Then $u$ is fixed by $\sigma_U$, and the eigenvalues of the induced action on $\fm_u/\fm_u^2$ are $\zeta_4$ and $\zeta_4^3$, where $\zeta_4$ is a primitive fourth root of unity.    
\end{claim}
\begin{proof}\renewcommand{\qedsymbol}{$\lozenge$}
Since $u$ is the only singular point of $U_{\beta_i}$, we must have $\sigma_U(u)=u$. The two irreducible components of the fiber give a basis $v_1,v_2$ of $\fm_u/\fm_u^2$. Denote by $\overline{\sigma}_u$ the induced automorphism of $\fm_u/\fm_u^2$. We have shown above that $\overline{\sigma}_u$ exchanges the sub-spaces generated by $v_1$ and $v_2$: so in this basis it is given by a matrix
        $$\begin{pmatrix}
        & \beta \\
        \alpha &
        \end{pmatrix}.$$
Then $\overline{\sigma}_u^2$ is given by the scalar matrix $\alpha\beta\cdot I_2$, which must be of order $2$: so $\alpha\beta=\pm 1$. If $\alpha\beta=1$ then $\overline{\sigma}_u$ and therefore $\sigma_U$ are of order $2$, which is impossible in our situation. So $\alpha\beta=-1$. It follows that the characteristic polynomial of the matrix is $t^2+1$, and so its eigenvalues are $\zeta_4$ and $-\zeta_4=\zeta_4^3$.
\end{proof}

In the above notations, notice that $\overline{\sigma}^2_u$ sends $v_j$ to $-v_j$: this proves that $\sigma^2_U$ does not restrict to the identity on either component of $U_{\beta_i}$.

The fiber $U_{\gamma_i}$ must be regular, as $\gamma_i$ is not fixed by $\sigma_B$. If $\sigma_U^2$ does not restrict to the identity on $U_{\gamma_i}$, then it has exactly two fixed points by \autoref{prop:involutions_on_P1}, which must be the intersection points of $E_U$ with $U_{\gamma_i}$. The proof is complete.
\end{proof}

We refine our notations as follows:
    \begin{itemize}
        \item we denote by $\gamma_i$ (resp.\ $\overline{\gamma_i}$) the strictly $2$-periodic points of $\sigma_B$ such that $\mu_4$ acts faithfully on $U_{\gamma_i}$ (resp.\ factors through a $\mu_2$-action on $U_{\overline{\gamma}_i}$);
        \item we denote by $\beta_i$ (resp.\ $\overline{\beta_i}$) the fixed points of $\sigma_B$ such that $U_{\beta_i}$ is singular (resp.\ $U_{\overline{\beta}_i}$ is regular);
        \item we denote by $u_i$ (respectively $v_i,w_i$) the $\sigma_U$-fixed points of the fiber $U_{\beta_i}$ (resp.\ of $U_{\overline{\beta}_i}$).
    \end{itemize}
The following picture represents the possible fibers of $f$ above the ramification points of $\sigma_B$. The black points are fixed by the $\mu_4$-actions, and the points and components in cyan are strictly $2$-periodic. We have drawn the components of $E_U$ horizontally in teal. The dotted arrows indicate the action of $\sigma_U$ when there is an ambiguity.

    \medskip
    \begin{center} 
    \begin{tikzpicture}
    \draw[-] (-0.5,0) -- (2,0) ;
    \draw[-] (2.5,0) -- (5,0) ;
    \draw[-] (5.5,0) -- (6.5,0) ;
    \draw[-] (7,0) -- (8,0) ;
    \node at (9,0) {$B$} ;
    \node at (9, 2) {$U$} ;
    \draw[->] (9,1.4) -- (9, 0.6) node[midway, right] {$f$} ;

    \node at (-1, 1.5) {$E_{U,1}$} ;
    \draw[-, teal, thick] (-0.5, 1.5) -- (-0.1,1.5) ;
    \draw[-, teal, thick] (0.1, 1.5) -- (1.4,1.5) ;
    \draw[-, teal, thick] (1.6, 1.5) -- (2,1.5) ;
    \draw[-, teal, thick] (2.5, 1.5) -- (2.9, 1.5) ;
    \draw[-, teal, thick] (3.1, 1.5) -- (4.4,1.5) ;
    \draw[-, teal, thick] (4.6, 1.5) -- (5,1.5) ;
    \draw[-, teal, thick] (5.5, 1.5) -- (6.15,1.5) ;
    \draw[-, teal, thick] (6.35, 1.5) -- (6.5,1.5) ;
    \draw[-, teal, thick] (7, 1.5) -- (7.4,1.5) ;
    \draw[-, teal, thick] (7.6, 1.5) -- (8,1.5) ;

    \node at (-1, 2.5) {$E_{U,2}$} ;
    \draw[-, teal, thick] (-0.5, 2.5) -- (-0.1,2.5) ;
    \draw[-, teal, thick] (0.1, 2.5) -- (1.4,2.5) ;
    \draw[-, teal, thick] (1.6, 2.5) -- (2,2.5) ;
    \draw[-, teal, thick] (2.5, 2.5) -- (2.9, 2.5) ;
    \draw[-, teal, thick] (3.1, 2.5) -- (4.4,2.5) ;
    \draw[-, teal, thick] (4.6, 2.5) -- (5,2.5) ;
    \draw[-, teal, thick] (5.5, 2.5) -- (6.15,2.5) ;
    \draw[-, teal, thick] (6.35, 2.5) -- (6.5,2.5) ;
    \draw[-, teal, thick] (7, 2.5) -- (7.4,2.5) ;
    \draw[-, teal, thick] (7.6, 2.5) -- (8,2.5) ;
   
    \filldraw[cyan] (0,0) circle (2pt) node[below]{{\color{black}\footnotesize{$\gamma_1$}}} ; 
    \draw[-] (0,1)--(0,3) ;
    \filldraw[cyan] (0,1.5) circle (2pt) ; 
    \filldraw[cyan] (0,2.5) circle (2pt) ; 

    \filldraw[cyan] (1.5,0) circle (2pt) node[below]{{\color{black}\footnotesize{$\gamma_2$}}} ; 
    \draw[-] (1.5,1)--(1.5,3) ;
    \filldraw[cyan] (1.5,1.5) circle (2pt) ; 
    \filldraw[cyan] (1.5,2.5) circle (2pt) ;

    \draw[<->, dotted] (0.1,2.4) -- (1.4,1.6) ;
    \draw[<->, dotted] (0.1,1.6) -- (1.4,2.4) ;

    \filldraw[cyan] (3,0) circle (2pt) node[below]{{\color{black}\footnotesize{$\overline{\gamma}_1$}}} ; 
    \draw[-, cyan, thick] (3,1)--(3,3) ;
    \filldraw[cyan] (3,1.5) circle (2pt) ; 
    \filldraw[cyan] (3,2.5) circle (2pt) ; 

    \filldraw[cyan] (4.5,0) circle (2pt) node[below]{{\color{black}\footnotesize{$\overline{\gamma}_2$}}} ; 
    \draw[-, cyan, thick] (4.5,1)--(4.5,3) ;
    \filldraw[cyan] (4.5,1.5) circle (2pt) ; 
    \filldraw[cyan] (4.5,2.5) circle (2pt) ;

    \draw[<->, dotted] (3.1,2.4) -- (4.4,1.6) ;
    \draw[<->, dotted] (3.1,1.6) -- (4.4,2.4) ;

    \filldraw[black] (6,0) circle (2pt) node[below]{\footnotesize{$\beta_i$}} ;
    \draw[-] (5.8,1.6) -- (6.5,3) ;
    \draw[-] (5.8,2.4) -- (6.5,1) ;
    \filldraw[black] (6,2) circle (2pt) node[left]{\footnotesize{$u_i$}};
    \filldraw[cyan] (6.25, 2.5) circle (2pt) ;
    \filldraw[cyan] (6.25, 1.5) circle (2pt) ;

    \filldraw[black] (7.5,0) circle (2pt) node[below]{\footnotesize{$\overline{\beta}_i$}} ;
    \draw[-, cyan, thick] (7.5,1) -- (7.5,3) ;
    \filldraw[cyan] (7.5,1.5) circle (2pt) ; 
    \filldraw[cyan] (7.5,2.5) circle (2pt) ;
    \filldraw[black] (7.5, 1.8) circle (2pt) node[left]{\footnotesize{$v_i$}};
    \filldraw[black] (7.5, 2.2) circle (2pt) node[left]{\footnotesize{$w_i$}};
    \end{tikzpicture}
    \end{center}

Now let $v\colon U'\to U$ be the blow-up of all the isolated strictly $2$-periodic points, and $E_{U'}$ be the strict transform of $E_U$. Notice that all the points we blow-up belong to $E_U$. So $K_{U'}+E_{U'}\sim 0$ and the action of $\mu_4$ lifts to $(U',E_{U'})$. The fixed locus $\Fix_{U'}(\mu_2)$ of the action of the sub-group $\mu_2\subset \mu_4$ on $U'$ is the union of the $v$-exceptional divisors (cf.\ \autoref{prop:simplify_mu2_sing}), of the (strict) transforms of the divisorial components of $\Fix_U(\mu_2)$ and of the points $u_i',v_i',w_i'$ that lie respectively above the points $u_i,v_i,w_i$. Notice that the divisorial locus $\Fix^1_{U'}(\mu_2)$ of $\Fix_{U'}(\mu_2)$ is invariant under the $\mu_4$-action.

Let $q_1\colon U'\to V_1$ be the quotient by the sub-group $\mu_2$: it is fibered over $B_1=B/\langle \sigma_B^2\rangle \cong \bP^1$, so we obtain a commutative diagram
    $$\begin{tikzcd}
    U'\arrow[r, "q_1"] \arrow[d, "v"] & V_1 \arrow[dd, "f_1"] \\
    U\arrow[d, "f"] & \\
    B\arrow[r, "\mathfrak{q}_1"] & B_1.
    \end{tikzcd}$$
We use the following notations: $G_1$ is the image of $\Fix_{U'}^1(\mu_2)$, $C_1=q_1(E_{U'})$, $\gamma_i^1=\mathfrak{q}_1(\gamma_i)$ and $u_i^1=q_1(u_i')$, and similarly for $\overline{\gamma}_i^1$, $v_i^1$, etc. Then:
    \begin{itemize}
        \item $K_{V_1}+C_1+\frac{1}{2}G_1\sim_\bQ 0$ by the Riemann--Hurwitz formula,
        \item $C_1=C_{1,1}\sqcup C_{1,2}$ where each $C_{1,i}$ is a section of $f_1$ (by \autoref{prop:restriction_linearly_red_qt}),
        \item by \cite[Corollary 2.43]{kk-singbook} the pair $(V_1,C_1+\frac{1}{2}G_1)$ is lc, and by \autoref{lemma:sing_qt} the only singular points of the surface $V_1$ are the $q(u_i')$, none of which is contained in $C_1+G_1$. 
    \end{itemize}
The $\mu_4$-action on $\left(U',E_{U'}+\Fix^1_{U'}(\mu_2)\right)$ and $B$ descends to a $\mu_4/\mu_2\cong \mu_2$-action on $(V_1,C_1+G_1)$ and $B_1$, and the fibration $f_1\colon V_1\to B_1$ is $\mu_2$-equivariant. The action on $V_1$ exchanges the two components of $C_1$.

Below we draw the possible fibers of $f_1$ over the images of the ramification points of $\sigma_B$, in a way that corresponds to the previous picture. The components of $G_1$ are colored in olive, and we have drawn one component of $C_1$ and its intersection with the fibers. The dotted arrows indicate the effect of the $\mu_2$-action.

    \medskip 
    \begin{center}
    \begin{tikzpicture}
    \draw[-] (-0.5,0) -- (2,0) ;
    \draw[-] (2.5,0) -- (5,0) ;
    \draw[-] (6,0) -- (7,0) ;
    \draw[-] (7.5,0) -- (8.5,0) ;

    \node at (-2,0.7) {$C_{1,1}$} ;
    \draw[-, brown, thick] (-1.5,0.7) -- (-0.7,0.7) ;
    \draw[-, brown, thick] (-0.9,0.7) -- (0.6,0.7) ;
    \draw[-, brown, thick] (0.8,0.7) -- (2,0.7) ;
    \draw[-, brown, thick] (2.5,0.7) -- (2.9,0.7) ;
    \draw[-, brown, thick] (3.1,0.7) -- (4.4,0.7) ;
    \draw[-, brown, thick] (4.6,0.7) -- (5,0.7) ;
    \draw[-, brown, thick] (6,0.7) -- (6.25,0.7) ;
    \draw[-, brown, thick] (6.45,0.7) -- (7,0.7) ;
    \draw[-, brown, thick] (7.5,0.7) -- (7.9,0.7) ;
    \draw[-, brown, thick] (8.1,0.7) -- (8.5,0.7) ;
    
    \node at (9.75,0) {$B_1$} ;
    5\node at (9.75, 2) {$V_1$} ;
    \draw[->] (9.75,1.4) -- (9.75, 0.6) node[midway, right]{$f_1$} ;

    \filldraw[black] (0,0) circle (2pt) node[below]{\footnotesize{$\gamma_{1}^1$}} ; 
    \draw[-] (0,1)--(0,3) ;
    \draw[-, olive, thick] (-1,0.5) -- (0.3,1.8) ;
    \draw[-, olive, thick] (-0.7,3.3) -- (0.3,2.3) ; 
    \filldraw[brown] (-0.8,0.7) circle (2pt) ;

    \filldraw[black] (1.5,0) circle (2pt) node[below]{\footnotesize{$\gamma_{2}^1$}} ; 
    \draw[-] (1.5,1)--(1.5,3) ;
    \draw[-, olive, thick] (0.5,0.5) -- (1.8,1.8) ;
    \draw[-, olive, thick] (0.8,3.3) -- (1.8,2.3) ; 
    \filldraw[brown] (0.7,0.7) circle (2pt) ;

    \draw[<->, dotted] (0.2,-0.1) .. controls (0.6,-0.5) and (0.9, -0.5) .. (1.3,-0.1) ;
    \draw[<->, dotted] (0.25, 1.55) .. controls (0.9,1.6) and (1,2.45) .. (1.1, 2.9) ;

    \filldraw[black] (3,0) circle (2pt) node[below]{\footnotesize{$\overline{\gamma}_1^1$}} ; 
    \draw[-, olive, thick] (3,0.5)--(3,3.3) ;
    \filldraw[brown] (3,0.7) circle (2pt) ;

    \filldraw[black] (4.5,0) circle (2pt) node[below]{\footnotesize{$\overline{\gamma}_2^1$}} ; 
    \draw[-, olive, thick] (4.5,0.5)--(4.5,3.3) ;
    \filldraw[brown] (4.5,0.7) circle (2pt) ;

    \draw[<->, dotted] (3.2,-0.1) .. controls (3.6,-0.5) and (3.9, -0.5) .. (4.3,-0.1) ;
   
    \filldraw[black] (6.5,0) circle (2pt) node[below]{\footnotesize{$\beta_i^1$}} ;
    \draw[-] (6.3,1.6) -- (7,3) ;
    \draw[-] (6.3,2.4) -- (7,1) ;
    \filldraw[violet] (6.5,2) circle (2pt) node[right]{{\color{black}\footnotesize{$u_i^1$}}};
    \draw[-, olive, thick] (6.3, 3.4)-- (6.85, 2.3) ;
    \draw[-, olive, thick] (6.25, 0.5) -- (6.85,1.7) ;
    \filldraw[brown] (6.35,0.7) circle (2pt) ;

    \draw[<->, dotted] (6.25, 2.35) .. controls (6.15, 2.1) and (6.15, 1.9) .. (6.25, 1.65) ;
    \draw[<->, dotted] (6.3, 3) .. controls (5.6, 2.4) and (5.6, 1.6) .. (6.3, 1) ;

    \filldraw[black] (8,0) circle (2pt) node[below]{\footnotesize{$\overline{\beta}_i^1$}} ; 
    \draw[-, olive, thick] (8,0.5)--(8,3.3) ;
    \filldraw[violet] (8, 1.4) circle (2pt) node[right]{{\color{black}\footnotesize{$v_i^1$}}} ;
    \filldraw[violet] (8, 2.4) circle (2pt) node[right]{{\color{black}\footnotesize{$w_i^1$}}} ;
    \filldraw[brown] (8,0.7) circle (2pt) ;
    
    \end{tikzpicture}
    \end{center}

\begin{lemma}
The only fixed points of the $\mu_2$-action on $V_1$ are the $u_i^1,v_i^1$ and $w_i^1$.
\end{lemma}
\begin{proof}
By equivariance of $f_1$, the fixed points on $V_1$ must dominate the fixed points on $B_1$. There are the $\beta_i^1$s and the $\overline{\beta}_i^1$, and the possible fibers above them are pictured above. Clearly the unique fixed point of $f_1^{-1}(\beta_i^1)$ is $u_i^1$. The $\mu_2$-action on $f_1^{-1}\left(\overline{\beta}_i^1\right)$ corresponds to the action of $\sigma_U$ on $U_{\overline{\beta}_i}$, which had $v_i$ and $w_i$ as fixed points.
\end{proof}

Now let $q_2\colon V_1\to V_2$ be the geometric quotient by the $\mu_2$-action. The equivariant morphism $f_1$ descends to a fibration $f_2\colon V_2\to B_2=B_1/\mu_2$. As before, we write $G_2=q_2(G_1)$ and $C_2=q_2(C_1)$, $u_i^2=q_2(u_i^1)$ and similarly for $v_i^2,w_i^2$. The quotient $B_1\to B_2$ identifies the $\gamma_i^1$ (or the $\overline{\gamma}_i^1$, depending on the situation): we let $\gamma^2$ (resp.\ $\overline{\gamma}_2$) stands for their common image. The points $\beta_i^1$ and $\overline{\beta}_i^1$ are not identified, and we denote their images by $\beta_i^2$ and $\overline{\beta}_i^2$. Then:
    \begin{itemize}
        \item By the above lemma, $q_2$ is \'{e}tale in away from $u_i^1,v_i^1$ and $w_i^1$. So \autoref{prop:sing_qt_pair} shows that $G_2+C_2$ is snc and contained in the regular locus of $V_2$, and the Riemann--Hurwitz formula implies $K_{V_2}+C_2+\frac{1}{2}G_2\sim_\bQ 0$.
        \item By \autoref{prop:restriction_linearly_red_qt} we see that $C_2\cong \bP^1$ has degree $2$ over $B_2$.
        \item $G_2$ is contained in exactly three fibers.
        \item The singular points of $V_2$ are the points $u_i^2,v_i^2$ and $w_i^2$, and none of them is not contained in $G_2+C_2$. The singularity at $u_i^2$ is the quotient singularity $(u_i\in U)/\mu_4$. By \autoref{lemma:eigenvalues_at_fixed_pt} and the Reid--Tai criterion \cite[Theorem 3.21]{kk-singbook}, the singularities $u_i^2\in V_2$ are canonical.
    \end{itemize}
Let $w\colon V_2'\to V_2$ be the minimal resolution of the canonical singularities $u_i^2$. If we let $C_2'=w^{-1}_*C_2$ and $G_2'=w^{-1}_*G_2$, then we have
    $$K_{V_2'}+C_2'+\frac{1}{2}G_2'\sim_\bQ 0$$ 
and the pair $(V_2',C_2'+\frac{1}{2}G_2')$ is lc. Let $g\colon V_2'\to \widetilde{V}_2$ be the minimal model over $B_2$: the morphism $g$ contracts some components in the fibers above $\gamma^2$ and $\beta_i^2$ (because these are the only reducible fibers), and at the end the fibers of $\widetilde{V}_2\to B_2$ are all irreducible. Notice that the singular points of $\widetilde{V}_2$ belong to the fibers above $\overline{\beta}^2_i$.

Let $\widetilde{C}_2=g_*C_2'$ and $\widetilde{G}_2=g_*G_2'$. The log CY condition on $V_2'$ implies that
        \begin{equation}\label{eqn:log_CY_m=4}
        K_{\widetilde{V}_2}+\widetilde{C}_2+\frac{1}{2}\widetilde{G}_2\sim_\bQ 0
        \end{equation}
and that $w\colon (V_2',C_2'+\frac{1}{2}G_2')\to (\widetilde{V}_2, \widetilde{C}_2+\frac{1}{2}\widetilde{G}_2)$ is a crepant morphism. Therefore $(\widetilde{V}_2, \widetilde{C}_2+\frac{1}{2}\widetilde{G}_2)$ is log CY. Note that $\widetilde{C}_2\neq 0$ is an irreducible rational curve that avoids the singularities of $\widetilde{V}_2$, and that the morphism $\widetilde{C}_2\to B_2$ is dominant of degree $2$.

\begin{lemma}\label{lemma:exactly_three_fibers}
$\widetilde{G}_2$ is the sum of three fibers of $\widetilde{V}_2\to B_2$.
\end{lemma}
\begin{proof}
By construction, the components of $G_2'$ are contained in exactly three different fibers of $V_2'\to B_2$: say $G_2'=F_1'+F_2'+F_3'$ where each $F_i'$ is a (possibly reducible) divisor contained in a single fiber. It is sufficient to prove that $w_*F_i'\neq 0$ for each $i$. By contradiction, assume, say, that $w_*F_1'=0$. Fix a component $F$ of $F_1'$: since $w$ is crepant, we have
        $$-\frac{1}{2}=a\left(F;\widetilde{V}_2, \widetilde{C}_2+\frac{1}{2}w_*(F_2'+F_3')\right).$$
Since the center of $F$ on $\widetilde{V}_2$ is disjoint from $w_*(F_2'+F_3')$, we actually have
        $$-\frac{1}{2}=a\left(F;\widetilde{V}_2, \widetilde{C}_2\right).$$
On the other hand, observe that $F$ does not belong to a fiber above $\overline{\beta}_i$, as $g$ is an isomorphism around such a fiber. Thus the center of $F$ on $\widetilde{V}_2$ belongs to the regular locus of that surface. In particular, $K_{\widetilde{V}_2}+\widetilde{C}_2$ is Cartier at the center of $F$: but then $a(F;\widetilde{V}_2, \widetilde{C}_2)$ is an integer and cannot be equal to $-1/2$. We have reached a contradiction, and so the lemma is proved.
\end{proof} 

We can now move towards the conclusion: let us perform adjunction for the pair $(\widetilde{V}_2, \widetilde{C}_2+\frac{1}{2}\widetilde{G}_2)$ along $\widetilde{C}_2^\nu\cong \bP^1$. By \autoref{eqn:log_CY_m=4} we obtain
        $$2=\deg \Diff_{\widetilde{C}_2^\nu}\left( \frac{1}{2}\widetilde{G}_2\right)
        = \deg N + \frac{1}{2}\deg \left(\widetilde{G}_2|_{\widetilde{C}_2}\right)$$
where $N\geq 0$ is the preimage of the singularities of $\widetilde{C}_2$ (which can only be nodal by \cite[Theorem 2.31]{kk-singbook}). By \autoref{lemma:exactly_three_fibers} and the fact that $\widetilde{C}_2$ has degree $2$ over $B_2$, we get that $\frac{1}{2}\deg \left(\widetilde{G}_2|_{\widetilde{C}_2}\right)=3$: but then the equalities displayed above become $2=\deg N+ 3$, which cannot hold.

\medskip
Therefore the assumption that $m=4$ leads to a contradiction, and we obtain:

\begin{corollary}\label{cor:m=4_does_not_occur}
Let $(S,D)$ be as in \autoref{notation:4A1_case}: then the Cartier index of $K_S+D$ is $2$.   
\end{corollary}